\documentclass[a4paper]{article}

\usepackage[english]{babel} 
\usepackage[utf8x]{inputenc} 
\usepackage{amsmath}
\usepackage{amssymb}
\usepackage{amsfonts}
\usepackage{amsthm}
\usepackage{multirow}
\usepackage{multicol}
\usepackage{graphicx}
\usepackage{color}
\usepackage[colorinlistoftodos]{todonotes}
\usepackage{enumerate}
\usepackage{tikz,xcolor}
\usepackage{float}
\usepackage{wrapfig}
\usepackage{subcaption}
\usepackage{mathtools}
\usepackage{pgfplots}
\pgfplotsset{compat=1.7}

4

\def\C{\mathbb{C}}
\def\R{\mathbb{R}}

\newtheorem{prop}{Proposition}
\newtheorem{defi}{Definition}
\newtheorem*{remark}{Remark}

\graphicspath{ {./images/} }

\title{Simultaneous reconstruction of conductivity, boundary shape and contact impedances in electrical impedance tomography}

\author{J. P. Agnelli\thanks{FaMAF-CIEM (CONICET), Universidad Nacional de C\'ordoba, Argentina (\texttt{jpagnelli@unc.edu.ar}).  }
\and V. Kolehmainen\thanks{Department of Applied Physics, University of Eastern Finland, Finland (\texttt{ville.kolehmainen@uef.fi}). }
\and M. Lassas \thanks{Department of Mathematics and Statistics, University of Helsinki, Finland (\texttt{matti.lassas@helsinki.fi}, 
\tt{petri.ola@helsinki.fi}, \tt{samuli.siltanen@helsinki.fi}).}
\and P. Ola\footnotemark[3]
\and S. Siltanen\footnotemark[3]}

\begin{document}
	
\maketitle

\begin{abstract}
The objective of electrical impedance tomography (EIT) is to reconstruct the internal conductivity of a physical body based on current and voltage measurements at the boundary of the body. In many medical applications the exact shape of the domain boundary and contact impedances are not available. This is problematic as even small errors in the boundary shape of the computation domain or in the contact impedance values can produce large artifacts in the reconstructed images which results in a loss of relevant information. A method is proposed that simultaneously reconstructs the conductivity, the contact impedances and the boundary shape from EIT data. 
The  approach consists of three steps: 
first, the unknown contact impedances and an anisotropic conductivity reproducing the measured EIT data in a model domain are computed. Second,  using isothermal coordinates, a deformation is constructed that makes the conductivity isotropic. The final step minimizes the error of true and reconstructed known geometric properties (like the electrode lengths) using conformal deformations.
The feasibility of the method is illustrated with experimental EIT data, with robust and accurate reconstructions of both conductivity and boundary shape. 
\end{abstract}



\section{Introduction}

In electrical impedance tomography (EIT) a set of electrodes is attached on the boundary of a physical body. Electric currents are injected to the body, and the resulting voltages on the electrodes are measured. The goal is to reconstruct the conductivity inside the  body from the boundary data. 
Medical applications of EIT include monitoring  heart and lung function~\cite{IMNS04,IMNS06}, detection of breast cancer~\cite{BKKK08}, imaging of brain function and stroke detection~\cite{ACLMSS20,APKSMH16,CHH19}.

Most EIT image reconstruction methods rely on an accurate model of the boundary of the body, which is not always available in medical applications. For example, consider using EIT for monitoring the lung function of an unconscious patient in intensive care~\cite{IMNS04,IMNS06}. 
A model for the boundary of the patient's chest could in principle be obtained using CT or MRI. However,  transporting a critically ill patient to a scanner may not be practical.
Further, chest shape changes with breathing and depends on patient positioning. This is problematic since even small errors in the domain model can cause large artifacts in the reconstructed conductivity images~\cite{AGB96,GHO96,NKKa}.

Another common feature in traditional EIT is to treat the electrode contact impedances as known parameters. They model the voltage drops caused by electrochemical effects  at the electrode-skin interface, and depend on local skin conditions and on temporal variations arising from perspiration and partial drying of the electrode gel. So the contact impedances are known only approximately in practice, which may lead to severe errors in the reconstructed images \cite{HVWKV02,KVKK97}.  
So-called four point measurements, where separate electrodes are used for current injection and voltage measurements, offer some help~\cite{HWWT93}. However, a better option is to consider the contact impedances as unknown parameters in the EIT problem ~\cite{HMS17}. 

One way to deal with imprecise boundary models is to use \emph{difference imaging}, aiming to reconstruct only the change in the conductivity between successive measurements \cite{BB84,SGA06,LKSS15}. Errors caused by (hopefully invariant) model inaccuracies cancel out to some extent when subtracting the two measurements. However, in this paper we focus on {\it absolute} imaging with the goal of reconstructing actual conductivity values. 

Let us review previous attempts to compensate for the domain modeling errors in absolute EIT. In~\cite{KLO05} it was shown that in an inaccurately modeled domain there is a unique minimally anisotropic conductivity matching the EIT data. The square root of the determinant of that conductivity gives a deformed image of the original conductivity defined in the true domain. This methodology was extended to include the estimation of unknown contact impedances in~\cite{KLO08}. On the other hand, in \cite{KLOS13} the original method introduced in~\cite{KLO05} was extended by using isothermal coordinates to transform the reconstructed anisotropic conductivity to an isotropic conductivity close to the original one. The idea of using isothermal coordinates for reducing anisotropic EIT to an isotropic model was introduced by Sylvester in \cite{sylvester1990anisotropic}.

The Bayesian approximation error method was applied to EIT with an imperfectly known boundary in~\cite{NKKa,NKKb}. 
A different approach was considered in~\cite{HKMS17} where the dependence of the electrode measurements on model properties were parametrized via polynomial collocation. 
In~\cite{DHSS13A, DHSS13b}, the authors  designed a reconstruction algorithm using the complete electrode model~\cite{SChI92} and  the Fréchet derivative of the current-to-voltage map with respect to the electrode locations and boundary shape.
A similar approach was considered in~\cite{CHH19} but for the problem of head imaging by EIT. Differently from~\cite{DHSS13A, DHSS13b}, to cope with the instability the so-called smoothened complete electrode model~\cite{HM17} was considered as the forward model.

\begin{figure}[t]
	\centering
	\includegraphics[width=7.5cm,height=5.8cm]{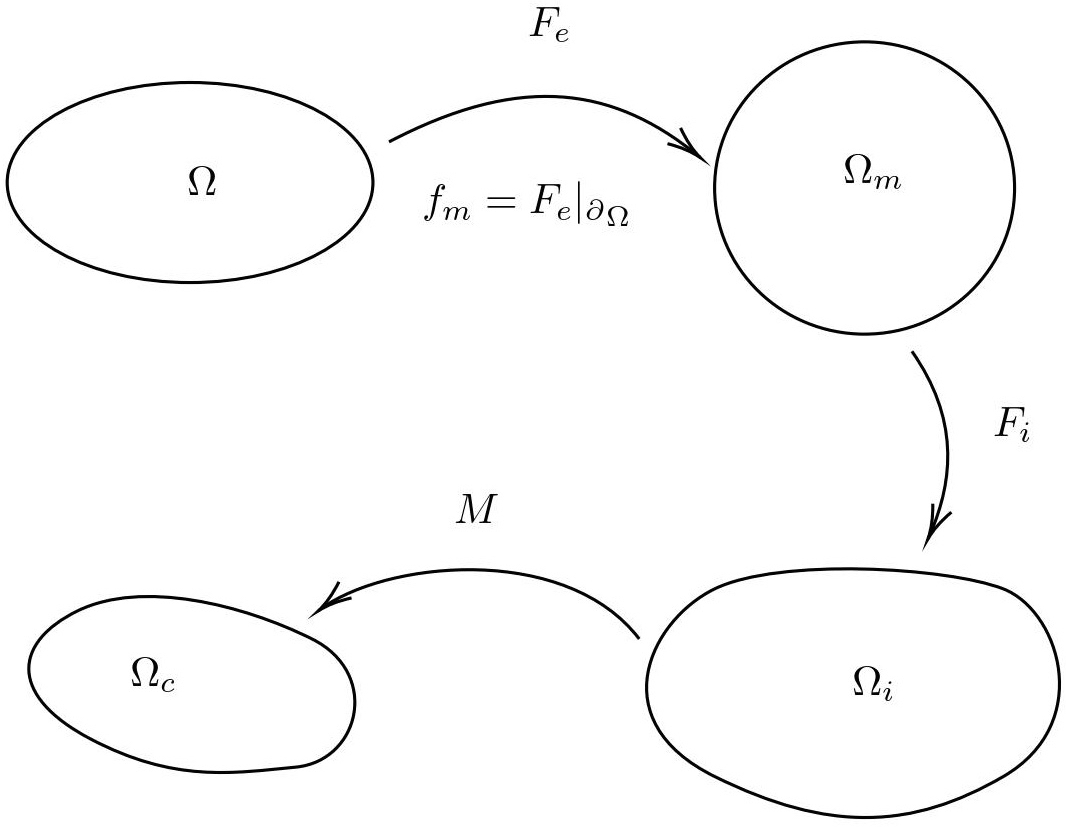}
	\caption{The different steps of the method and the corresponding maps.}
	\label{fig:mappings-diagram}
\end{figure}

In this paper we extend the methods  in~\cite{KLO05,KLO08,KLOS13} and introduce a way to simultaneously reconstruct the conductivity, the contact impedance, and the boundary shape in absolute EIT imaging in 2D. The reconstructed conductivity coincides with the original conductivity up to a conformal deformation. 
Briefly, the new method consists of the following steps (see Figure~\ref{fig:mappings-diagram}): 
\begin{itemize}
	\item[\emph{(i)}] Choose a model domain $\Omega_m$ that approximates the true domain $\Omega$. For instance, $\Omega_m$ could be a disc having approximately the same perimeter than $\Omega$.
	\item[\emph{(ii)}]  Estimate the unknown contact impedances and compute the conductivity $\gamma_a$, the least anisotropic of all conductivities producing the same voltage data in the model domain $\Omega_m$  that was measured on $\partial \Omega$. Mild assumptions ensure that $\gamma_a$ is  unique  (see section~\ref{sec:minimally-anisotropic}).
	\item[\emph{(iii)}] Transform the reconstructed conductivity $\gamma_a$ to an isotropic conductivity $\gamma_i$. This is done by solving a Beltrami equation whose coefficient is related to the conductivity $\gamma_a$. If we denote by $F_i$ the solution of Beltrami equation, then the map $x \mapsto F_i(x)$ can be interpreted as the isothermal coordinates in which the conductivity $\gamma_a$ can be represented in isotropic form.
	\item[\emph{(iv)}]  Apply a conformal map $M$ to $\Omega_i = F_i(\Omega_m)$
	such that $\Omega_c :=M(\Omega_i)$  and the conductivity $\gamma_c := \gamma_i \circ M^{-1}$ are close to known properties of $\Omega$ and of the true conductivity defined in it, respectively. Moreover, the reconstructed isotropic conductivity $\gamma_c$ defined in $\Omega_c$ is a conformally deformed image of the true isotropic conductivity defined in $\Omega$. This last step is essential for it guarantees that the reconstruction matches optimally with the original domain at the electrodes.
\end{itemize}

The rest of the paper is organized as follows. In section~\ref{sec:theory} we present the theory behind the process of shape-deforming reconstruction. Section~\ref{sec:forward-model} introduces the complete electrode model (CEM) and gives the conditions of applicability of the proposed method in the case of the CEM. The numerical implementation of the process of shape-deforming reconstruction is given in section~\ref{sec:implementation}. Then, in section~\ref{sec:experiments} we present several results based on experimental EIT data that was measured from a phantom tank. Finally, the conclusions are given in section~\ref{sec:conclusion}.

\section{Shape-deforming reconstruction: theory}\label{sec:theory}

The measurement domain is represented by a bounded plane domain $\Omega$ having a smooth enough boundary. The electric potential in $\Omega$ is denoted by $u$. More precisely, it is the unique $H^1(\Omega)$--solution of 
the conductivity problem  
\begin{align}
\nabla \cdot \gamma \nabla u &= 0    \quad  \text{in  } \Omega \label{eq:conductivity} \\
z \nu \cdot \gamma \nabla u + u &= h \quad \text{on } \partial \Omega, \label{eq:Robin}
\end{align}
where $h$ is the Robin-boundary value, $\nu$ is the unit normal vector at the boundary $\partial \Omega$ and $z$ is a function that models the contact impedance on $\partial \Omega$.
The voltage to current measurements are modeled by the Robin-to-Neumann map $R = R_{z,\gamma}$ given by
\begin{equation}\label{eq:R-to-N-map}
R: h \mapsto \nu \cdot \gamma \nabla u \vert_{ \partial \Omega},
\end{equation}
which maps the potential distribution at the boundary $\partial \Omega$ 
to the current through the boundary. The contact impedance function $z$ is assumed continuously differentiable and non-negative: $z \geq c > 0$. 

The use of the (incorrect) model domain $\Omega_m$
instead of the true measurement domain $\Omega$ can be interpreted mathematically as a deformation of the domain. This is a way to model our lack of information about the true shape of the object under imaging. Consider then what happens to the conductivity equation when $\Omega$ is deformed to $\widetilde{\Omega}$. 
Let $F:\Omega \to \widetilde{\Omega}$ be a sufficiently smooth diffeomorphism and denote $f =F\rvert_{\partial \Omega}$. 
If $u$ is a solution of the Robin boundary value problem~\eqref{eq:conductivity}-\eqref{eq:Robin}, then
$\widetilde{u}=u\circ F^{-1}$ and $\widetilde{h}=h\circ f^{-1}$ satisfy the conductivity equation
\begin{align}
\nabla \cdot \widetilde{\gamma} \nabla \widetilde{u} &= 0    \quad  \text{in  } \widetilde{\Omega} \label{eq:conductivity-tilde}, \\
\widetilde{z} \nu \cdot \widetilde{\gamma} \nabla \widetilde{u} + \widetilde{u} &= \widetilde{h} \quad \text{on } \partial \widetilde{\Omega}, \label{eq:Robin-tilde}
\end{align}
where $\widetilde{z}(x)=z(f^{-1}(x)) \| \tau \cdot \nabla(f^{-1})(x)\|$, $\tau$ is the unit tangent vector of $\partial \widetilde{\Omega}$ and $\widetilde{\gamma}$ the conductivity given by
\begin{equation}\label{eq:push-forward}
\widetilde{\gamma}(x):= \left. F_*\gamma(x)= \frac{F^{'}(y)\gamma(y)(F^{'}(y)^T)}{|\text{det} F^{'}(y)|}\right|_{y=F^{-1}(x)}, 
\end{equation}
where $F^{'}=DF$ is the Jacobian matrix of the map $F$. We say that $\widetilde{\gamma}$ is the push-forward of $\gamma$ by the diffeomorphism $F$. 
Observe that equation~\eqref{eq:push-forward} implies that even if the conductivity $\gamma$ in the true domain $\Omega$ is isotropic (scalar valued), the transformed conductivity $\widetilde{\gamma}$ in the deformed domain $\widetilde{\Omega}$ can be anisotropic (matrix valued).

Consider next what happens to the boundary measurements in the deformed domain $\widetilde{\Omega}$. If mapping $\widetilde{R}$ corresponds to conductivity $\widetilde{\gamma}$ and contact impedance $\widetilde{z}$ in $\widetilde{\Omega}$, then the boundary measurements are transformed by
\begin{equation}\label{eq:R-to-N-map-trasnform}
( \widetilde{R} h ) (x) = (R(h\circ f ))(y)\vert_{y=f^{-1}(x)}.
\end{equation}

In the proposed method, we consider the problem of estimating the (isotropic) conductivity $\gamma$ in $\Omega$ from measurements of the Robin-to-Neumann map at $\partial \Omega$. 
We assume that the shape of the true boundary $\partial \Omega$, contact impedance $z$ and map $R_{z,\gamma}$ are not known. 
Let $\Omega_m$, called the {\em model domain}, be our model for the domain $\Omega$ and let $f_m: \partial \Omega \to \partial\Omega_m$ ($m$ for model) be a diffeomorphism which models the inexact knowledge of the boundary. The data for the solution of the inverse problem consist of the boundary of the model domain $\partial \Omega_m$ and the Robin-to-Neumann $R_m$ given in~\eqref{eq:R-to-N-map-trasnform} with $f=f_m$. 

If we consider a traditional approach for solving the EIT inverse problem, then one tries to find an isotropic conductivity in $\Omega_m$ that minimizes
\begin{equation}\label{eq:EIT-traditional}
\|R_{m} - R_{z,\gamma}\|^2 + \alpha \, W(z,\gamma),
\end{equation}
where $z$ is the function that models the contact impedance, $\gamma$ is the isotropic conductivity inside the model domain $\Omega_m$, the functional $W$ is used for providing appropriate regularization, and $\alpha>0$ is the regularization parameter.

Due to the deformation done when going from the true domain $\Omega$ to the model domain $\Omega_m$, the measurement $R_m$  does not in general correspond to any isotropic conductivity in $\Omega_m$, implying that minimization of~\eqref{eq:EIT-traditional} will lead to an erroneous reconstruction.
For demonstrations of the domain model related reconstruction artifacts with practical EIT measurements, see e.g.~\cite{GHO96,KVKK97,NKKa}.

To overcome the reconstruction errors due to the unknown domain shape, we seek to reconstruct a conductivity that is close to the original conductivity up to a conformal deformation. The first step of the proposed approach is to find in $\Omega_m$ contact impedances and a minimally anisotropic conductivity that explain the measurement. 
We remark that the idea of using an anisotropic conductivity to recover an isotropic conductivity was first mentioned in~\cite{KV87}. However, the approach suggested there is different from the one considered here.

\subsection{Recovering contact impedances and minimally anisotropic conductivity}\label{sec:minimally-anisotropic}

\begin{defi}
	Let $[\gamma^{jk}(x)]_{j,k=1}^2$ be a symmetric positive definite matrix-valued conductivity with entries in $L^\infty(\Omega)$ and bounded from below almost everywhere by a positive constant.
	Denote its eigenvalues by $\lambda_1(x)$ and $\lambda_2(x)$, $\lambda_1(x) \geq \lambda_2(x)$. The anisotropy of $\gamma$ at $x$ is defined by
	\[
	A(\gamma,x) = \frac{\sqrt{\lambda(x)}-1}{\sqrt{\lambda(x)}+1}, \quad \lambda(x) =\frac{\lambda_1(x)}{\lambda_2(x)},
	\]
	and the maximal anisotropy of $\gamma$ is then defined by
	\[ 
	A(\gamma) = \sup_{x  \in \Omega} A(\gamma,x).
	\]
	The lower positive bound on the entries $\gamma^{jk}$ implies that always $A(\gamma)<1$.
\end{defi}

If $\gamma$ is such that the anisotropy $A(\gamma,x)$ is constant (with respect $x$), then $\gamma$ is a \emph{uniformly anisotropic conductivity}. A useful feature of uniformly anisotropic conductivities is that they can be expressed as 
\begin{equation}\label{eq:uni-anisotropic}
{\gamma}(x) = \eta(x) \, T_{\theta}(x) 
\begin{pmatrix}
\lambda^{1/2} & 0\\
0 & \lambda^{-1/2}
\end{pmatrix}
T_{\theta}^{-1}(x),
\end{equation}
where $\lambda\geq 1$ is a constant, $\eta(x)\in \R_{+}$ is a real-valued function and $T_{\theta}(x)$ is a rotation matrix 
\begin{equation*}
T_{\theta} = 
\begin{pmatrix}
\cos \theta & \sin \theta \\
-\sin \theta & \cos \theta
\end{pmatrix}
\end{equation*}
where the angle $\theta(x)$ defines the direction of the anisotropy. In what follows, conductivities of the form~\eqref{eq:uni-anisotropic} are denoted by ${\gamma} =  {\gamma}_{\lambda,\eta,\theta}.$ 

Based on a classical result of Strebel on the existence of the extremal quasiconformal map (see~\cite{BLMM98,S76}), it was shown in~\cite{KLO08} that among all anisotropic conductivities defined in $\Omega_m$ and that match the observed Robin-to-Neumann map $R_m$, there is a unique conductivity $\gamma_a$ that has the minimal anisotropy $A(\gamma_a)$ and that $\gamma_a$ is of the form~\eqref{eq:uni-anisotropic}. Indeed, consider all pairs $(\widetilde{z},\widetilde{\gamma})$ of a contact impedance $\widetilde{z} \! : \! \partial\Omega_m \to \R$ and an anisotropic conductivity $\widetilde{\gamma}$ in $\Omega_m$ for which the map $R_{\widetilde{z},\widetilde{\gamma}}$ matches the map $R_m$ and denote by $S$ the class of these pairs, that is
\begin{equation}
\begin{split}
S = \{ (\widetilde{z},\widetilde{\gamma})\, :\ &\widetilde{\gamma} \in   L^\infty(\Omega_m, \R^{2\times 2}), \; \widetilde{\gamma}\ge c_1I, \;
\widetilde{z}\!: \partial \Omega_m \to \R \; \text{is} \; C^1\text{-smooth}, \;
\widetilde{z}\ge c_2, \\
\; 
&\text{and} \;  R_{\widetilde{z},\widetilde{\gamma}} = R_m , \; \hbox{where} \; c_1,c_2>0 \}.
\end{split}
\end{equation}
The following result is proved in~\cite{KLO08}. 

\begin{prop}\label{prop:anisotropic}
	Let $\Omega$ be a bounded, simply connected $C^2$-domain. 
	Let $\gamma \in C^2(\Omega)$ be an isotropic conductivity, $z\!:\! \partial \Omega \to \R$ be the $C^1$-smooth contact impedance function and $R_{z,\gamma}$ its Robin-to-Neumann map. Assume that $\Omega_m$ is a model domain satisfying the same regularity
	assumptions as $\Omega$, and let $f_m \!:\! \partial \Omega \to \partial \Omega_m$ be a $C^{2}$-smooth diffeomorphism.
	Assume that we are given $\partial \Omega_m$ and $R_{m}$ defined by~\eqref{eq:R-to-N-map-trasnform} with $f=f_m$. 
	Then the minimization problem in \(\Omega_m\),
	\begin{equation}\label{eq:constrained-minimization-1}
	\min_{(\widetilde{z},\widetilde{\gamma})\in S} A(\widetilde{\gamma})
	\end{equation}
	has a unique minimizer $(\widetilde{z}_0,\widetilde{\gamma}_0)$. Moreover, let $\lambda \geq 1 $ be such that $A(\widetilde{\gamma}_0) = (\lambda^{1/2}-1) / (\lambda^{1/2}+ 1)$. Then, there are unique $\theta \in L^\infty(\Omega_m, [0,2\pi))$ and $\eta \in L^\infty(\Omega_m, \R_+)$ such that 
	\begin{equation}\label{eq:min-anisotropic}
	\widetilde{\gamma}_0 = \gamma_{\lambda,\eta,\theta}. 
	\end{equation}\
	Finally, in \(\Omega_m\) we have $\widetilde{z}_0 (x) = z(f_m^{-1}(x)) || \tau \cdot \nabla(f_m^{-1})(x)||$ and there is a unique map $F_e:\Omega \to \Omega_m$ ($e$ for extremal), depending only on $f_m$, such that $F_e \vert_{\partial \Omega} = f_m$ and 
	\begin{equation*}
	\eta(x) = \det(\widetilde{\gamma}_0(x))^{1/2} = \gamma( F_e^{-1}(x)), \qquad \text{for } \; x  \in \Omega_m. 
	\end{equation*}
\end{prop}
\begin{remark}
	If the contact impedances are known, then the assumptions about the smoothness of the conductivity, the domain and the model map can be relaxed to $\gamma \in L^{\infty}(\Omega)$,  $\Omega$ a $C^{1,\alpha}$-domain and $f_m$ a {$C^{1,\alpha}$}-smooth diffeomorphism with $\alpha>0$, respectively.
	Note also that in the above proposition the mappinng \(f_m\) is not generally known, and hence neither is \(F_e\) recoverable. However, we can recover the contact impedance \(\widetilde z\) and the (generally anisotropic) extremal conductivity matrix \(\widetilde \gamma\) in \(\Omega_m\).
\end{remark}

In the following, the anisotropic conductivity defined by equation~\eqref{eq:min-anisotropic} is denoted by $\gamma_a$. Proposition~\ref{prop:anisotropic} can be interpreted such that one can find
unique contact impedances defined in $\partial \Omega_m$ and a unique conductivity defined in the model domain $\Omega_m$ that is as
close as possible to being isotropic and the square root of the determinant of this conductivity gives a deformed image of the original conductivity of the true measurement domain $\Omega$ in the model domain $\Omega_m$. Further, the deformation $F_e$ of the conductivity image depends only on the map $f_m$, that is, on the error in the boundary model, not on the original conductivity $\gamma$ in the true domain $\Omega$.

In the practical implementation of the algorithm for EIT measurements, the equality constrained problem~\eqref{eq:constrained-minimization-1} is approximated by the regularized least squares minimization problem
\begin{equation}\label{eq:constrained-minimization-2}
\min_{z>0,\lambda \geq 1,\eta>0,\theta} \big\{ \|R_m - R_{z,\gamma_{\lambda,\eta,\theta}} \|^2_{ L(\small{H^{1/2}(\partial \Omega_m) , H^{-1/2}(\partial \Omega_m) } ) } + \alpha \, W(z,\lambda,\eta,\theta)   \big\},
\end{equation}
where $W$ is an appropriate regularization functional and $\alpha$ the regularization parameter. Details of the discretization and minimization of problem~\eqref{eq:constrained-minimization-2} in the case experimental EIT data with a finite number of measurements is given in section~\ref{sec:minimally-anisotropic-discrete}.

\subsection{Transforming the anisotropic conductivity to an isotropic conductivity}

Once the anisotropic conductivity $\gamma_a$ has been found, the next step is to find coordinates where $\gamma_a$ can be represented as an isotropic conductivity $\gamma_i$. In order to do that, we identify $\C$ and $\R^2$ and then extend $\gamma_a$ by an isotropic unit conductivity to the whole $\C$ and define $F_i: \C \to \C$ ($i$ for isotropization) to be the unique solution of the problem
\begin{align}
\overline{\partial}F_i(x)& = \mu(x)\partial F_i(x), \quad x \in \C, \label{eq:beltrami-1} \\
F_i(x) & = x+h(x), \label{eq:beltrami-2}\\
h(x)  & \to 0 \; \text{as } |x|\to \infty. \label{eq:beltrami-3}
\end{align}
where $\overline{\partial}$ denotes the so-called d-bar operator and $\mu$ the Beltrami coefficient given by
\begin{equation}
\mu(x) = \frac{\gamma_{a}^{11} - \gamma_{a}^{22} + 2i \gamma_{a}^{12}}{\gamma_{a}^{11} + \gamma_{a}^{22} + 2 \sqrt{\text{det} \gamma_a} }, \qquad \gamma_a(x) = [\gamma_{a}^{jk}(x)]_{j,k=1}^2. 
\end{equation}
The problem~\eqref{eq:beltrami-1}-\eqref{eq:beltrami-3} has unique solution since $|\mu(x)|\leq c_0 < 1$ and $\mu$ vanishes outside $\Omega_m$, see~\cite{A66}. 
The map $x \mapsto F_i(x)$ can be interpreted as the isothermal coordinates where $\gamma_a$ is characterized by being an isotropic conductivity.  We say that $\gamma_a$ is isotropized by 
\begin{equation}\label{eq:isotropization}
\gamma_i  := (F_i)_* \gamma_a
\end{equation}
according to the equation~\eqref{eq:push-forward}.

In~\cite{KLO08} is proved that the Robin-to-Neumann map $R_{z,\gamma}$ determines uniquely the contact impedance $z$. Thus, the  knowledge of $R_{z,\gamma}$ is equivalent to the knowledge of the Dirichlet-to-Neumann map. By this result and by Proposition 1.1 in~\cite{KLO10} we have the following result: 

\begin{prop}\label{prop:isotopization}
	Let $\Omega$ be a bounded, simply connected $C^2$-domain.
	Assume that $\gamma \in C^2(\Omega)$ is an isotropic conductivity and that $R_{z,\gamma}$ is its Robin-to-Neumann map. Let $\Omega_m$ be a model domain satisfying the same regularity
	assumptions as $\Omega$, and let $f_m : \partial \Omega \to \partial \Omega_m$ be a $C^2$-smooth orientation preserving diffeomorphism.
	Assume that we are given $\partial \Omega_m$ and $R_{m}$ defined by~\eqref{eq:R-to-N-map-trasnform} with $f=f_m$.
	Let $\gamma_a$ be the unique solution of the minimization problem~\eqref{eq:constrained-minimization-1}, let $F_i$ be the unique solution of problem~\eqref{eq:beltrami-1}-\eqref{eq:beltrami-3}, and let $\gamma_i$ 
	given by~\eqref{eq:isotropization}. Then, 
	\begin{equation}\label{eq:isotropization-2}
	\gamma_i(y) = \gamma(G^{-1}(y)), \qquad \text{for  } \; y \in \Omega_i = F_i(\Omega_m),
	\end{equation}
	where $G:=F_i \circ F_e: \Omega \to \Omega_i$ is a conformal map.
\end{prop}
\begin{remark}
	If the contact impedances are known, then the assumptions about the smoothness of the conductivity, the domain and the model map can be relaxed to $\gamma \in L^{\infty}(\Omega)$,  $\Omega$ a $C^{1,\alpha}$-domain and $f_m$ a {$C^{1,\alpha}$}-smooth diffeomorphism with $\alpha>0$, respectively.
\end{remark}

The previous proposition can be interpreted by saying that we can find an isotropic conductivity $\gamma_i$ that is a conformally deformed image of the true conductivity $\gamma$.

\subsection{Post-processing by using a M\"{o}bius transformation}	
By Proposition~\ref{prop:isotopization}, the map $G:\Omega\to \Omega_i$ is conformal, and thus the restriction of the map to the boundary, $g = G\rvert_{\partial \Omega}$, could stretch the length element on the boundary. In consequence, the lengths of the images of the electrodes are changed when compared to the true physical lengths. One could correct this stretching effect by applying a conformal map to $\Omega_i$. The simplest alternative would be to do post-processing of the obtained image by using a M\"{o}bius transformation of the complex plane as follows. 

Given $\Omega_i = F_i(\Omega_m)$ and electrodes $\widehat{e}_\ell$ on $\partial \Omega_i$, find a M\"{o}bius transformation $M$ such that 
\[   |M(\widehat{e}_\ell)| = |\widetilde{e}_\ell| \qquad \text{for  } \ell=1,\ldots,L, \]
where $\widetilde{e}_\ell$ denote the electrodes on $\partial \Omega_m$ (we assume that the model map for the boundary $f_m: \partial \Omega \to \partial \Omega_m$ is length preserving on the electrodes, hence the lengths of electrodes $\widetilde{e}_\ell$ on $\partial \Omega_m$ are equal to the true physical lengths, see section~\ref{sec:CEM} ).

Additionally, assume that some extra information about the true domain $\Omega$ has been given to us. For example, in a clinical situation one might be able to measure the width, height or perimeter of the patient's chest. Then, we can use this information and try to find a M\"{o}bius transformation $M$ such that 
\[d(\Omega) = d(M(\Omega_i)),\]
where $d(\Omega)$ denotes the given measured information of the true domain $\Omega$ and $d(M(\Omega_i))$ the corresponding measure of the approximated domain $M(\Omega_i)$.

Recall that a M\"{o}bius transformation is a function of the form
\[M(y) = \frac{ay+b}{cy+d} \]
where $a, b, c, d$ are any complex numbers satisfying $ad − bc \neq 0$.		
Then, the last step of the reconstruction algorithm consists in finding complex numbers $a, b, c$ and  $d$ minimizers of the following function
\begin{equation}\label{eq:minimization-mobius}
{\cal M}(a,b,c,d) =   \Big(d(\Omega) - d(M(\Omega_i))\Big)^2 +
\beta \sum_{\ell=1}^{L}  \Big( |M(\widehat{e}_\ell)| - |e_\ell| \Big) ^2,
\end{equation}
where $\beta \geq 0$ is a scaling factor.

In what follows we denote $\Omega_c = M(\Omega_i)$ and
\[ \gamma_c(y) = \gamma_i ( M^{-1} (y) ) = \gamma( \bar{G}^{-1} (y) ),\quad \text{for}\; y \in \Omega_c,\]
where $\bar{G} = M \circ G: \Omega \to \Omega_c$ is a conformal map.  That is, $\gamma_c$ is a conformally deformed image of the true conductivity $\gamma$.

\section{EIT forward model}\label{sec:forward-model}

\subsection{Complete electrode model}\label{sec:CEM}

In an EIT experiment, a set of $L$ electrodes are
attached at the boundary $\partial \Omega$ of the body and 
the EIT measurement data consists of a finite
number of voltage and current measurements taken on these electrodes.  
The most accurate and widely used mathematical model 
for the electrode measurements is the \emph{complete electrode model} (CEM)~\cite{SChI92}, which can be considered as finite-dimensional approximation of the Robin-to-Neumann map. 

Let $e_\ell \subset \partial \Omega$, $\ell=1,\dots,L$ be disjoint open paths
modelling the electrodes. The CEM is defined by the elliptic boundary value problem: 
\begin{align}
\nabla \cdot \gamma \nabla u &= 0,   \qquad  \text{in  } \Omega \label{eq:conductivity-CEM} \\
z_\ell  \nu \cdot \gamma \nabla u + u  &= U_\ell, \quad \, \ \text{on  } e_{\ell}, \;  \ell=1,\ldots, L  \label{eq:Robin-CEM} \\
\nu \cdot \gamma \nabla u &= 0, \qquad  \text{on }  \partial \Omega \setminus \cup_{\ell = 1}^{L}e_\ell  \label{eq:Neumann-CEM} 
\end{align}
where $U_\ell$ is (constant) representing electric potential on electrode $e_\ell$, $z_\ell$ is the contact impedance at electrode $e_\ell$ 
and the normal current density outside the electrodes is zero.

In this model, the currents on the electrodes are defined by
\[ I_\ell = \int_{e_\ell} \nu \cdot \gamma \nabla u(x) \, ds(x), \qquad \ell =1,\dots,L,\]
and the relation between the electrode currents and voltages are modelled 
by the map $E: \R^L \to \R^L$: 
\[ E(U_1, \ldots, U_L) = (I_1, \ldots, I_L),\]
where $E$ is called the electrode measurement matrix for $(\partial \Omega, \gamma, e_1, \dots,e_L,z_1,\dots,z_L)$.
The existence and uniqueness of the solution $(u,U)$, where $u \in H^1(\Omega)$ and $U =(U_1, \ldots, U_L )^T \in \R^L$ is guaranteed by imposing the charge conservation
$\sum_{\ell=1}^{L} I_\ell = 0$
and by fixing the ground level of the potentials
$\sum_{\ell=1}^{L} U_\ell = 0,$
for details see~\cite{SChI92}. 

Let $\Omega$ and $\widetilde{\Omega}$ $C^{1,\alpha}$-smooth domains. We say that $f:\partial \Omega \to \partial \widetilde{\Omega}$ is length preserving on $\cup_{\ell = 1}^{L}e_\ell$ if $\|\tau \cdot \nabla f(x) \|=1$ for $x \in \cup_{\ell = 1}^{L}e_\ell$, where $\tau$ is the unit tangent vector $\partial \Omega$.

The conditions for the 
applicability of the minimization problem~\eqref{eq:constrained-minimization-2} in the case of the complete electrode model~\eqref{eq:conductivity-CEM}-\eqref{eq:Neumann-CEM} are given in Proposition~\ref{prop:R-to-N-discrete} (the proof can be found in Proposition 4.1 in~\cite{KLO05}). 

\begin{prop}\label{prop:R-to-N-discrete}
	Let $\Omega$ and $\widetilde{\Omega}$ be $C^{1,\alpha}$-smooth domains and $F :  \overline{\Omega} \to \overline{\widetilde{\Omega}}$ be a $C^{1,\alpha}$ diffeomorphism, $e_\ell \subset \partial \Omega$ be disjoint open sets, and $\gamma$ be a conductivity in $\Omega$. Let $f \! = \! F\vert_{\partial \Omega}$, $\widetilde{e}_\ell = f(e_\ell)$ and $\widetilde{\gamma} = F_*\gamma$. 
	Assume that f is length preserving on $\cup_{\ell = 1}^{L}e_\ell$. Then, 
	the electrode measurement matrices $E$ for $(\partial \Omega, \gamma, e_1, \dots,e_L,z_1,\dots,z_L)$ and $\widetilde{E}$ for $(\partial \widetilde{\Omega}, \widetilde{\gamma}, \widetilde{e}_1, \dots,\widetilde{e}_L, z_1,\dots,z_L)$ coincide. 
\end{prop}

Particularly, if $\widetilde{\Omega}$ is the model domain $\Omega_m$ and $f = f_m : \partial \Omega \to \partial \Omega_m$ is the model map for the boundary, then the assumption that $f$ is length preserving on electrodes implies 
that the size of the electrodes has to be known correctly. This is a highly feasible assumption, as the sizes of the electrodes can be measured precisely.
In this case, by Proposition~\ref{prop:R-to-N-discrete}, the electrode discretization $E$ of the Robin-to-Neumann map $R$ equals the corresponding discretization $\widetilde{E}$ of $\widetilde{R}$ given by~\eqref{eq:R-to-N-map-trasnform} and $f=f_m$. 
Thus, if the boundary shape is modeled incorrectly but the lengths of the electrodes are modelled correctly, then the electrode measurements do not change.

\subsection{Discretization and notation}

The numerical solution of problem~\eqref{eq:conductivity-CEM}-\eqref{eq:Neumann-CEM} is computed using the finite element method (FEM). The corresponding weak formulation and FEM discretization of the CEM in the case of anisotropic conductivities has been presented in~\cite{KLO08}. In this section we present only the notation that is used for the discretized problem. 

In the discretization, the domain $\Omega_m$ is divided into a set of $P$ disjoint image elements (square pixels) and functions $\eta$ and $\theta$ are approximated as piecewise constant of the form
\begin{equation}\label{eq:eta-theta-dicrete}
\eta(x) = \sum_{i=1}^{P}\eta_i\chi_i(x), \qquad
\theta(x) = \sum_{i=1}^{P}\theta_i\chi_i(x), 
\end{equation}
where $\chi_i$ is the characteristic function of the $i$-th pixel and $\eta_i$, $\theta_i$ are the pixel values of the unknown parameters. Using this notation, the finite dimensional approximations of $\eta$ and $\theta$ are identified with the coefficient vectors
\begin{equation}\label{eq:eta-theta-vector}
\eta = (\eta_1, \dots, \eta_P)^T \in \R^P, \qquad
\theta = (\theta_1, \dots, \theta_P)^T \in \R^P, 
\end{equation}
and $\lambda$ is a scalar parameter.

In practice, the EIT measurements are often made such that known currents
are injected into the domain $\Omega$ using some of the electrodes at $\partial \Omega$, and the electrode voltages needed to maintain the currents are measured. 
Sometimes, to avoid contact impedance related problems, the voltages may be measured using a four-point measurement where the voltage readings are recorded
only on those electrodes that are not used to inject current for that particular current injection. Thus, the measured data may contain only partial information of the matrix $E$.
To accommodate the possibility of a partial measurement into the forward model, the following notation is used for the discretized problem. Assume that the EIT experiment consists of a set of $J$ voltage vectors $V^{(j)}$ 
obtained as a response to $J$ different current patterns $I^{(j)} \in \R^L$, $j=1,\dots, J$, each fullfilling the charge conservation $\sum_{\ell=1}^{L} I_{\ell}^{(j)} = 0$. 
Typically, the elements of vector $V^{(j)}$ are the voltages (potential differences) between pairs of neighboring electrodes. Assume that each $V^{(j)}$ contains $K$ voltage readings, that is, we have $V^{(j)} \in \R^ K$. Then the forward model for a single current injection becomes
\[ V^{(j)} = P_j E^{-1}I^{(j)} + \epsilon{(j)},\] 
where $E$ is the electrode measurement matrix, $\epsilon{(j)}$ is a random vector which
models the measurement errors, and $P_j: \R^L \to \R^K$ is a measurement operator that
maps the electrode potentials to the measured voltages.

For the solution of the inverse problem, all the measurement vectors $V^{(1)} , V^{(2)}, \\ \ldots, V^{(J)}$ are collected into a single (column) vector
\[V = (V^{(1)} , V^{(2)}, \ldots, V^{(J)})^T \in \R^N, \; N=JK.\]
For the corresponding forward problem, we use the notation
\[ V_{\text{FEM}}(z,\eta,\theta,\lambda) = (V_{\text{FEM}}^{(1)}(z,\eta,\theta,\lambda) , V_{\text{FEM}}^{(2)}(z,\eta,\theta,\lambda), \ldots, V_{\text{FEM}}^{(J)}(z,\eta,\theta,\lambda) )^T \in \R^N \]
where
\[ V_{\text{FEM}}^{(j)}(z,\eta,\theta,\lambda) = P_j E^{-1}(z,\eta,\theta,\lambda)I^{(j)} \; \in  \R^K \]
corresponds to the measurement vector with current pattern $I^{(j)}$, vector of contact impedances $z$ and uniformly anisotropic conductivity ${\gamma}_{\eta,\theta,\lambda}$.

\section{Shape-deforming reconstruction: implementation}\label{sec:implementation}

\subsection{Inverse problem}
Given a vector $ V = \left( V^{(1)}, \ldots, V^{(J)} \right)^T \in \R^N$ of $J$ voltage measurements made on the electrodes on $\partial \Omega$, corresponding to known injected currents $I^{(j)}$, $j=1,\dots, J$, our aim is to simultaneously recover an unknown isotropic conductivity $\gamma$, the unknown boundary shape $\partial \Omega$ and the unknown vector of contact impedances $z \in \R^L$, based on these  current-to-voltage data. 

Recall that we assume that instead of the true domain $\Omega$ we are given an approximate model domain $\Omega_m$. Then, when we search for an isotropic conductivity 
in $\Omega_m$, we can not find any that matches the EIT measurements.
Therefore, we propose to use the theory and methods introduced in section~\ref{sec:theory}.  In what follows, we explain how to implement the methodology described in section~\ref{sec:theory} but in the context of practical EIT experiments with finite number of measurements.

\subsection{Recovering contact impedances and minimally anisotropic conductivity}\label{sec:minimally-anisotropic-discrete}

In this subsection we present how the constrained minimization problem~\eqref{eq:constrained-minimization-2} is discretized and solved numerically.  

Another useful property of uniformly anisotropic conductivities of the form~\eqref{eq:uni-anisotropic} is that ${\gamma}_{\eta,\theta,\lambda} = {\gamma}_{ \eta,\theta^{'},\lambda^{'} }$, where $ \lambda^{'} = 1 / \lambda$ and $\theta^{'}(x) = \theta(x) + \pi / 2,$. Therefore, problem~\eqref{eq:constrained-minimization-2} can be recast as finding an anisotropic conductivity such that $\lambda$ gets values $\lambda>0$ instead of $\lambda \geq 1$. Hence, the discrete version of~\eqref{eq:constrained-minimization-2} can be stated as finding the minimizer of 
\begin{equation}\label{eq:min-problem-discrete}
\min_{z>0,\eta>0,\lambda>0,\theta} J(z,\eta,\lambda,\theta), 
\end{equation}
where
\begin{equation*}
J(z,\eta,\lambda,\theta) = \| V - V_{FEM}(z,\eta,\lambda,\theta) \|^2 +W_{z}(z) + W_{\eta}(\eta)  + W_{\lambda}(\lambda) + W_{\theta}(\theta), 
\end{equation*}
and the regularization functionals are given by 
\begin{align}
W_{z}(z) &= \alpha_0 \sum_{\ell=1}^{L} z_\ell^2 + \alpha_1 \sum_{\ell=1}^{L} \sum_{j\in {\cal N}_\ell } |z_\ell - z_j |^2,  \label{eq:regularization-z} \\
W_{\eta}(\eta) &= \alpha_2 \sum_{k=1}^{P} \eta_k^2 + \alpha_3 \sum_{k=1}^{P} \sum_{j\in {\cal N}_k } |\eta_k - \eta_j |^2,  \label{eq:regularization-eta} \\
W_{\theta}(\theta) &=   \alpha_4 \sum_{k=1}^{P} \theta_k^2 + \alpha_5 \sum_{k=1}^{P} \sum_{j\in {\cal N}_k } |e^{i\theta_k} - e^{i\theta_j} |^2,  \label{eq:regularization-theta}\\
W_{\lambda}(\lambda) &= \alpha_6 (\lambda-1)^2,  \label{eq:regularization-lambda}
\end{align}
where ${\cal N}_k$ denotes the four-point neighborhood system for pixel $k$ in the pixel grid and $\alpha_0, \dots,\alpha_6$ are non-negative, scalar valued regularization parameters. Note that $i$ in~\eqref{eq:regularization-theta} denotes the imaginary unit.

The minimization of~\eqref{eq:min-problem-discrete} is done with the Gauss-Newton method combined with a line search strategy, for details see~\cite{NW06}.
In addition, to simplify the minimization problem~\eqref{eq:min-problem-discrete}, we compute the solution in two stages:
\begin{itemize}
	\item $1^{st}$ Stage: assume that $\eta$ is constant in $\Omega_m$ and minimize problem~\eqref{eq:min-problem-discrete} with respect to parameters $\eta \in \R$, $\lambda \in \R$, $\theta \in \R^{P}$ and $z \in \R^{L}$ until convergence is reached.
	
	\item $2^{nd}$ Stage: set the values of contact impedances equal to computed values $\widehat{z}$, change  $\eta$ to be the piecewise constant approximation~\eqref{eq:eta-theta-dicrete} and minimize problem~\eqref{eq:min-problem-discrete} with respect to parameters $\eta \in \R^{P}$, $\lambda \in \R$, $\theta \in \R^{P}$  using as initial condition the estimated values $\widehat{\eta}, \widehat{\theta},\widehat{\lambda}$ computed in the previous stage.
	
\end{itemize} 

To ensure the positivity of the constraints we consider the interior point method~\cite{W05}. That is, 
the functional $J$ in~\eqref{eq:min-problem-discrete} is augmented with the barrier function 
\[ B_s(z,\eta,\lambda) = \tau_s \left( \sum_{\ell=1}^L \frac{1}{z_{\ell}} + \sum_{k=1}^P \frac{1}{\eta_k} + \frac{1}{\lambda} \right), \]
where $\tau_s$ is a positive coefficient. Then, a solution of the original constrained problem~\eqref{eq:min-problem-discrete} is computed by solving a sequence of unconstrained problems of the form
\begin{equation} 
\min \, J(z,\eta,\lambda,\theta) + B_s(z,\eta,\lambda),
\end{equation}
with $\tau_s > \tau_{s+1}$ and $\tau_s \to 0$ as $s\to 0.$ 
The stopping criteria is based on the magnitude of the gradient and the decrease in the augmented functional. Once these are under predetermined thresholds or maximum number of iterations is reached, the Gauss-Newton iteration is terminated. For a detailed explanation on how to compute the Jacobian matrix of the forward map $V_{FEM}(z,\eta, \lambda, \theta)$ we refer to~\cite{KKSV00}.

\subsection{Transforming the anisotropic conductivity to an isotropic conductivity}

The next step of our method is the isotropization of the anisotropic conductivity $\gamma_a$. In order to do that, we need to find the isotropization map $F_i$ by solving the problem~\eqref{eq:beltrami-1}-\eqref{eq:beltrami-3} numerically. 
Substituting~\eqref{eq:beltrami-2} to~\eqref{eq:beltrami-1} gives
\begin{equation}\label{eq:neumann-series}
h(x) = {\cal P} [I-\mu {\cal S}]^{-1}\mu(x),
\end{equation}	
where ${\cal P}$ is the solid Cauchy transform and ${\cal S}$ is the Beurling transform. Recall that ${\cal P}$ is the inverse operator of $\overline{\partial}$ and that ${\cal S}$ transforms $\overline{\partial}$ derivatives into $\partial$ derivatives,  see~\cite{AP06,N96}.

The inverse operator in~\eqref{eq:neumann-series} is well-defined as it can be expressed as a convergent Neumann series using the fact that $|\mu(x)|< 1$. However, equation~\eqref{eq:neumann-series} is defined in the whole plane $\R^2$, therefore some sort of truncation is needed to compute a numerical solution. In order to do that, we use a periodization technique introduced in~\cite{GV97} where a periodic version of~\eqref{eq:neumann-series} is considered and its solution can be then used to compute a solution of~\eqref{eq:neumann-series}.
A detailed description of how to numerically solve problem~\eqref{eq:beltrami-1}-\eqref{eq:beltrami-3} is given in~\cite{KLOS13}.

\subsection{Post-processing by using a M\"obius transformation}

The last step of the proposed methodology consists in finding a M\"{o}bius transformation $M$ such that $\Omega_c= M(\Omega_i)$ is close to the true domain $\Omega$.

As the data for finding the transformation $M$,  we are given the lengths of electrodes $\widetilde{e}_\ell$ on $\partial \Omega_m$. Recall that the model map for the boundary $f_m: \partial \Omega \to \partial \Omega_m$ is length preserving on the electrodes, hence the lengths of electrodes $\widetilde{e}_\ell$ are equal to the true physical lengths, see section~\ref{sec:CEM}. Additionally, we can also assume that we are given some \emph{a priori} measured information $d(\Omega)$ of the true domain $\Omega$. 

Note that complex numbers $a, b, c, d$ and $\xi a, \xi b, \xi c, \xi d$ with $\xi \neq 0$ define the same M\"{o}bius map. Therefore, we restrict our attention to normalized M\"{o}bius transformations, that is $ad − bc = 1$. In this case, parameters $a,b,c$ determine parameter $d$, hence we write $d = d(a,b,c)$.
If we denote $a=m_1+i\,m_2$, $b=m_3+i\,m_4$ and $c=m_5+i\,m_6$, then problem~\eqref{eq:minimization-mobius} can be written in the form 
\begin{equation}\label{eq:minimization-mobius-1}
\min_{m\in \R^6} {\cal M}(m). 
\end{equation}
To solve the unconstrained minimization problem~\eqref{eq:minimization-mobius-1} we have considered the modified version of the Barzilai-Borwein method proposed in~\cite{R97}, where a nonmonotone line search technique that guarantees global convergence is combined with classical the Barzilai-Borwein method.

\section{Reconstructions from experimental data}\label{sec:experiments}

In this section we evaluate the feasibility of the presented method with experimental EIT data from tank measurements.

\begin{figure}[t]
  \centering
    \includegraphics[width=3.7cm,height=3.3cm]{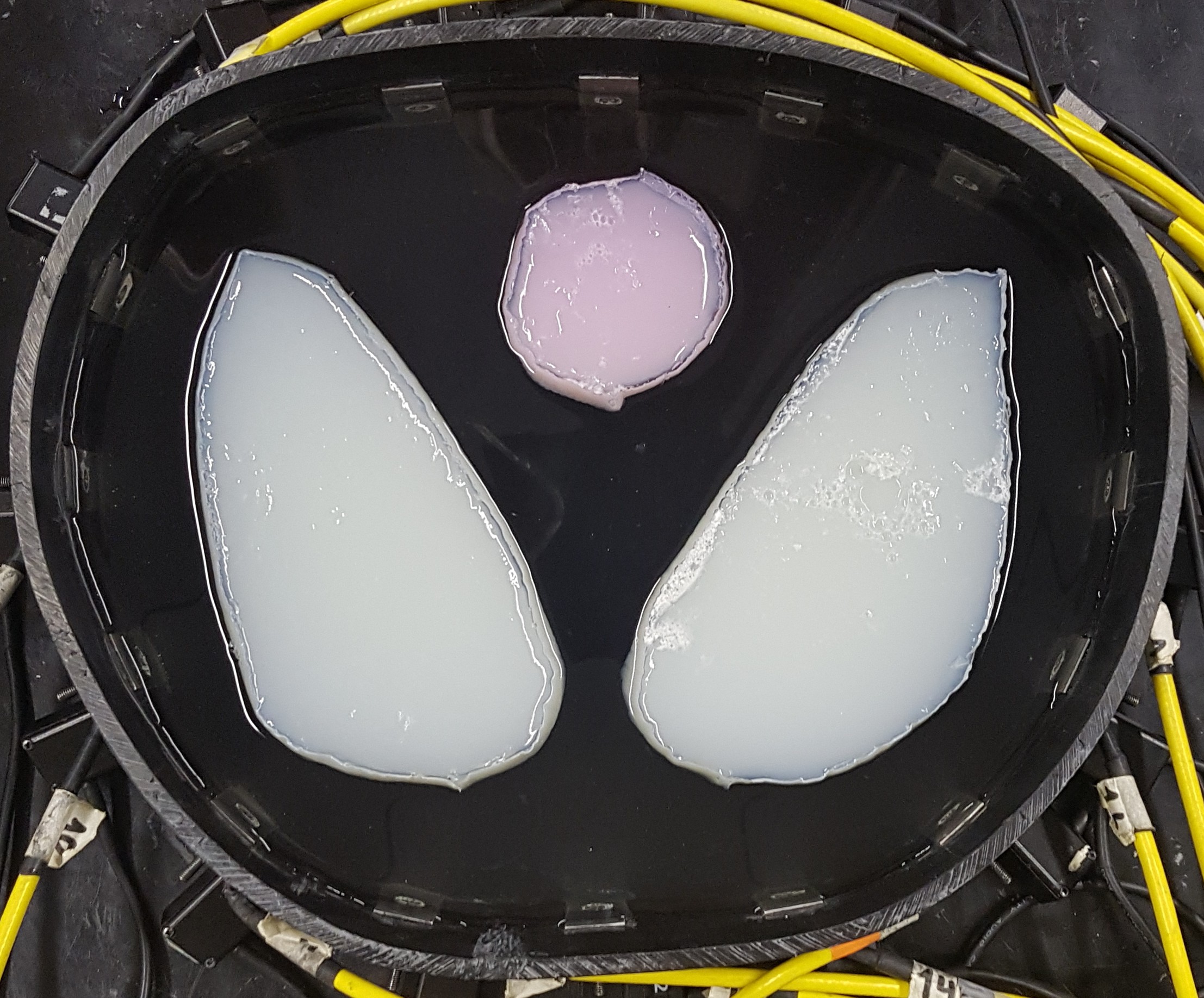}
    \quad
	\includegraphics[width=3.7cm,height=3.3cm]{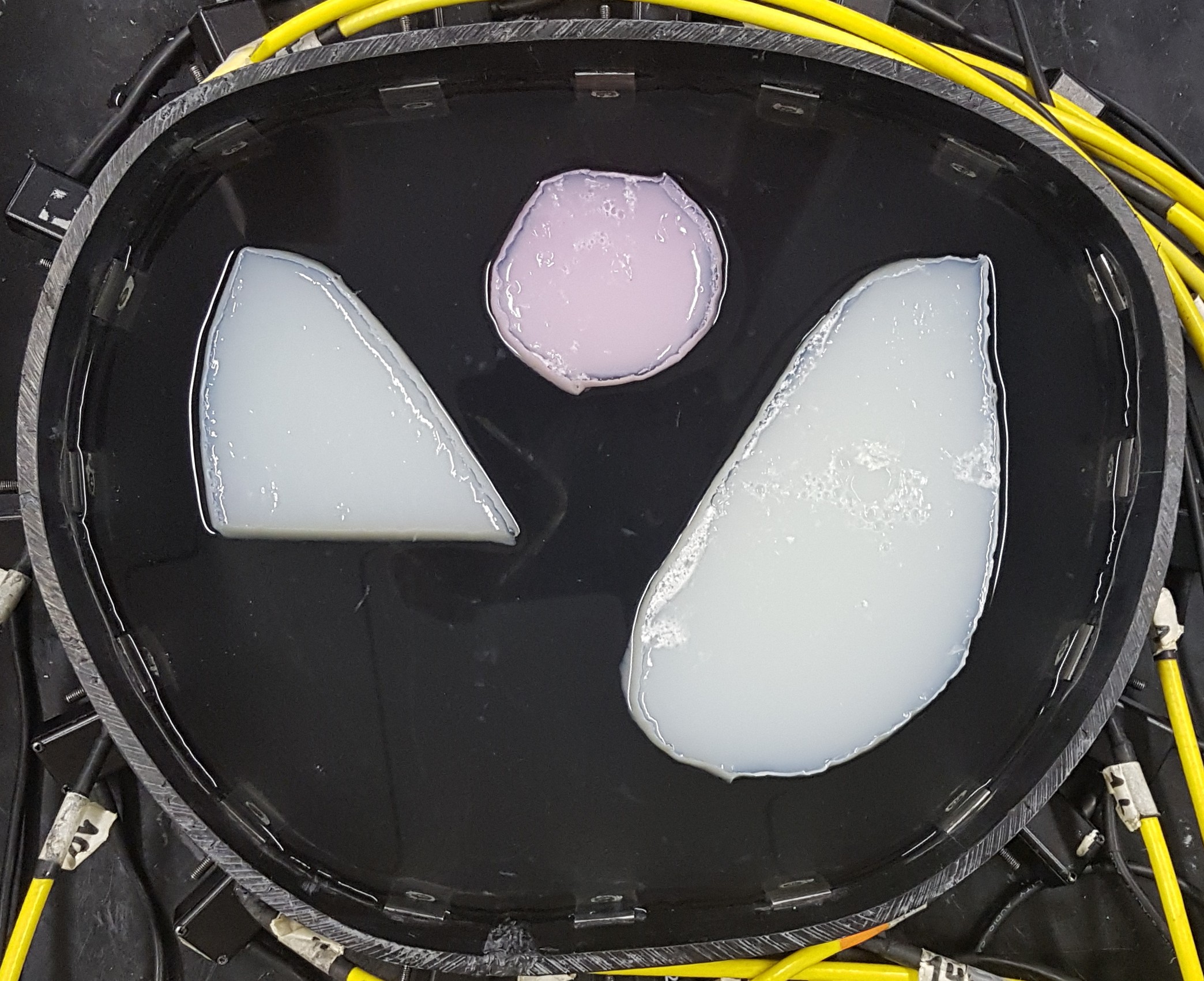}
	\quad
	\includegraphics[width=3.7cm,height=3.3cm]{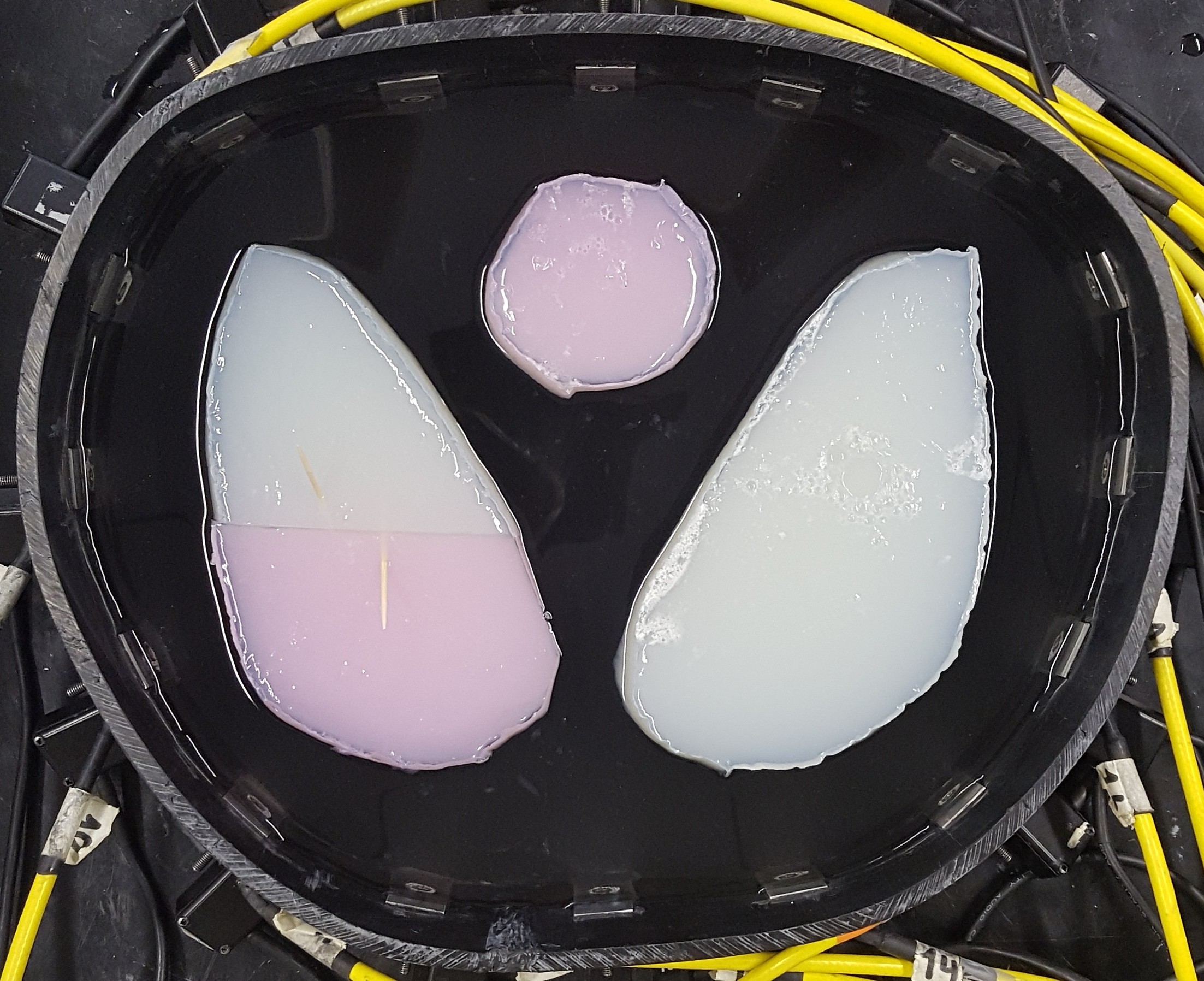}
  	\caption{The three different experimental tank setups and agar targets (pink high conductivity and white low conductivity).}
  \label{fig:chest-tank}
\end{figure}

The data sets were collected using the Kuopio impedance tomography (KIT4) equipment~\cite{KIT4} from a vertically symmetric chest shaped tank as shown in Figure~\ref{fig:chest-tank}. The chest shaped tank has $L=16$ electrodes of length $2$cm located at almost equally spaced positions on its boundary $\partial \Omega$.

The EIT measurements were collected using the adjacent (skip-$0$) method obtaining a total of $N=L^2$ voltage measurements for each EIT experiment. Therefore, having  a device with $L=16$ electrodes we have a measurement vector $V \in \R^{256}$.
The amplitude of the injected currents was $3$mA with frequency $10$kHz.

Heart and lung shaped inclusions made of agar were used to simulate the true isotropic conductivity $\gamma$. The conductivity of the saline was roughly $55\%$ higher compared to the conductivity of the lung target while the conductivity of the heart target was roughly $100\%$ higher compared to the saline.

Three different experimental data sets are considered: experimental data 1 simulates the thorax of a healthy patient, experimental data 2 simulates a human chest with an injury in the left lung (the bottom portion was extracted completely) and experimental data 3 that also simulates a human chest with an abnormality in the left lung (a higher conductivity piece of agar was included in the bottom portion),
see Figure~\ref{fig:chest-tank}.

For the three experiments two different reconstructions were computed:
\begin{enumerate}
	\item[$(i)$]  Reconstruction of isotropic conductivity in the model domain $\Omega_m$ using a traditional approach.
	\item[$(ii)$] Reconstruction of isotropic conductivity in the reconstructed domain $\Omega_c = M (F_i(\Omega_m)) \\ = \bar{G}(\Omega) $ using the proposed methodology.
\end{enumerate}

The conventional reconstruction $(i)$  was obtained by solving the problem
\begin{equation}\label{eq:min-tradional}
\min_{z>0,\gamma>0}\| V - V_{FEM}(z,\gamma) \| +W_{z}(z) + W_{\gamma}(\gamma), 
\end{equation}
where the penalty term $W_z$ is given by~\eqref{eq:regularization-z} and $W_{\gamma}$ is as in~\eqref{eq:regularization-eta}.
The solution of~\eqref{eq:min-tradional} is computed using a Gauss-Newton implementation similar to the one considered for~\eqref{eq:min-problem-discrete}.
In the first stage, we assumed $\gamma(x)=\gamma_0$ and the only unknowns were $\gamma_0 \in \R$ and $z \in \R^L$. 
In the second stage, we assume $z=\widehat{z}$, that is the contact impedances are constant and equal to the estimated values $\widehat{z}$ while $\gamma$ is written as piecewise constant function similar to~\eqref{eq:eta-theta-dicrete} on the model domain $\Omega_m$. Then, problem~\eqref{eq:min-tradional} is solved minimizing with respect to conductivity $\gamma \in \R^P$.

In all the examples of this section we assume that the perimeter of the true domain $\Omega$ is given. That is, the measured information $d(\Omega)$ is the perimeter of
chest tank and is $102.33$cm.

To measure the quality of the reconstructions of the unknown domain $\Omega$ we consider the following error
\begin{equation}\label{eq:relative-error}
E(D) = \frac{ |\Omega \setminus D| + |D \setminus \Omega| }{|\Omega|} \, \cdot \,  100\%,
\end{equation} 
where $|D|$ denotes the area of a set $D\subset \R^2$.

\subsection{Recovering: conductivity $\gamma$, boundary shape $\partial \Omega$ and contact impedances $z$}\label{subsec:disc-model-domain}

In this section we consider the problem of simultaneously recovering an unknown isotropic conductivity $\gamma$, the unknown boundary shape $\partial \Omega$ and the unknown vector of contact impedances $z \in \R^L$ from current-to-voltage data measurements made on the electrodes on $\partial \Omega$. 

\begin{figure}[htb!]
	\centering
	\setlength{\unitlength}{1cm}
	\begin{picture}(10.0,6.8) 
	\put(0.0,0.0){\includegraphics[width=9.0cm,height=6.8cm]{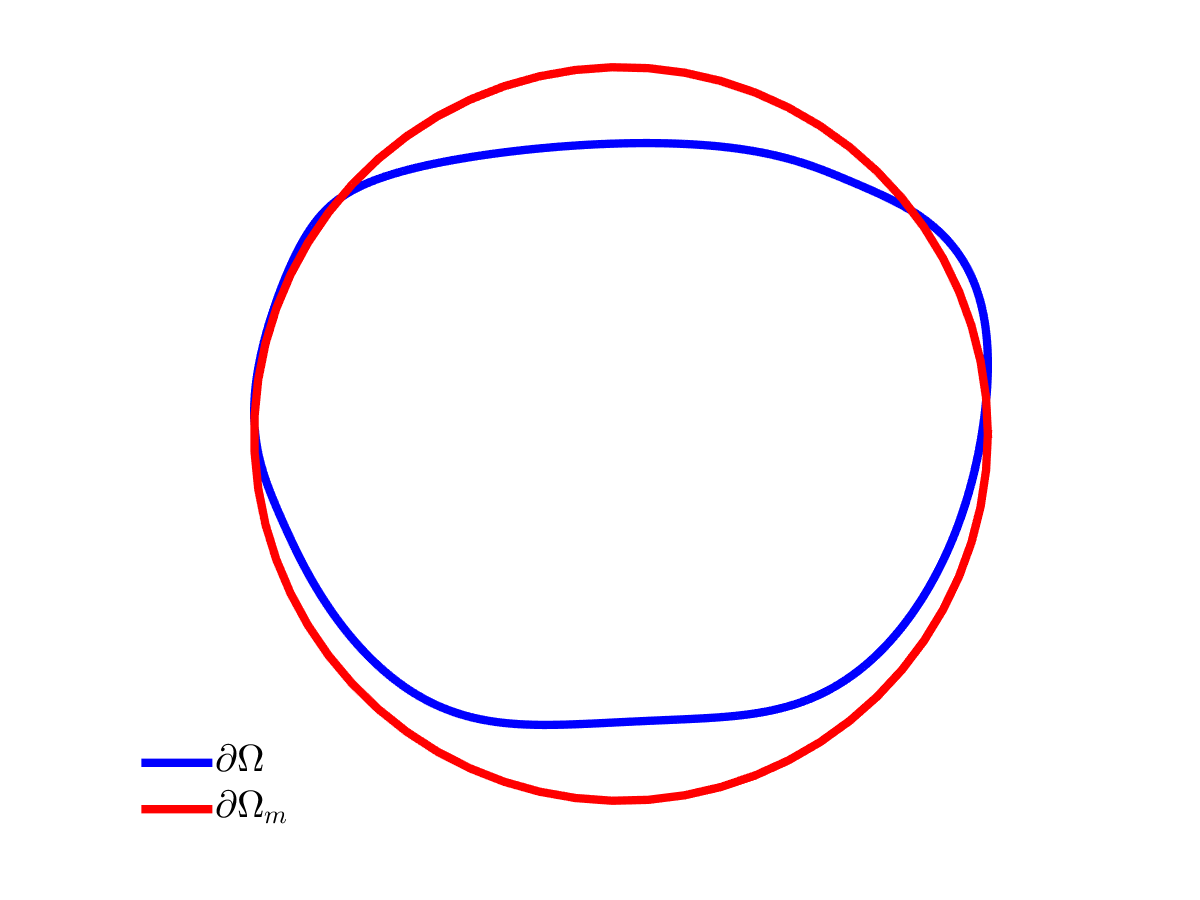}}
	\put(3.5,4.0){{Error $21.12\%$}}
	\end{picture}    
	\caption{Boundary of the true domain $\Omega$ (blue) and boundary of the model domain $\Omega_m$ (red). The relative error~\eqref{eq:relative-error} for $\Omega_m$ is $E(\Omega_m) = 21.12\%$ }
	\label{fig:model_domain}
\end{figure}

In all examples of this section the model domain $\Omega_m$ is a circle with radius $17.5$cm, see Figure~\ref{fig:model_domain}. 
The relative error~\eqref{eq:relative-error} for the model domain is $\Omega_m$ is $E(\Omega_m) = 21.12\%$. 

\begin{figure}[t]
	\centering
	\setlength{\unitlength}{1cm}
	\begin{picture}(20,4.9) 
	\put(3.0,0){\includegraphics[width=6.2cm,height=4.8cm]{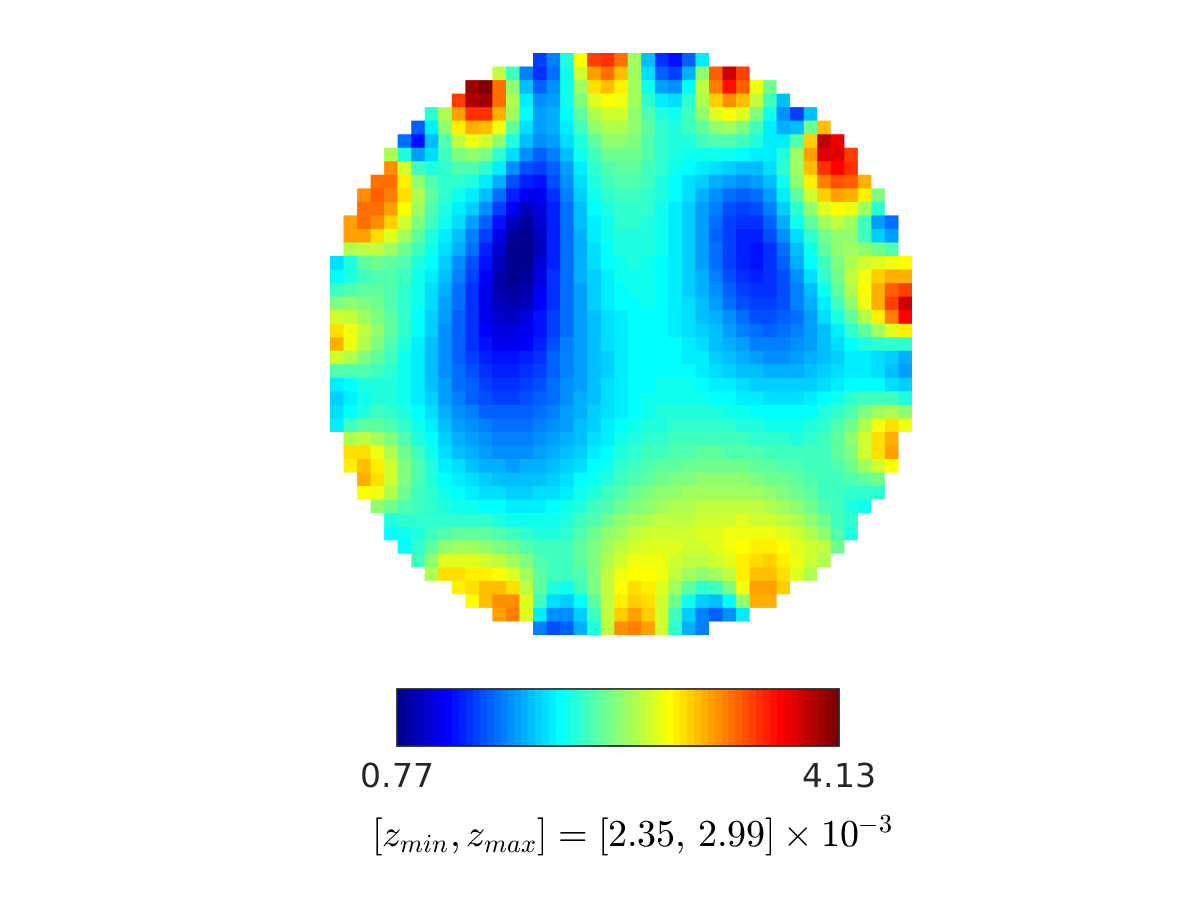}} 	
	\put(-1.5,0){\includegraphics[width=6.1cm,height=4.8cm]{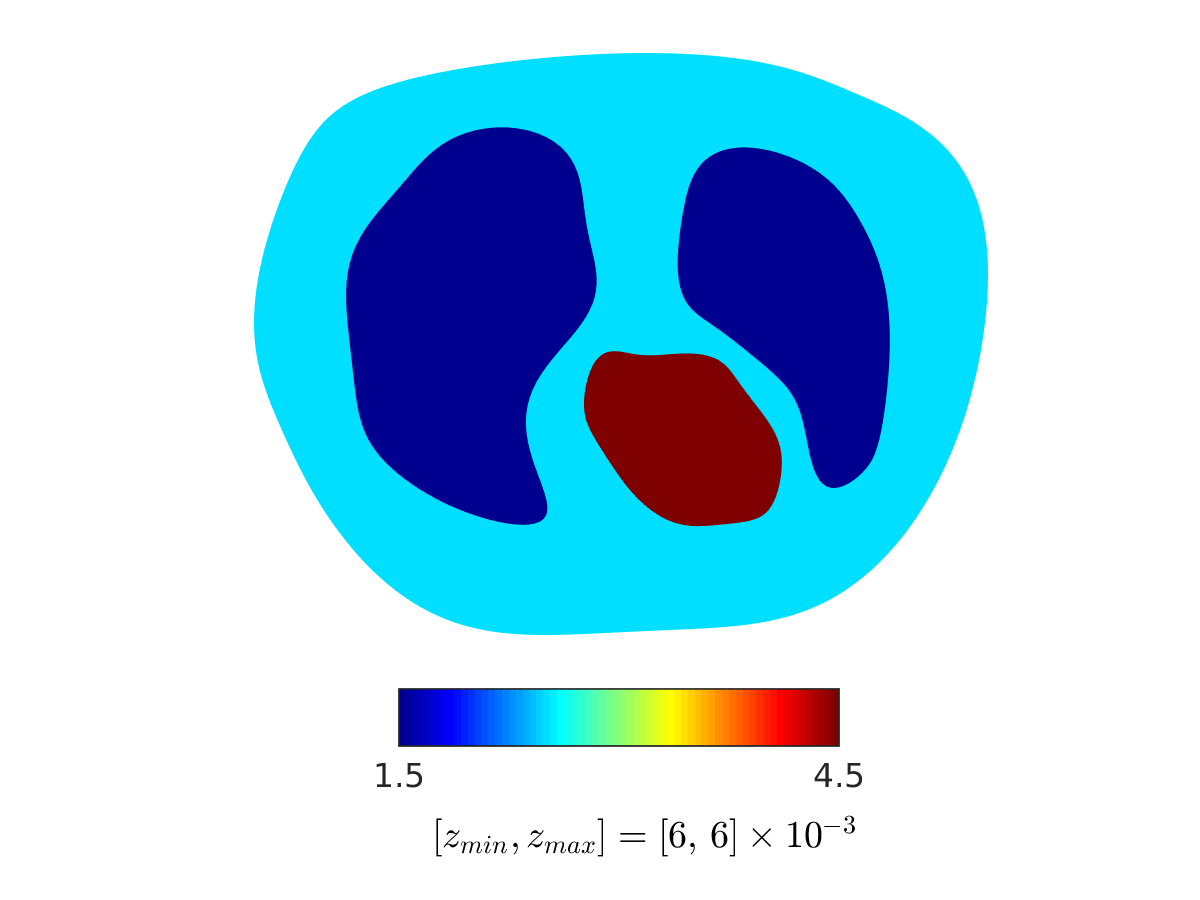}}
	\put(7.7,0){\includegraphics[width=6.3cm,height=4.9cm]{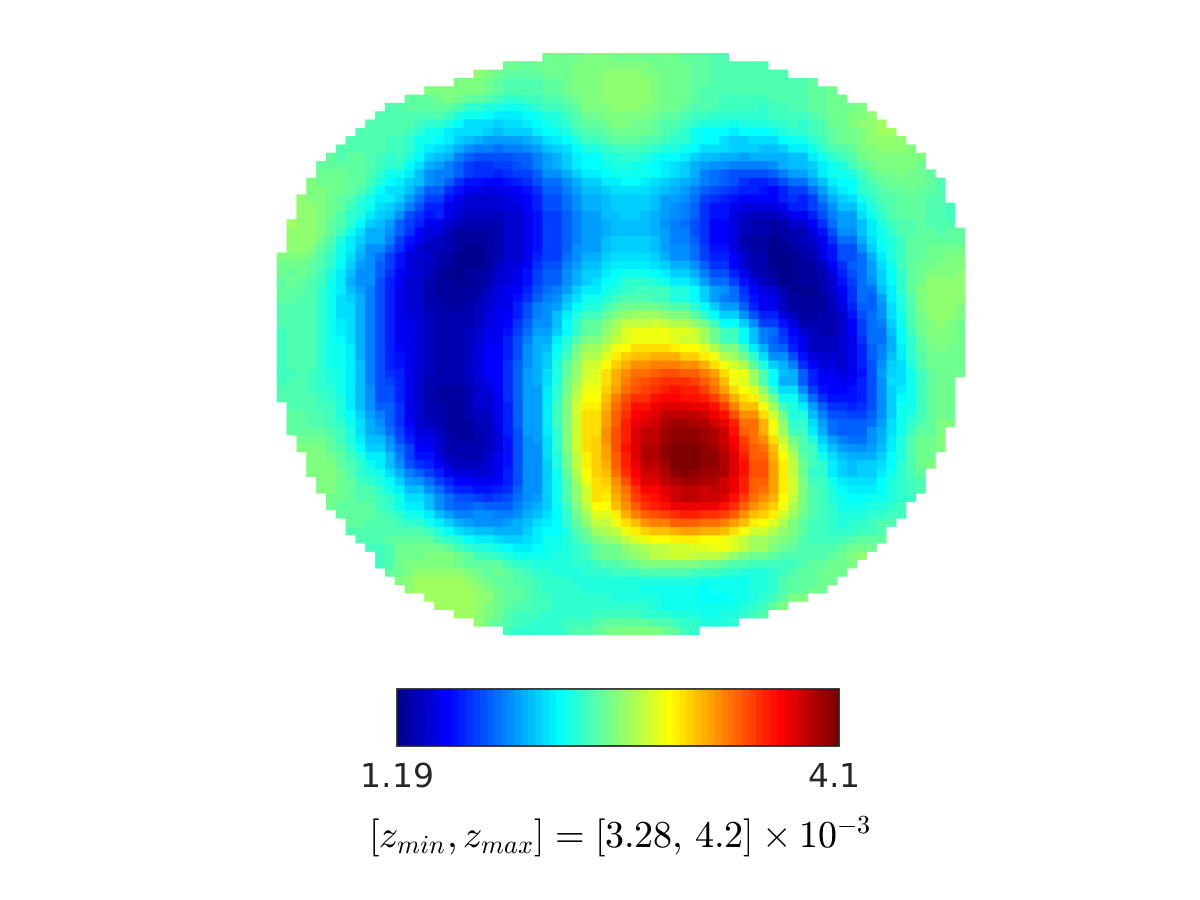}}
	\end{picture}  
	\caption{\textbf{Synthetic data}. Left: true conductivity and true domain $\Omega$. Middle: traditional reconstruction of isotropic conductivity $\gamma$ using incorrect model domain $\Omega_m$. Right: reconstruction with the proposed method using the incorrect model domain $\Omega_m$. The displayed quantity is  $\gamma_c(y) = \gamma(\bar{G}^{-1}(y))$, for $y \in \Omega_c = M (F_i(\Omega_m))$.
	The minimum and maximum value of the true and the estimated vector of contact impedances are shown below each figure, respectively.}
	\label{fig:synthetic}
\end{figure}

To calibrate the regularization parameters associated to the penalty functionals in the minimization problems~\eqref{eq:min-problem-discrete} and~\eqref{eq:min-tradional}, we first use simulated EIT measurement using the complete electrode model~\eqref{eq:conductivity-CEM}-\eqref{eq:Neumann-CEM} with $L=16$ electrodes attached to the boundary $\partial \Omega$, where $\Omega$ represents a cross-section of a human chest, and $\gamma$ simulates the lungs and heart (see Figure~\ref{fig:synthetic}). 
The parameters were tuned manually in this example for best (visual) reconstruction quality and the values of these parameters were fixed for all the experimental cases.

The simulated EIT measurements were computed using FEM. The chest shaped domain $\Omega$ was discretized using a mesh with $N_e = 19459$ triangular elements and with $N_n = 10501$ nodes for the numerical approximation of the potential $u$. 
Where did these numbers come from? The convergence of the FEM solution was studied with respect to a solution in a very dense mesh and the discretization was selected so that the error with respect to the reference FEM solution was negligible. The selected mesh was then also verified to yield high quality reconstruction from the experimental data (with the conventional reconstruction in correct domain).

In the reconstruction process, the model domain $\Omega_m$ was divided to $N_e = 11398$ triangular elements and $N_n = 6084$ node points. On the other hand, to represent the conductivity, $\Omega_m$ was discretized using $P = 2732$ square pixels of size $6\times 6$mm$^2$, leading to unknown $\gamma \in \R^{2732}$ in the minimization of\eqref{eq:min-problem-discrete}.

The image in the middle of Figure~\ref{fig:synthetic} shows the traditional reconstruction of the isotropic conductivity $\gamma$ by minimization of~\eqref{eq:min-tradional} in the incorrect model geometry $\Omega_m$.
The regularization parameters were $\alpha_0 = 10$ and $\alpha_1 = 50$ for penalty functional $W_z(z)$ and $\alpha_2 = 5\cdot 10^{-7}$ and $\alpha_3= 5 \cdot 10^{-6}$ for functional $W_{\gamma}(\gamma)$. The estimated contact impedance $\widehat{z}$ were in the range $[2.35, 2.99]\times10^{-3}$.
The image in the right of Figure~\ref{fig:synthetic} shows the reconstruction with the proposed methodology considering the incorrect model domain $\Omega_m$.
The regularization parameters for the first stage of the algorithm were $\alpha_0 = 10, \alpha_1 =  50, \alpha_4 =10^{-5} ,\alpha_5=10^{-1}  $and $\alpha_6 = 1$. 
The values of estimated contact impedance $\widehat{z}$ were in the interval $[3.28, 4.2]\times10^{-3}$.
The regularization parameters for the second stage of the algorithm were $\alpha_2 = 0,\alpha_3=10^{-7},\alpha_4=0,\alpha_5=10^{-7}$ and $\alpha_6= 10^{-5}$. Finally, the regularization parameter for the post-processing step was $\beta=0$ and therefore only the perimeter of the true domain $\Omega$ is considered in the minimization of the function given by~\eqref{eq:minimization-mobius}. In summary, only $3$ non zero parameters were considered during the second stage. The error~\eqref{eq:relative-error} corresponding to the recovery domain was $E(\Omega_c) = 5.17\%$.

\begin{figure}[t!]
	\centering
	\setlength{\unitlength}{1cm}
	\begin{picture}(20,10.0) 
	\put(0.5,6.3){\includegraphics[width=4.5cm,height=3.5cm]{Fig_ct-1.jpg}}
	\put(5.1,4.8){\includegraphics[width=7.5cm,height=5.4cm]{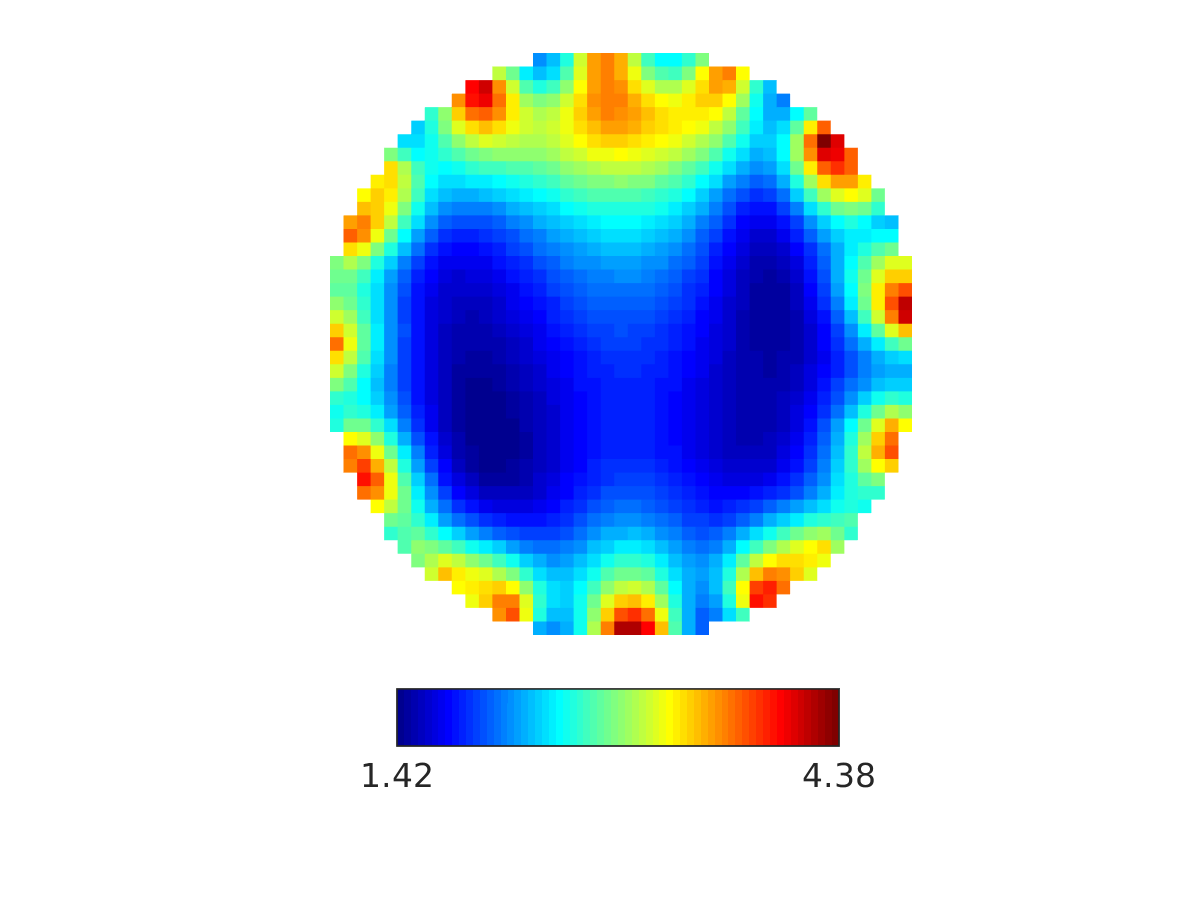}}
	\put(-1.2,-0.5){\includegraphics[width=7.7cm,height=5.7cm]{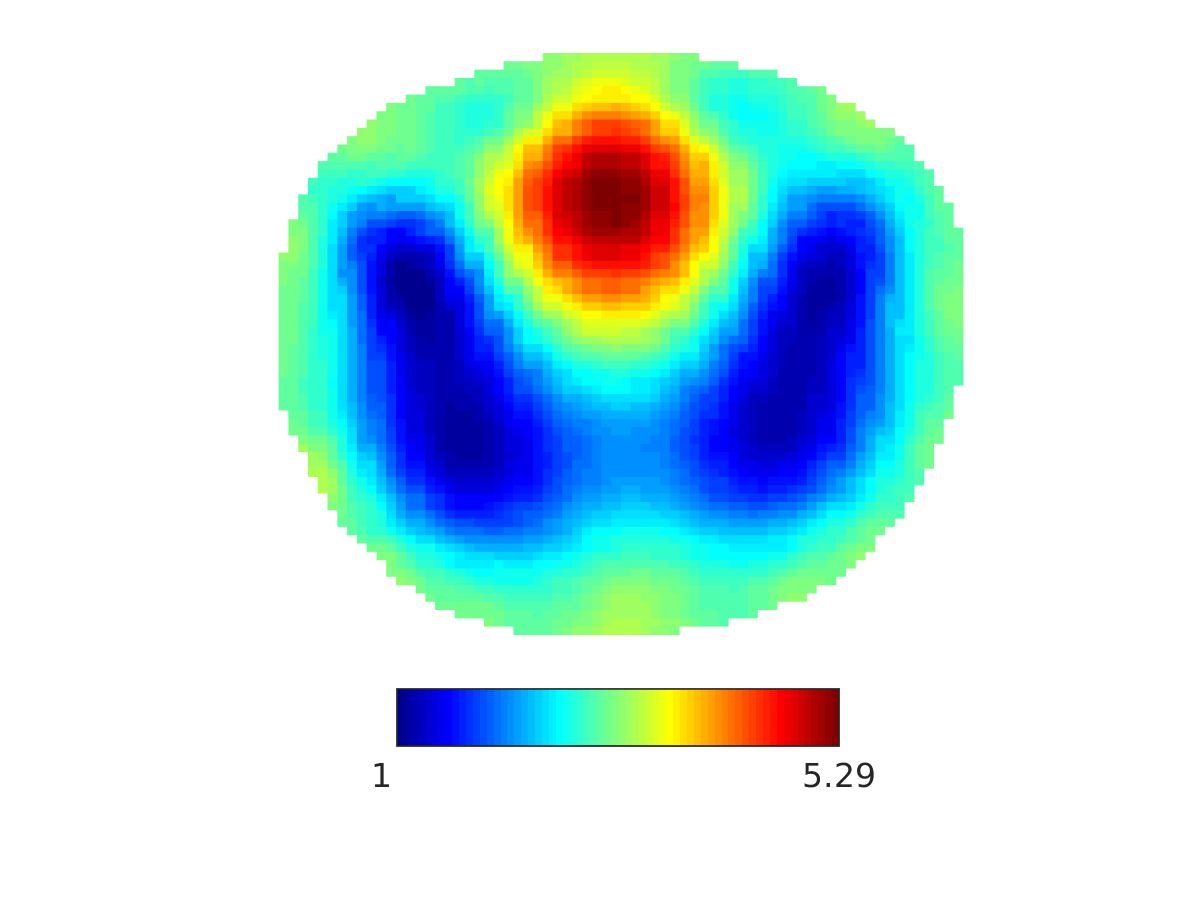}}
	\put(5.7,0.5){\includegraphics[width=6.3cm,height=4.7cm]{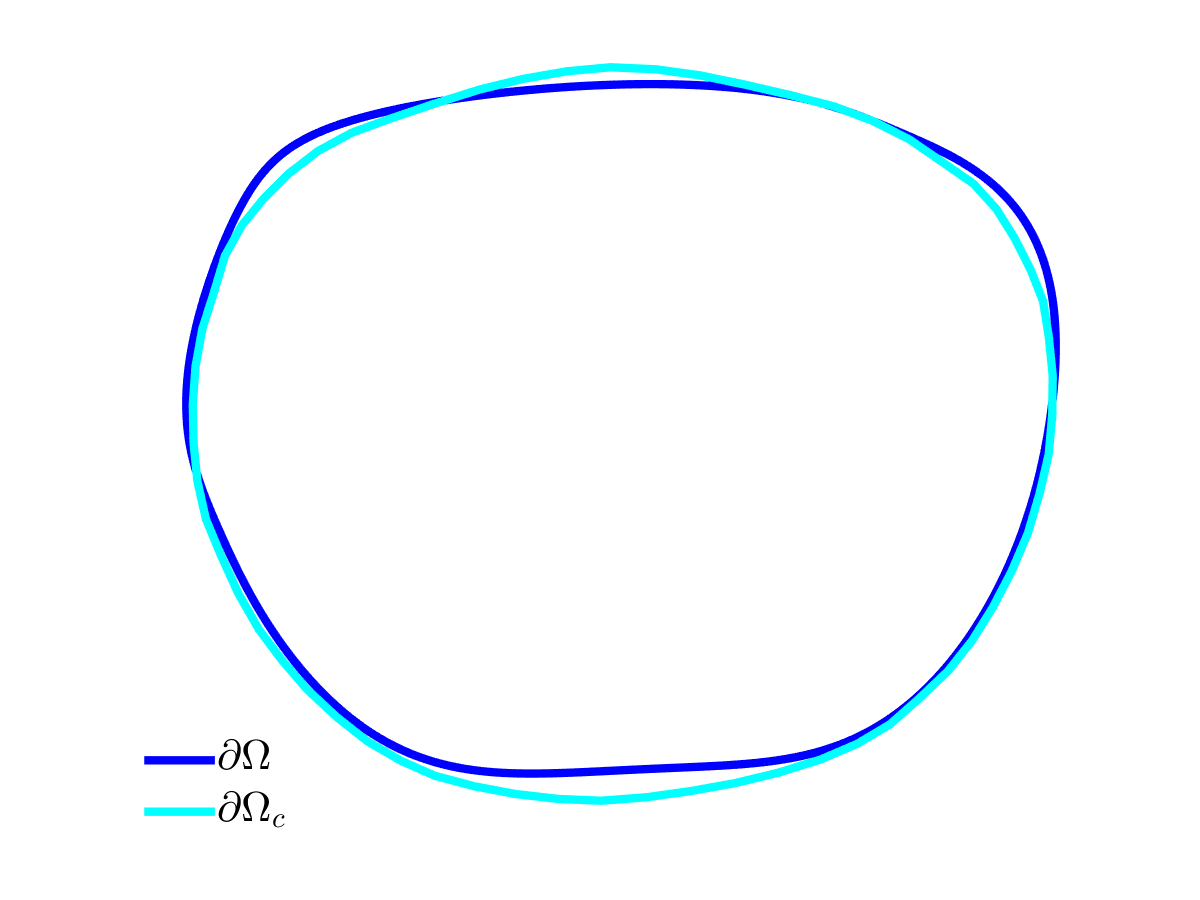}}
	\put(8.0,3.3){Error 5.18\%}
	\end{picture}    
	\caption{\textbf{Experimental data 1.} Top left: measurement setup. 
		Top right: traditional reconstruction of isotropic conductivity $\gamma$ utilizing incorrect model domain $\Omega_m$. 
		Bottom left: reconstruction with the proposed approach utilizing incorrect model domain $\Omega_m$. The displayed quantity is  $\gamma_c(y) = \gamma(\bar{G}^{-1}(y))$, for $y \in \Omega_c$.
		Bottom right: boundaries of the true domain $\Omega$ and recovery domain $ \Omega_c = M (F_i(\Omega_m)) = \bar{G}(\Omega)$. The relative error~\eqref{eq:relative-error} for the recovery domain is $E(\Omega_c) = 5.18\%$. }
	\label{fig:exp1_Oc}
\end{figure}

In Figure~\ref{fig:exp1_Oc}, we present the results corresponding to  experimental data 1.  
The measurement domain $\Omega$ and the target conductivity are shown in the top left in Figure~\ref{fig:exp1_Oc}. The top right image shows the traditional reconstruction of isotropic conductivity by minimization of~\eqref{eq:min-tradional} and considering the incorrect model domain $\Omega_m$ The values of estimated contact impedance $\widehat{z}$ are contained in the interval $[1.53, 2.73]\times10^{-3}$.
Bottom left shows the reconstruction with the proposed approach considering 
the incorrect model domain $\Omega_m$. The image shows the isotropic conductivity $\gamma_c(y) = \gamma(\bar{G}^{-1}(y))$, for $y \in \Omega_c = M (F_i(\Omega_m))$. The bottom right shows the boundary of the true domain $\Omega$ and the recovery domain $\Omega_c$. 
In this example $\widehat{z} \in [2.26, 4.84]\times10^{-3}$ and the error~\eqref{eq:relative-error} of the recovery domain is $E(\Omega_c) = 5.18\%$. 
While the conventional reconstruction of the isotropic conductivity utilizing the model domain $\Omega_m$ has serious artifacts, the proposed methodology provides an accurate reconstruction of the conductivity distribution as well as the boundary shape.

\begin{figure}[t!]
	\centering
	\setlength{\unitlength}{1cm}
	\begin{picture}(20,10.0) 
	\put(0.5,6.3){\includegraphics[width=4.5cm,height=3.5cm]{Fig_ct-2.jpg}}
	\put(5.1,4.8){\includegraphics[width=7.5cm,height=5.4cm]{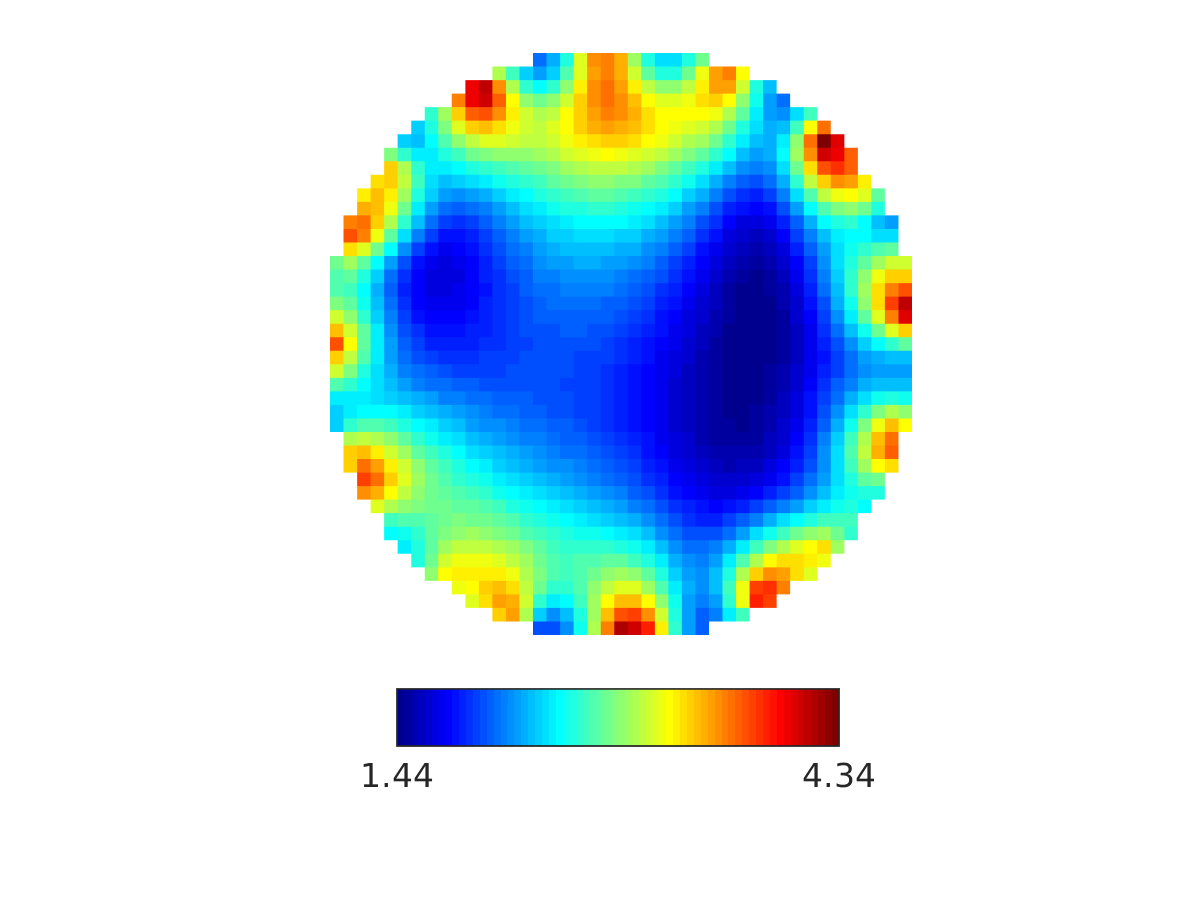}}
	\put(-1.2,-0.5){\includegraphics[width=7.7cm,height=5.7cm]{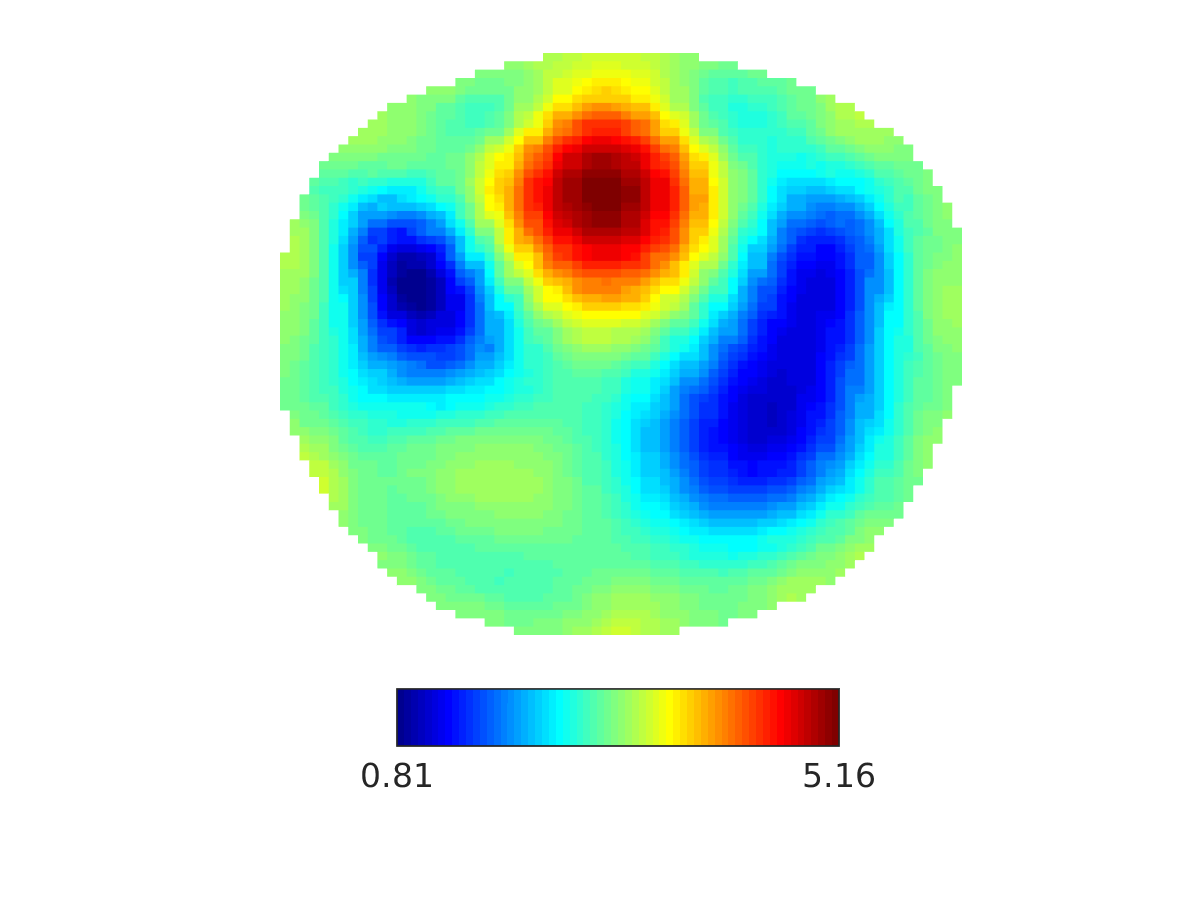}}
	\put(5.7,0.5){\includegraphics[width=6.3cm,height=4.7cm]{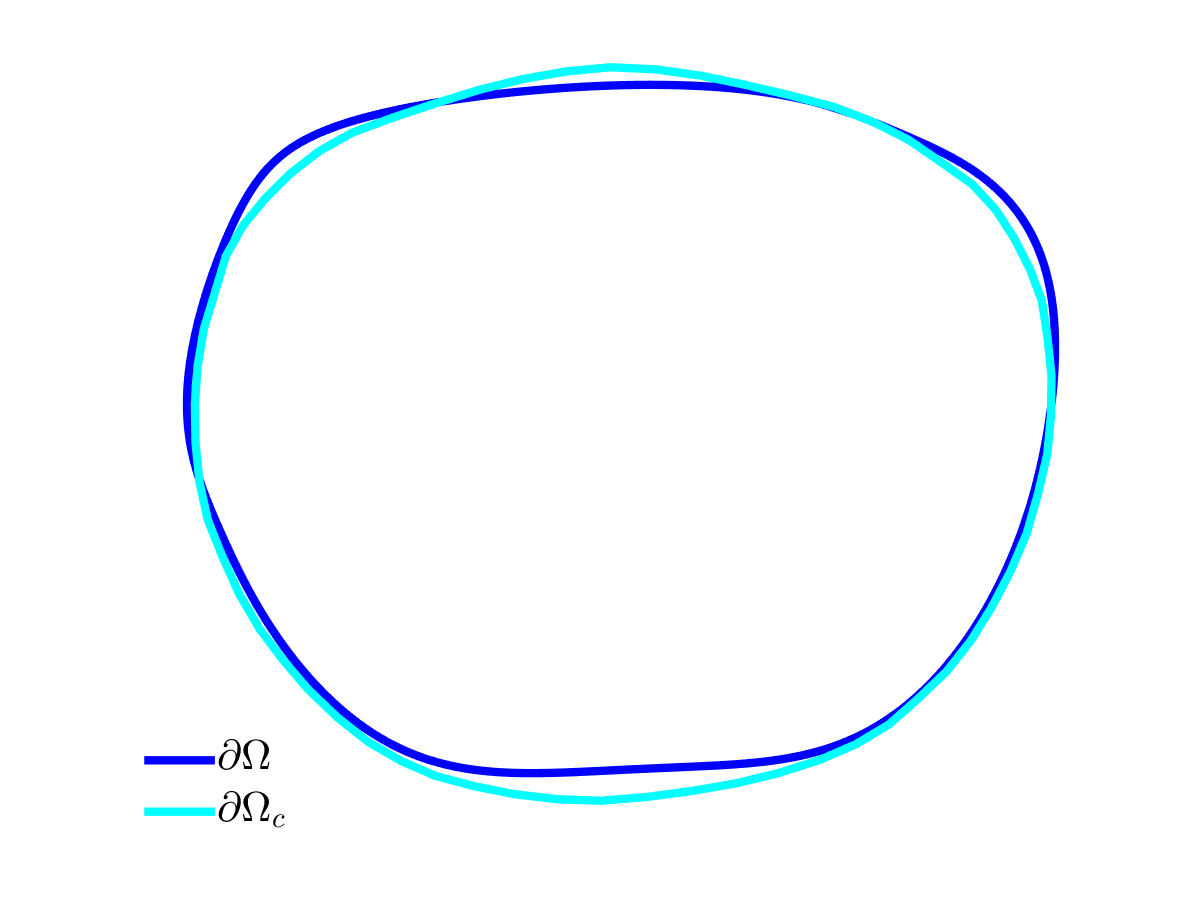}}
	\put(8.0,3.2){Error 5.27\%}
	\end{picture}    
	\caption{\textbf{Experimental data 2.} Top left: measurement setup. 
		Top right: traditional reconstruction of isotropic conductivity $\gamma$ considering incorrect model domain $\Omega_m$. 
		Bottom left: reconstruction with the proposed approach considering incorrect model domain $\Omega_m$. The displayed quantity is  $\gamma_c(y) = \gamma(\bar{G}^{-1}(y))$, for $y \in \Omega_c$.
		Bottom right: boundaries of the true domain $\Omega$ and recovery domain $ \Omega_c = M (F_i(\Omega_m))= \bar{G}(\Omega)$. The relative error~\eqref{eq:relative-error} for the recovery domain is $E(\Omega_c) = 5.27\%$.}
	\label{fig:exp2_Oc}
\end{figure}

The corresponding results for the experimental data 2 are shown in Figure~\ref{fig:exp2_Oc}. 
The values of $\widehat{z}$ were in the range $[1.61, 2.74]\times10^{-3}$
for the traditional reconstruction and in the range  $[2.49, 4.89]\times10^{-3}$ using the proposed methodology. The relative error~\eqref{eq:relative-error} for the recovery domain was $E(\Omega_c) = 5.27\%$

\begin{figure}[t!]
	\centering
	\setlength{\unitlength}{1cm}
	\begin{picture}(20,10.0) 
	\put(0.5,6.3){\includegraphics[width=4.5cm,height=3.5cm]{Fig_ct-3.jpg}}
	\put(5.1,4.8){\includegraphics[width=7.5cm,height=5.4cm]{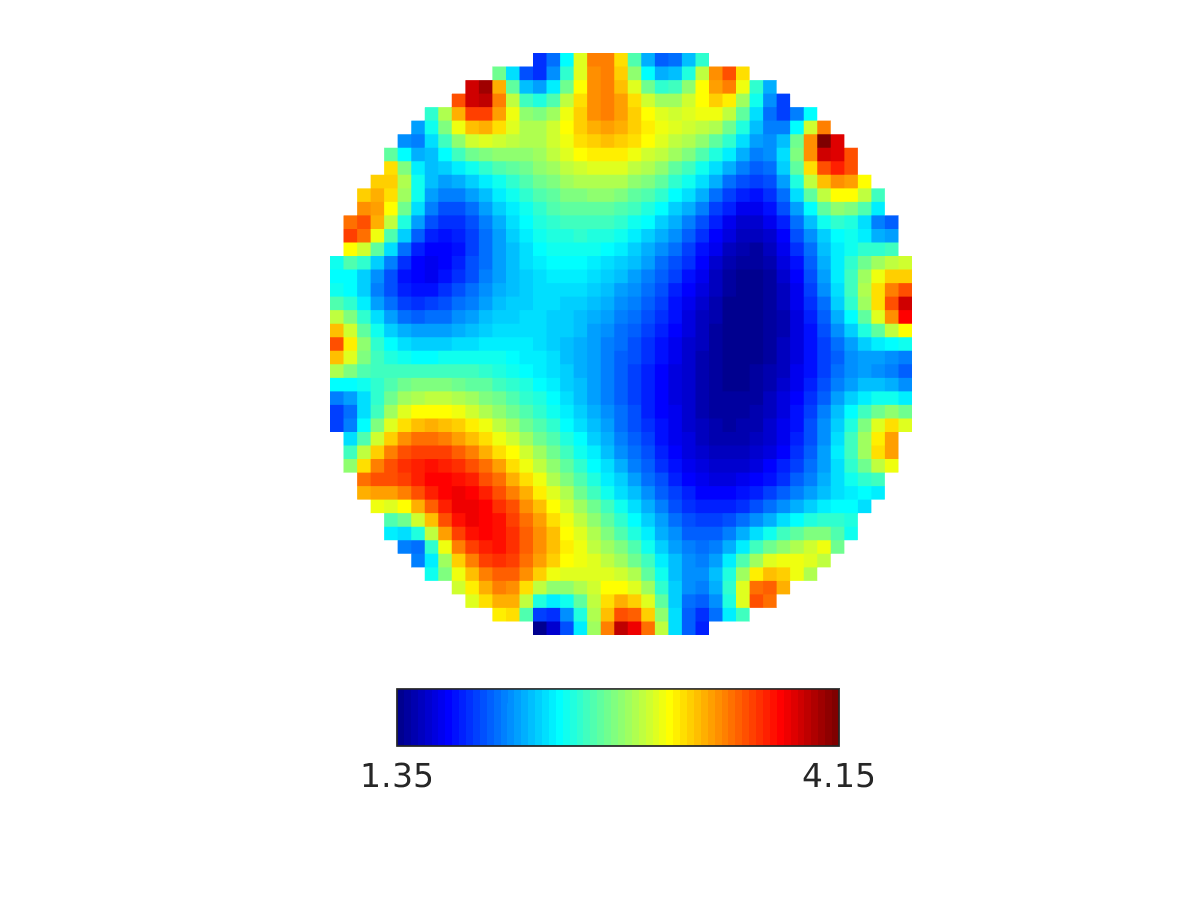}}
	\put(-1.2,-0.5){\includegraphics[width=7.7cm,height=5.7cm]{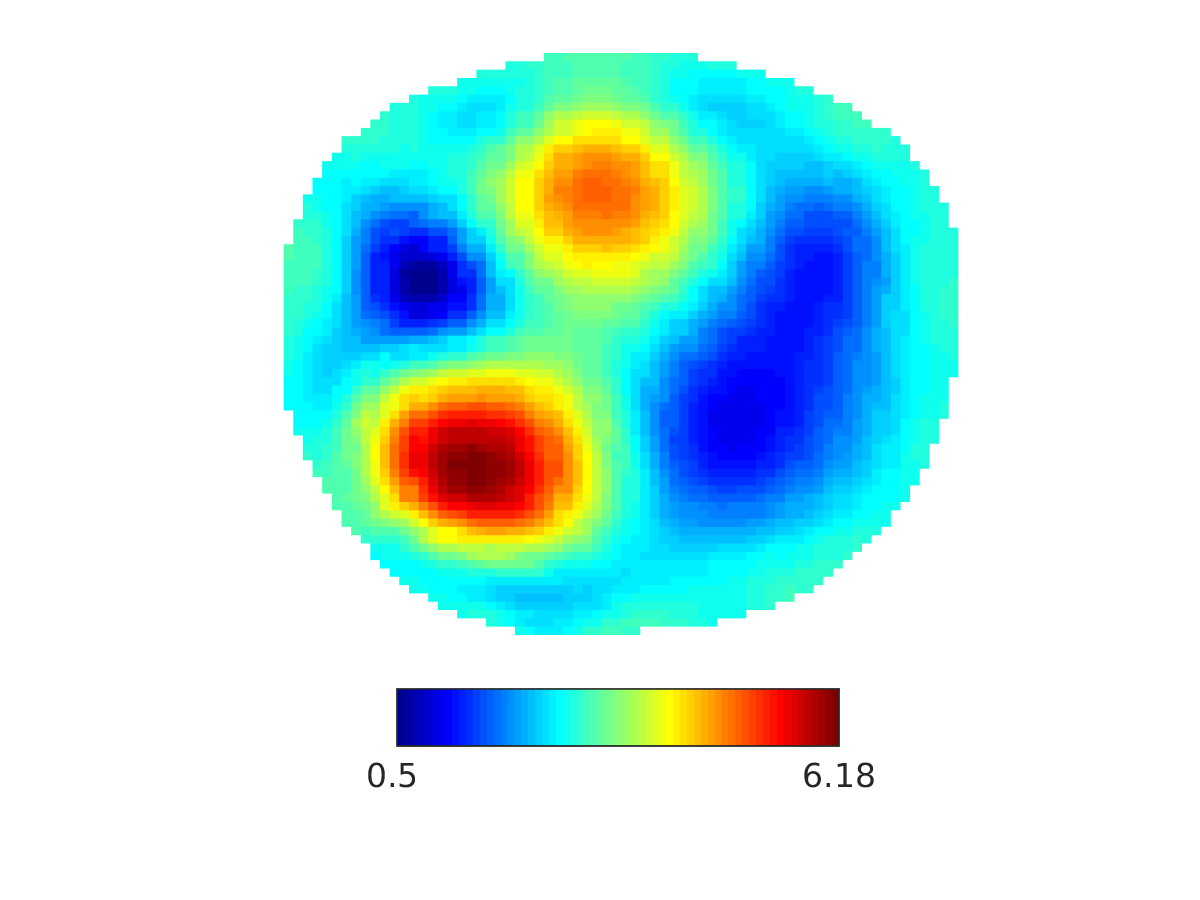}}
	\put(5.7,0.5){\includegraphics[width=6.3cm,height=4.7cm]{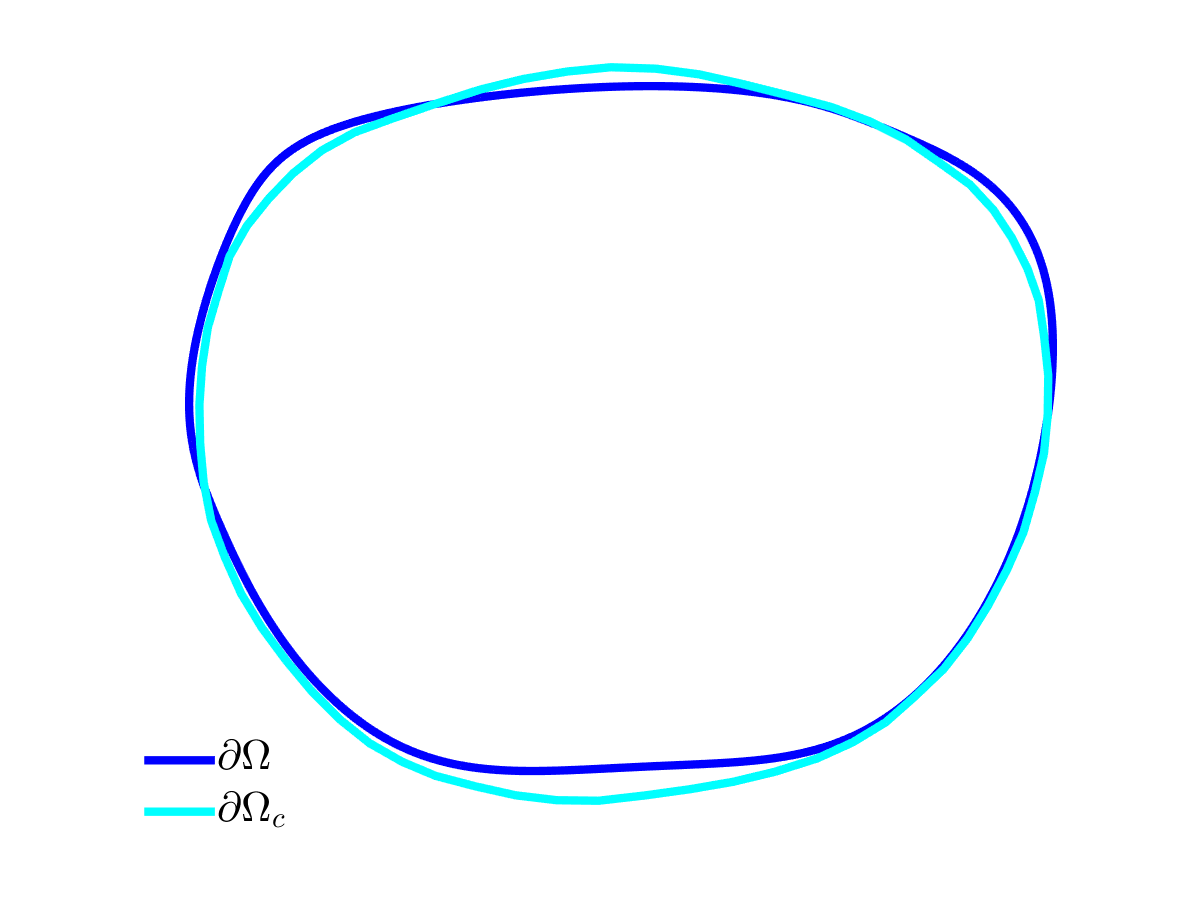}}
	\put(8.0,3.2){Error 5.57\%}
	\end{picture}    
	\caption{\textbf{Experimental data 3.} Top left: measurement setup. 
		Top right: traditional reconstruction of isotropic conductivity $\gamma$ utilizing incorrect model domain $\Omega_m$. 
		Bottom left: reconstruction with the proposed approach utilizing incorrect model domain $\Omega_m$. The displayed quantity is  $\gamma_c(y) = \gamma(\bar{G}^{-1}(y))$, for $y \in \Omega_c$.
		Bottom right: boundaries of the true domain $\Omega$ and recovery domain $ \Omega_c = M (F_i(\Omega_m))= \bar{G}(\Omega)$. The relative error~\eqref{eq:relative-error} for the recovery domain is $E(\Omega_c) = 5.57\%$.}
	\label{fig:exp3_Oc}
\end{figure}

For the experimental data 3 the results of the reconstructions are shown in Figure~\ref{fig:exp3_Oc}. The estimated contact impedance $\widehat{z}$ were in the range $[1.8, 3.37]\times10^{-3}$ for the traditional reconstruction and in the range  $[2.74, 5.47]\times10^{-3}$ using the proposed methodology. The relative error~\eqref{eq:relative-error} for the recovery domain was $E(\Omega_c) = 5.57\%$.

It is worth mentioning that the reconstructions obtained here are only qualitative since they are computed using the CEM in 2D and employing data collected from a vertically symmetric 3D target.

\subsection{Recovering: conductivity $\gamma$ and contact impedances $z$ with incorrectly modeled electrode locations}\label{subsec:eldis}

\begin{figure}
	\centering
	\setlength{\unitlength}{1cm}
	\begin{picture}(10,5.9) 
	\put(4.6,0){\includegraphics[width=7.5cm,height=5.9cm]{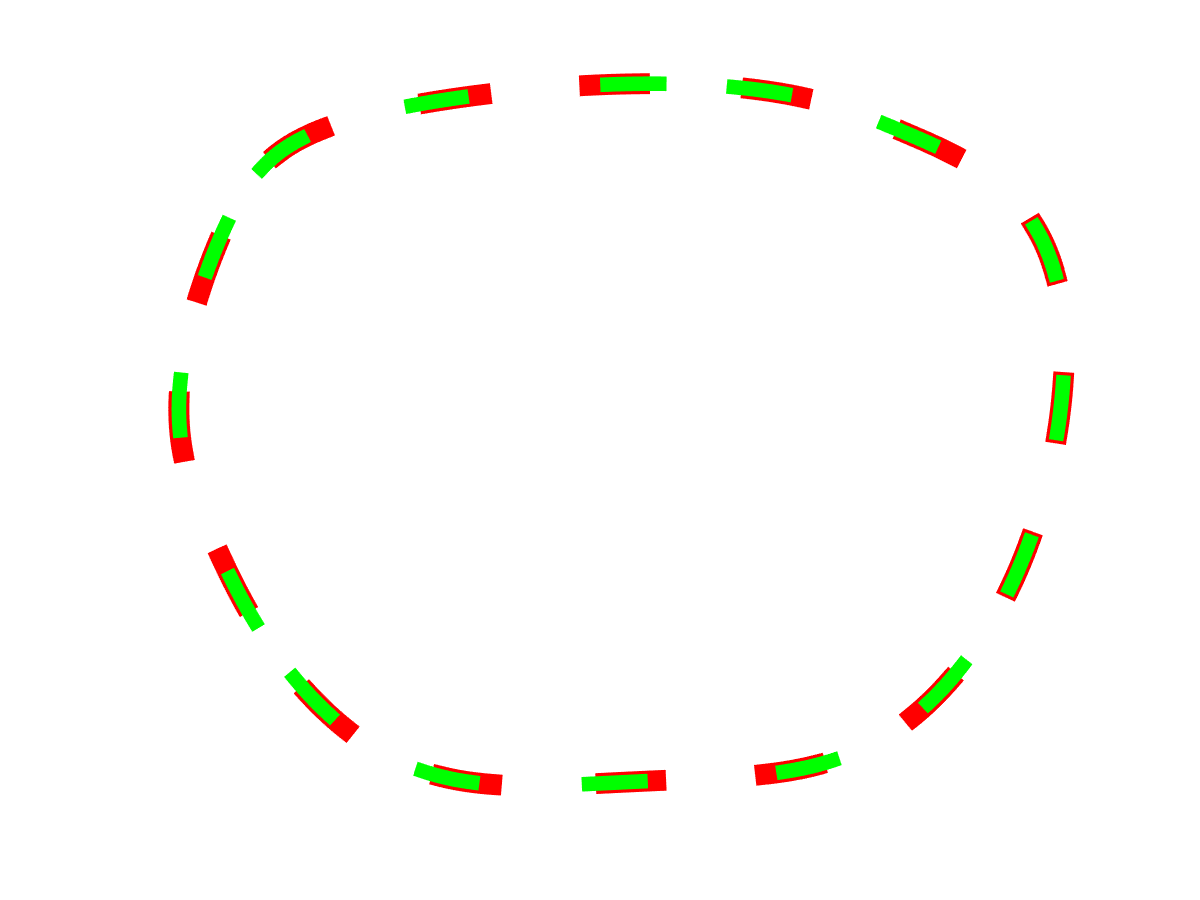}}
	\put(-2,0){\includegraphics[width=7.5cm,height=5.9cm]{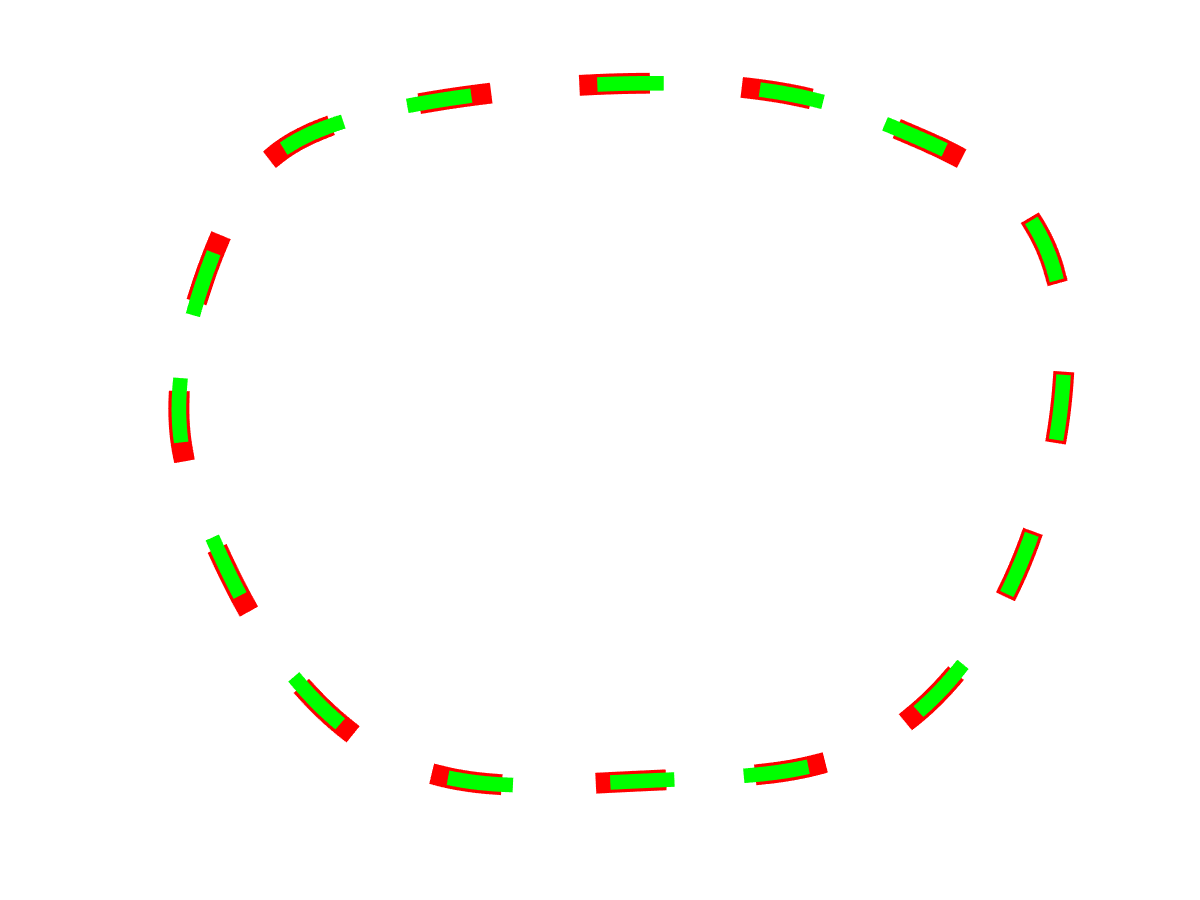}}
	\put(11.2,4.3){ \footnotesize{$e_2$}}
	\put(11.3,3.2){ \footnotesize{$e_1$}}
	\put(11.1,2.1){ \footnotesize{$e_{16}$}}
	\end{picture}    
	\caption{Left: Case 1 electrodes are displaced approximately $25\%$ of the true physical length. Right: Case 2 electrodes are displaced approximately $35\%$ of the true physical length. True electrodes are shown in red. Displaced electrodes are shown in green. In both cases, no displacement was added to electrodes $e_1$, $e_2$ and $e_{16}$.} 
	\label{fig:electrodes_eldis_35_50}
\end{figure}
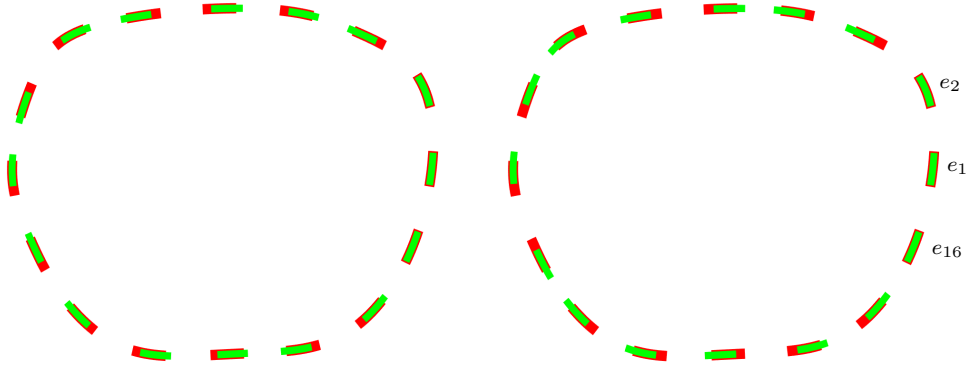

In this section we consider the problem of recovering an unknown isotropic conductivity $\gamma$ from current-to-voltage data, assuming that true domain $\Omega$ is known but the vector of contact impedances $z \in \R^L$ is unknown and the electrode locations are modelled incorrectly. This would be a realistic scenario in practice when the domain shape would be known but the electrodes would be positioned, for example, around the chest of the patient manually without access to an auxiliary measurement of their locations.   
The measurement electrodes were displaced with randomly selected sign of the displacement (with the exception that no displacement was added to electrodes $1, 2$ and $16$). We consider two case studies (see Figure~\ref{fig:electrodes_eldis_35_50} ):
\begin{itemize}
	\item[] \textbf{Case 1:} electrodes were displaced approximately $25\%$ of the true physical length.
	\item[] \textbf{Case 2:} electrodes were displaced approximately $35\%$ of the true physical length.
\end{itemize}

For each case study we present the results using the same set of experimental data used in the previous section. Similarly, for each experiment, two different reconstructions were computed:
\begin{itemize}
	\item[]$(i)$  Traditional reconstruction of isotropic conductivity in the true domain $\Omega$ using incorrect electrode locations.
	\item[]$(ii)$ Reconstruction of isotropic conductivity using the proposed methodology taking as model domain the true domain, that is $\Omega_m = \Omega$, but with incorrect electrode locations. 
\end{itemize}

The true domain $\Omega$ was discretized using mesh with $N_e = 15809$ triangular elements and $N_n = 8676$ node points. Then, to represent the conductivity, $\Omega$ was divided to $P = 2303$ square pixels of size $6\times 6$mm$^2$, leading to unknown $\gamma \in \R^{2303}$.

\begin{figure}
	\centering
	\setlength{\unitlength}{1cm}
	\begin{picture}(20,13.2) 	
	\put(-0.8,9.7){\includegraphics[width=4.0cm,height=3.2cm]{Fig_ct-1.jpg}}
	\put(3.2,8.7){\includegraphics[width=6.1cm,height=4.5cm]{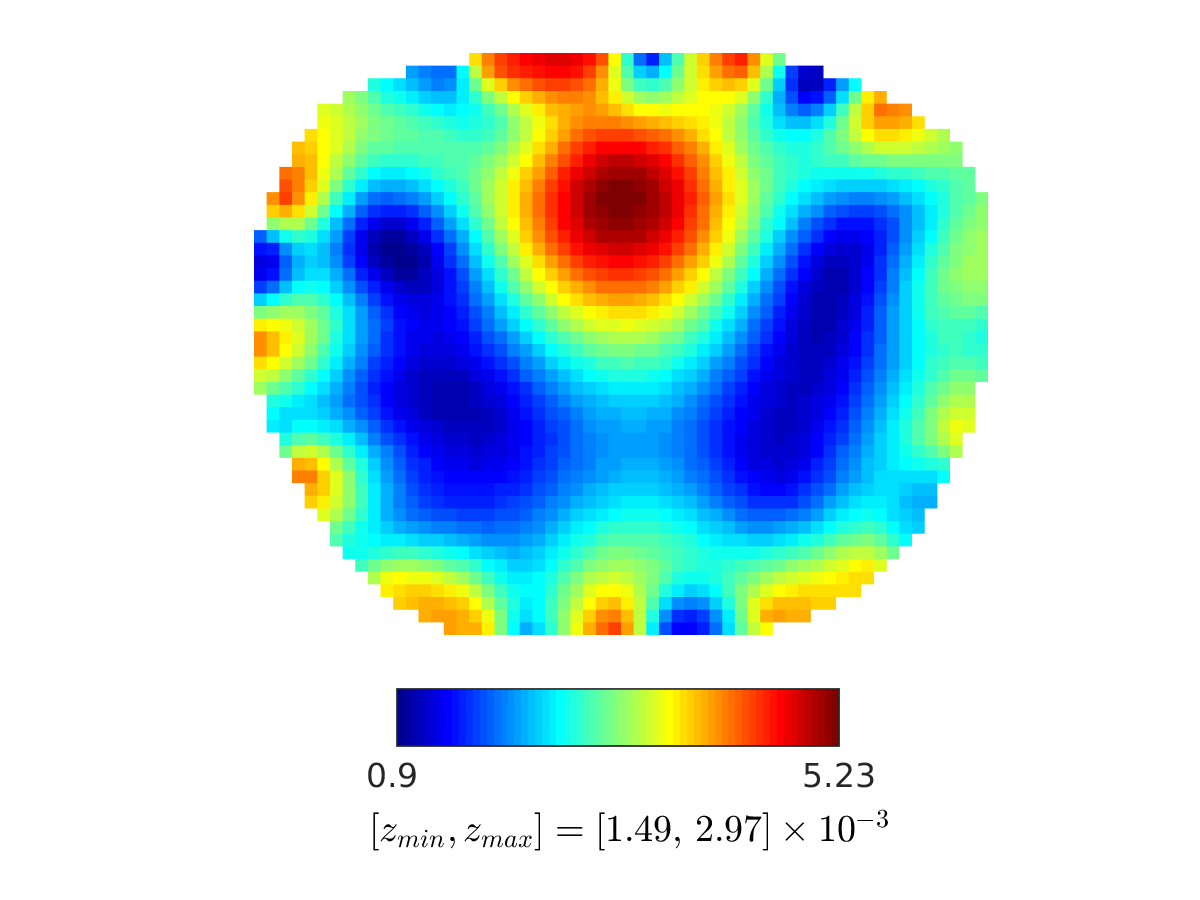}}
	\put(8.2,8.7){\includegraphics[width=6.1cm,height=4.5cm]{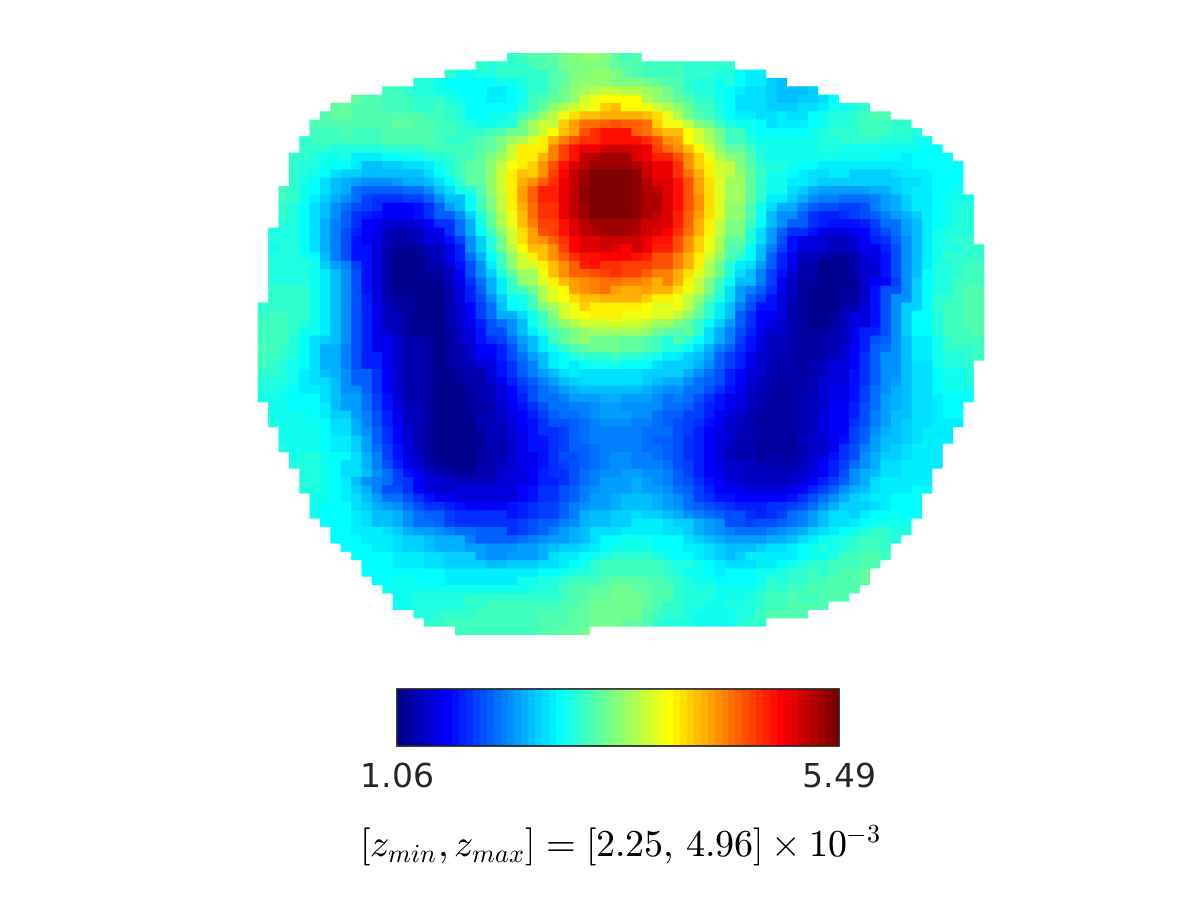}}
	\put(-0.8,5.2){\includegraphics[width=4.0cm,height=3.2cm]{Fig_ct-2.jpg}}
	\put(3.2,4.2){\includegraphics[width=6.1cm,height=4.5cm]{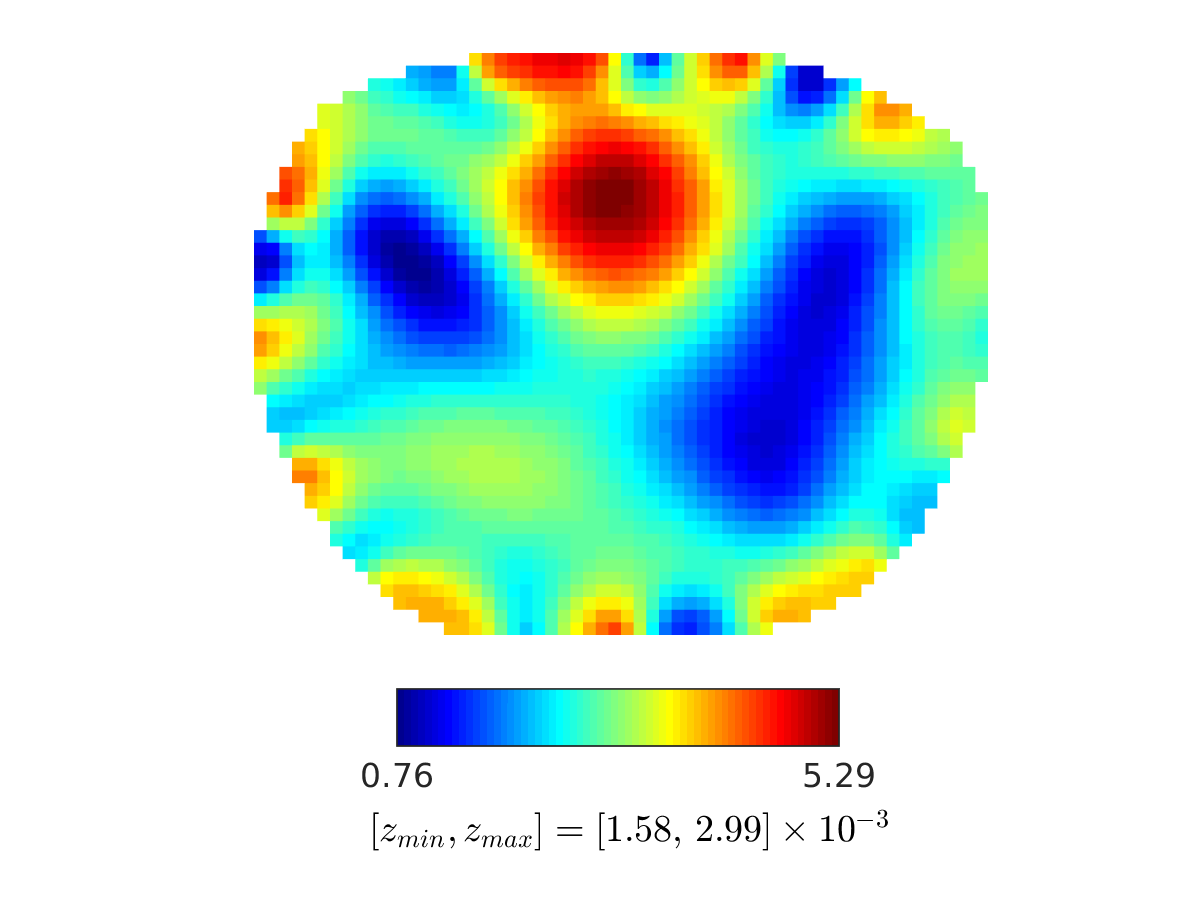}}
	\put(8.2,4.2){\includegraphics[width=6.1cm,height=4.5cm]{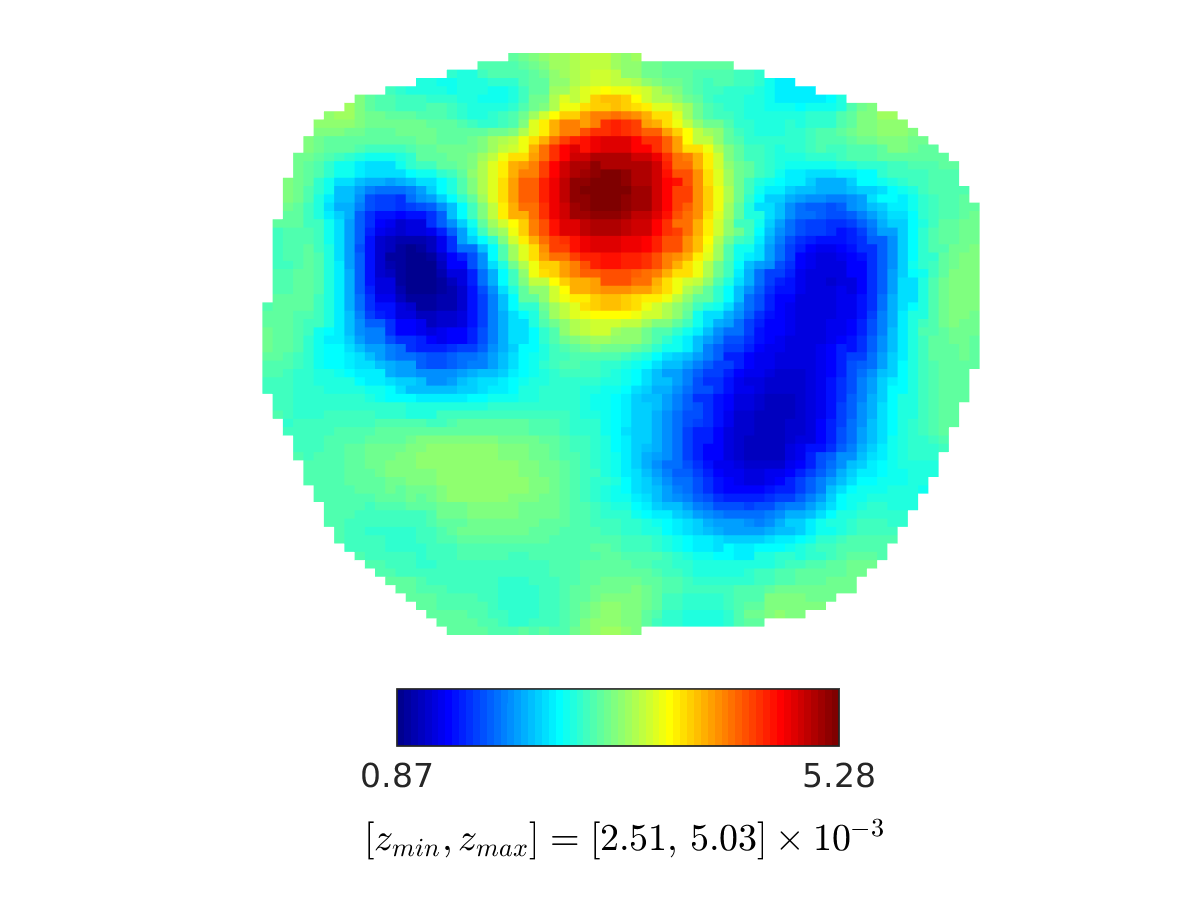}}
	\put(-0.8,0.7){\includegraphics[width=4.0cm,height=3.2cm]{Fig_ct-3.jpg}}
	\put(3.2,-0.3){\includegraphics[width=6.1cm,height=4.5cm]{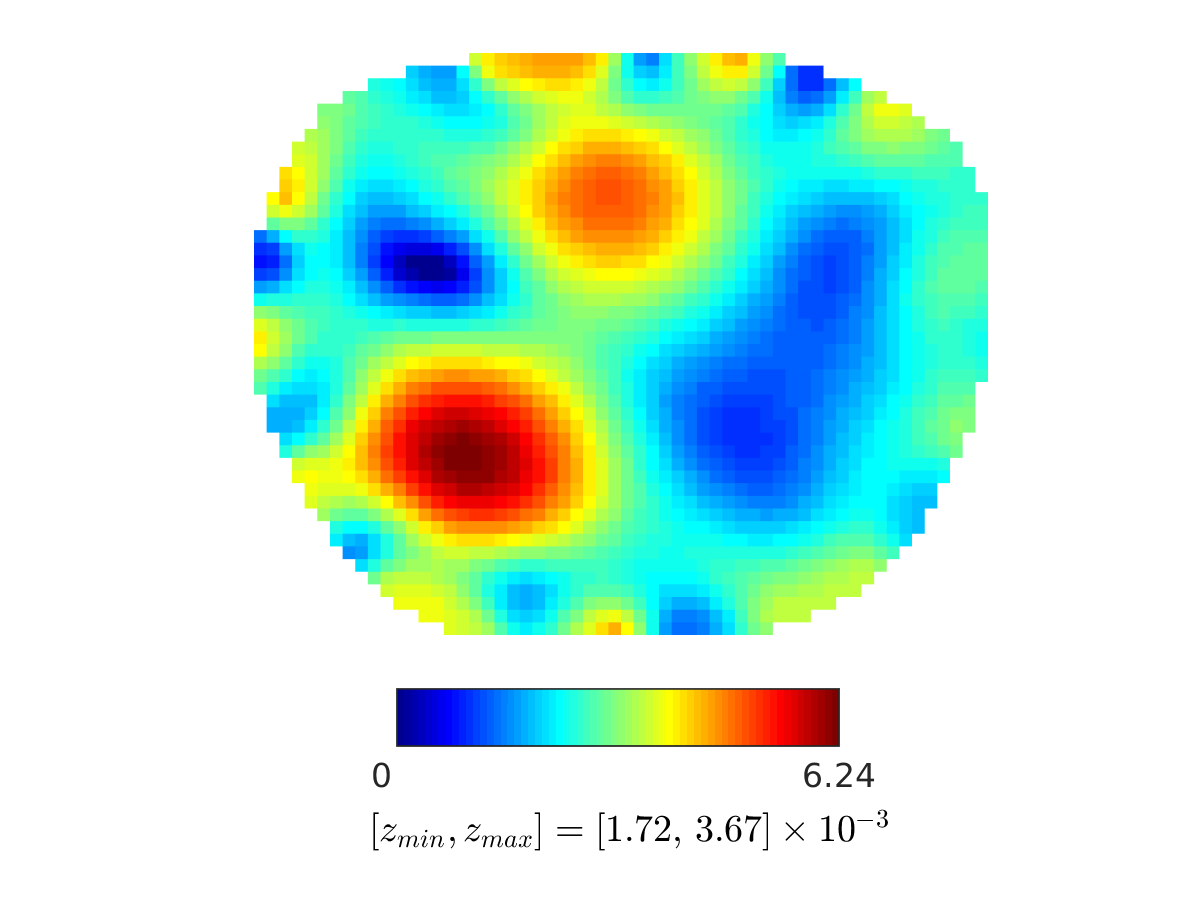}}
	\put(8.2,-0.3){\includegraphics[width=6.1cm,height=4.5cm]{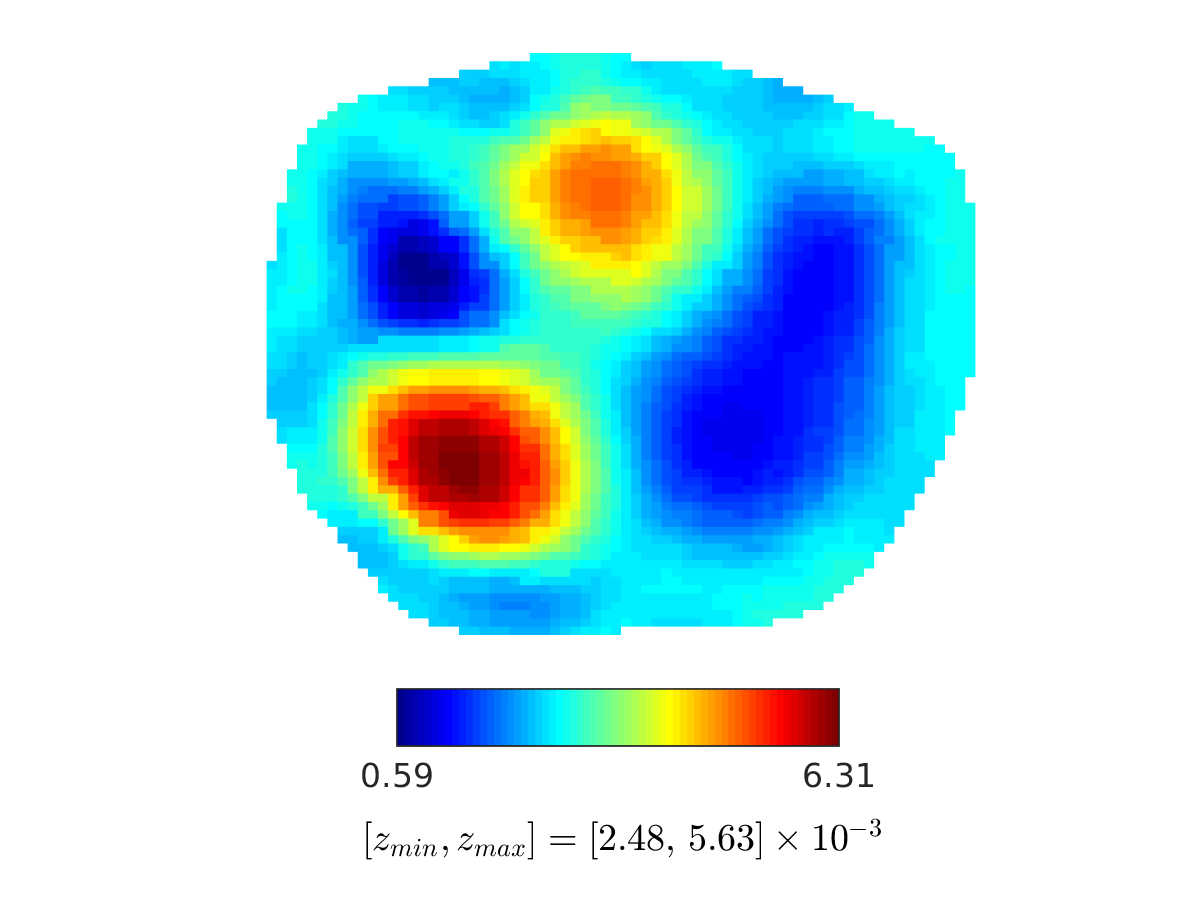}}	
	\end{picture}    
	\caption{\textbf{Case 1:} electrodes were displaced approximately $25\%$ of the true physical length. Top row shows the results of experimental data 1, images in middle row correspond to experimental data 2 and bottom row shows the results of experimental data 3. Left column: the 3 different experimental setups. Middle column: traditional reconstructions of isotropic conductivity in $\Omega$ using incorrect electrode locations. Right column: reconstructions with the proposed method using $\Omega_m = \Omega$ but incorrect electrode locations. The displayed quantity is  $\gamma_c(y) = \gamma(\bar{G}^{-1}(y))$, for $y \in \Omega_c = M (F_i(\Omega_m))$. The minimum and maximum value of the estimated contact impedance $\widehat{z}$ are shown below each figure.}
	\label{fig:case1_eldis_35}
\end{figure}

The results for the Case study 1 are shown in Figure~\ref{fig:case1_eldis_35}. Top row shows the results of experimental data 1, images in middle row correspond to experimental data 2 and bottom row shows the results of experimental data 3. 
The effect of incorrect modelling of electrode locations is clearly evident when using the conventional reconstruction method. Spurious details, mostly close to the boundary of the model domain, deteriorate the quality of the reconstruction.  
On the other hand, the proposed method produces reasonably good reconstructions even if the 
electrode locations are imperfectly known. The relative error~\eqref{eq:relative-error} for the recovery domain was $E(\Omega_c) = 1.79\%$, $E(\Omega_c) = 1.19\%$ and $E(\Omega_c) = 2.16\%$, for experimental data 1, 2 and 3 respectively.

\begin{figure}
	\centering
	\setlength{\unitlength}{1cm}
	\begin{picture}(20,13.2) 	
	\put(-0.8,9.7){\includegraphics[width=4.0cm,height=3.2cm]{Fig_ct-1.jpg}}
	\put(3.2,8.7){\includegraphics[width=6.1cm,height=4.5cm]{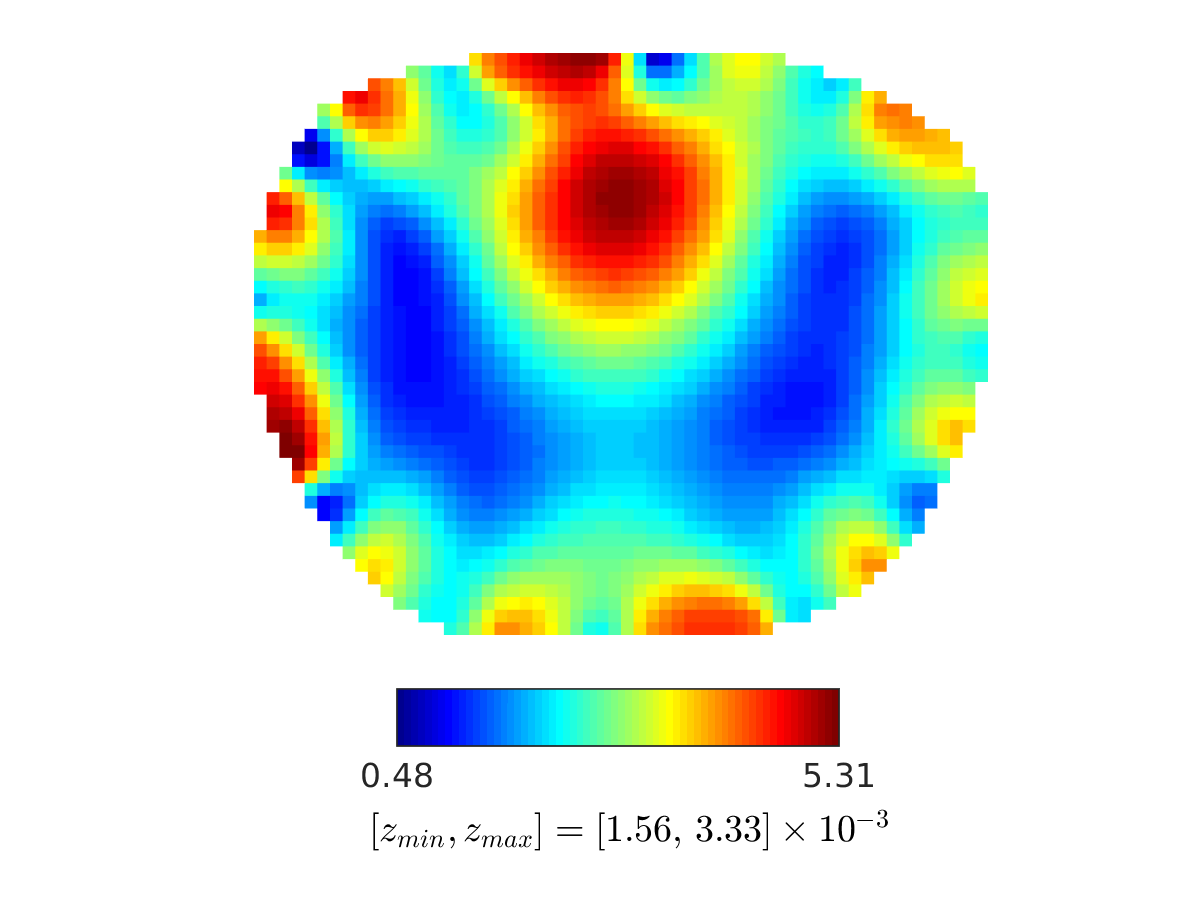}}
	\put(8.2,8.7){\includegraphics[width=6.1cm,height=4.5cm]{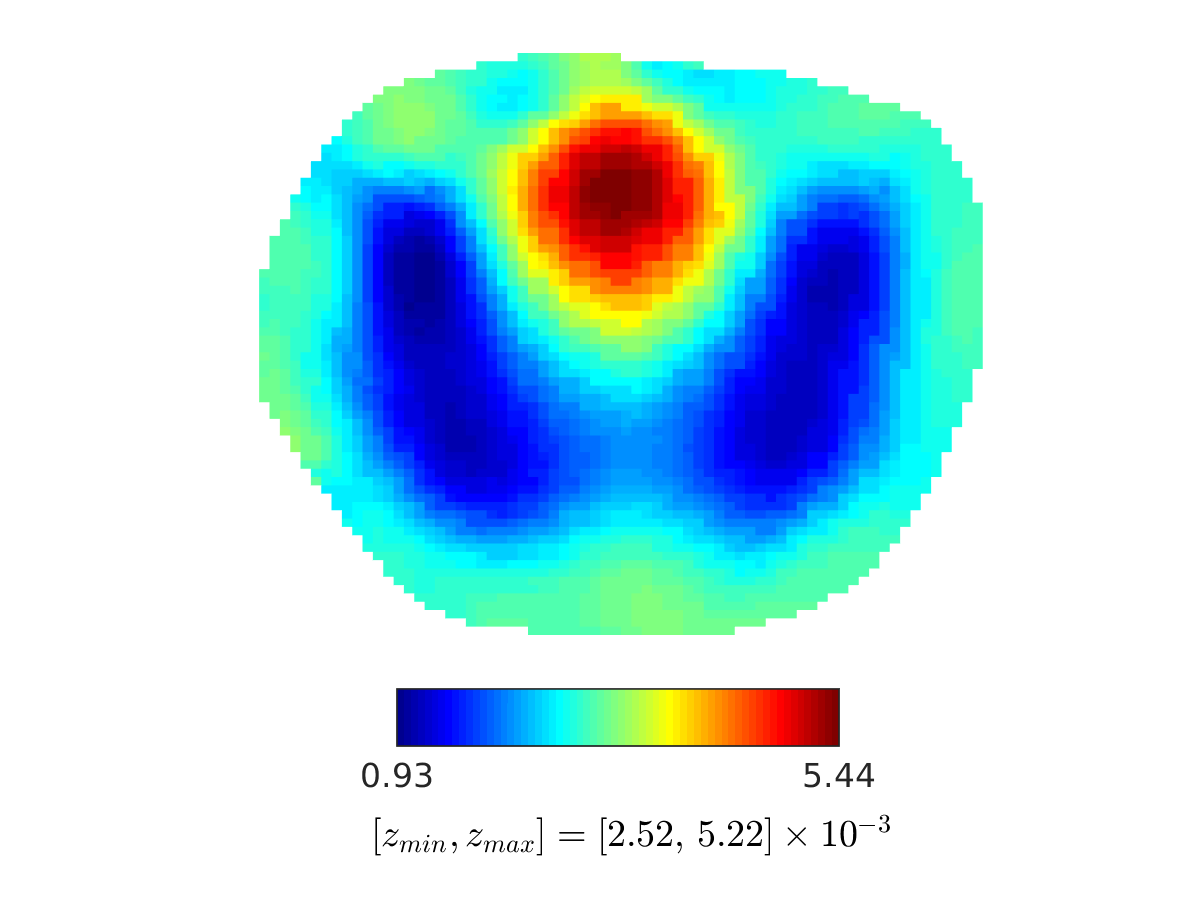}}
	\put(-0.8,5.2){\includegraphics[width=4.0cm,height=3.2cm]{Fig_ct-2.jpg}}
	\put(3.2,4.2){\includegraphics[width=6.1cm,height=4.5cm]{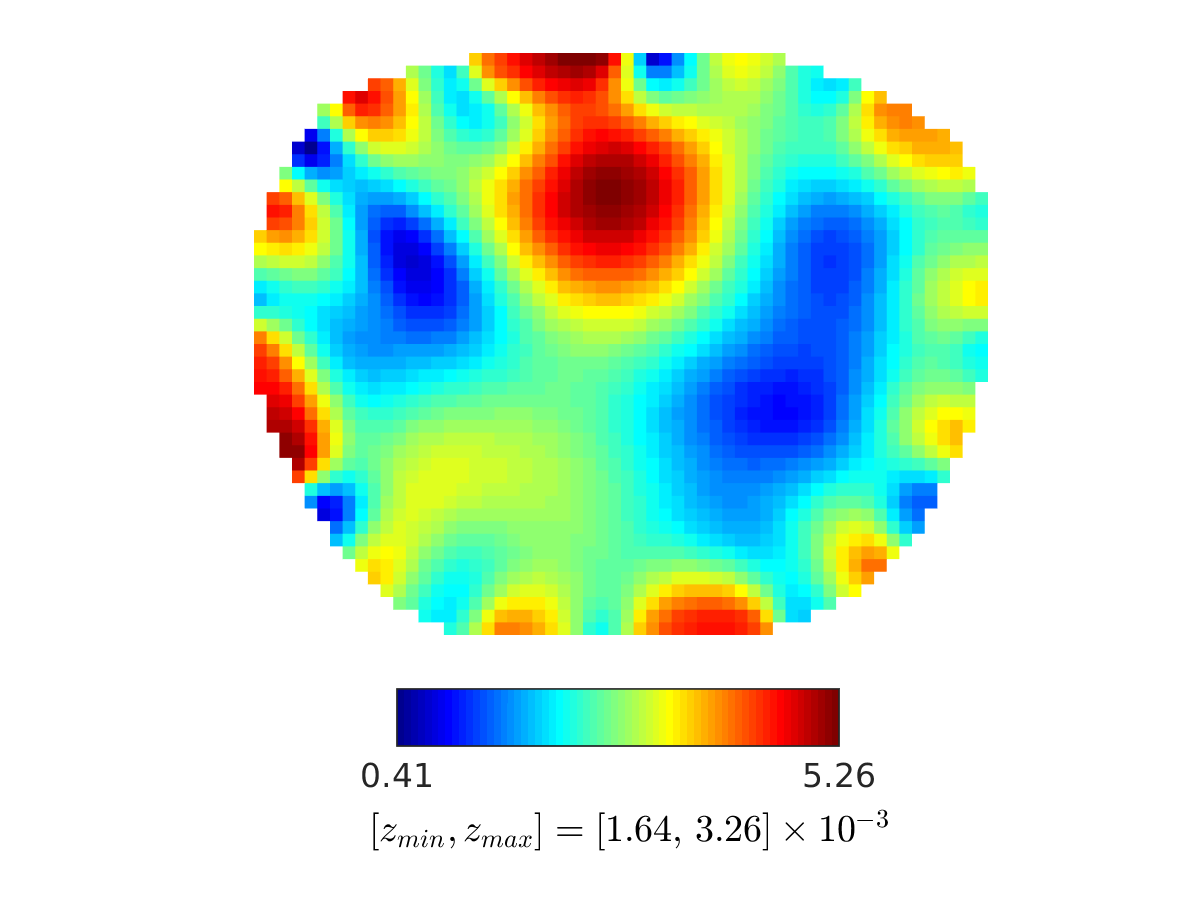}}
	\put(8.2,4.2){\includegraphics[width=6.1cm,height=4.5cm]{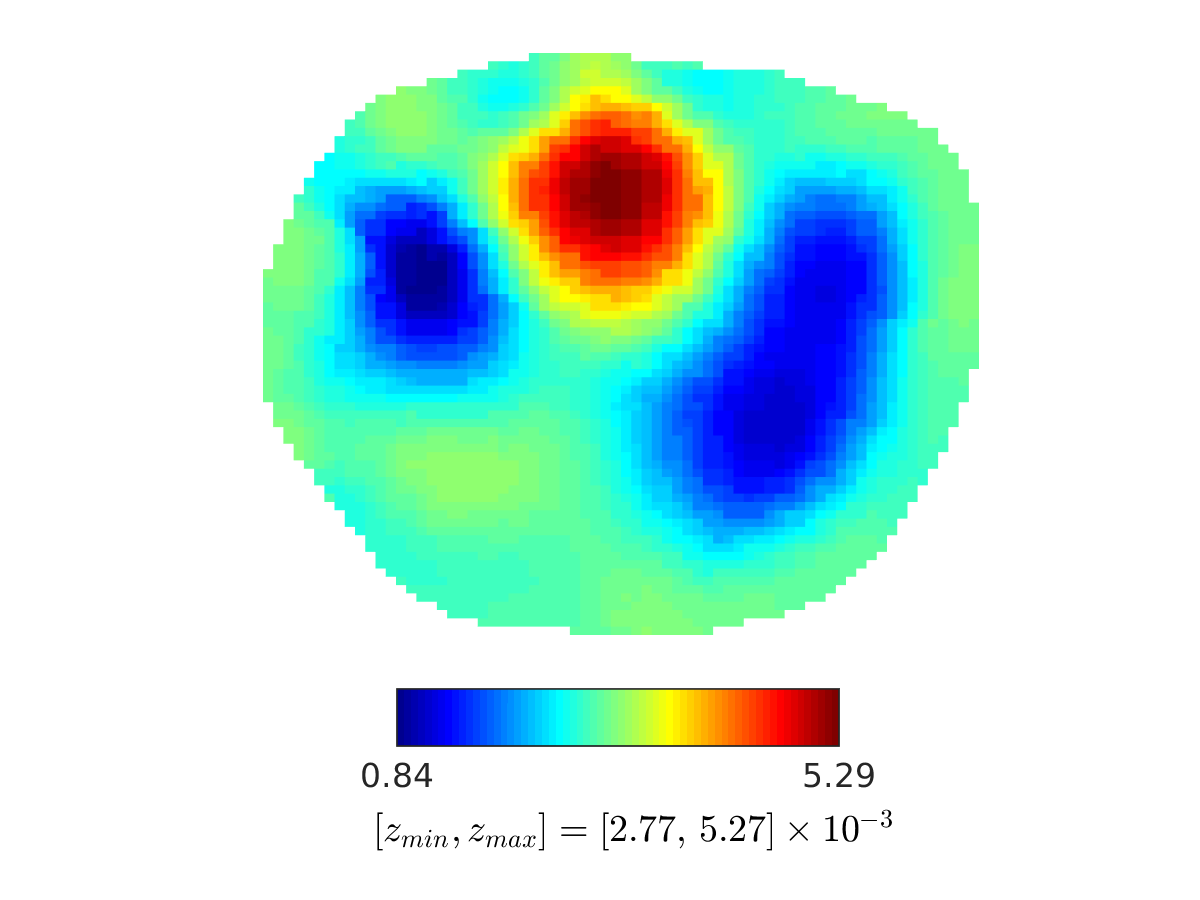}}
	\put(-0.8,0.7){\includegraphics[width=4.0cm,height=3.2cm]{Fig_ct-3.jpg}}
	\put(3.2,-0.3){\includegraphics[width=6.1cm,height=4.5cm]{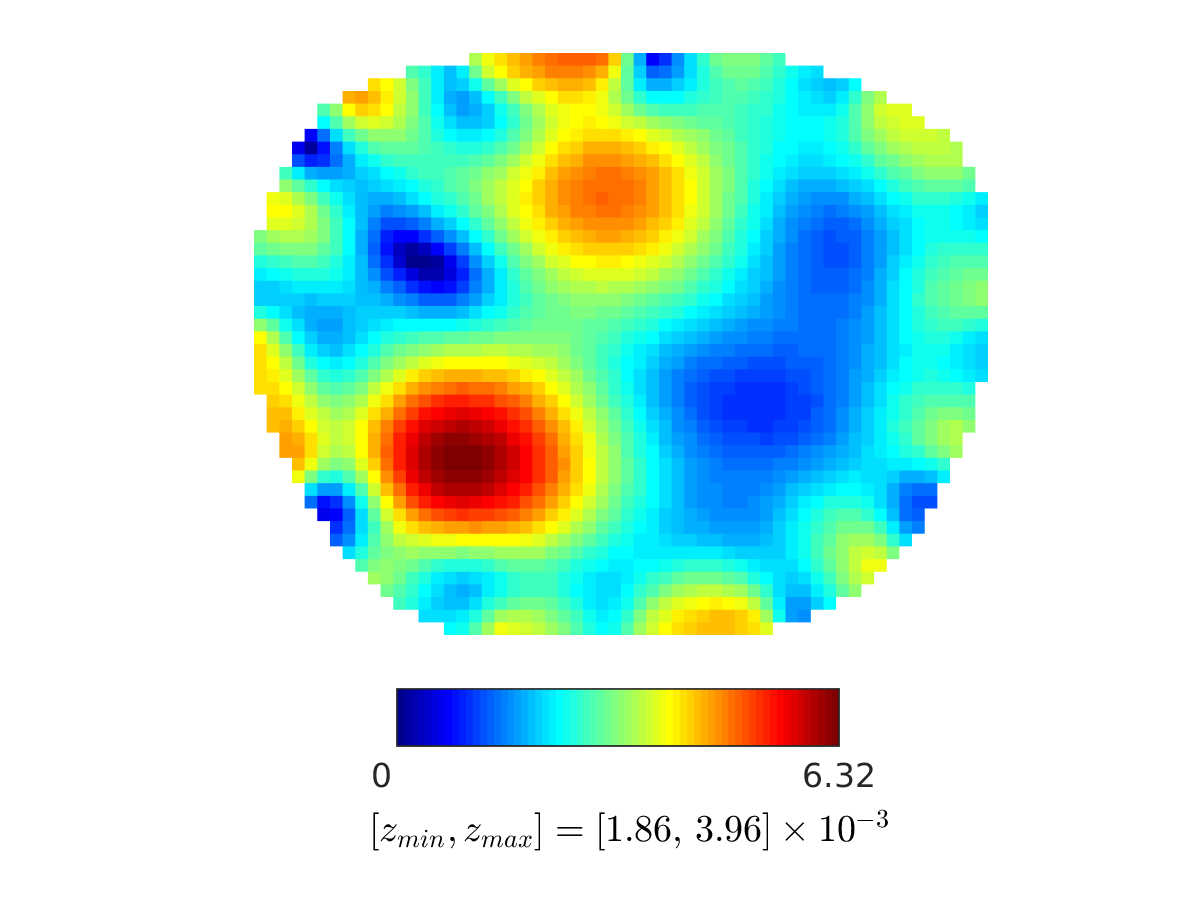}}
	\put(8.2,-0.3){\includegraphics[width=6.1cm,height=4.5cm]{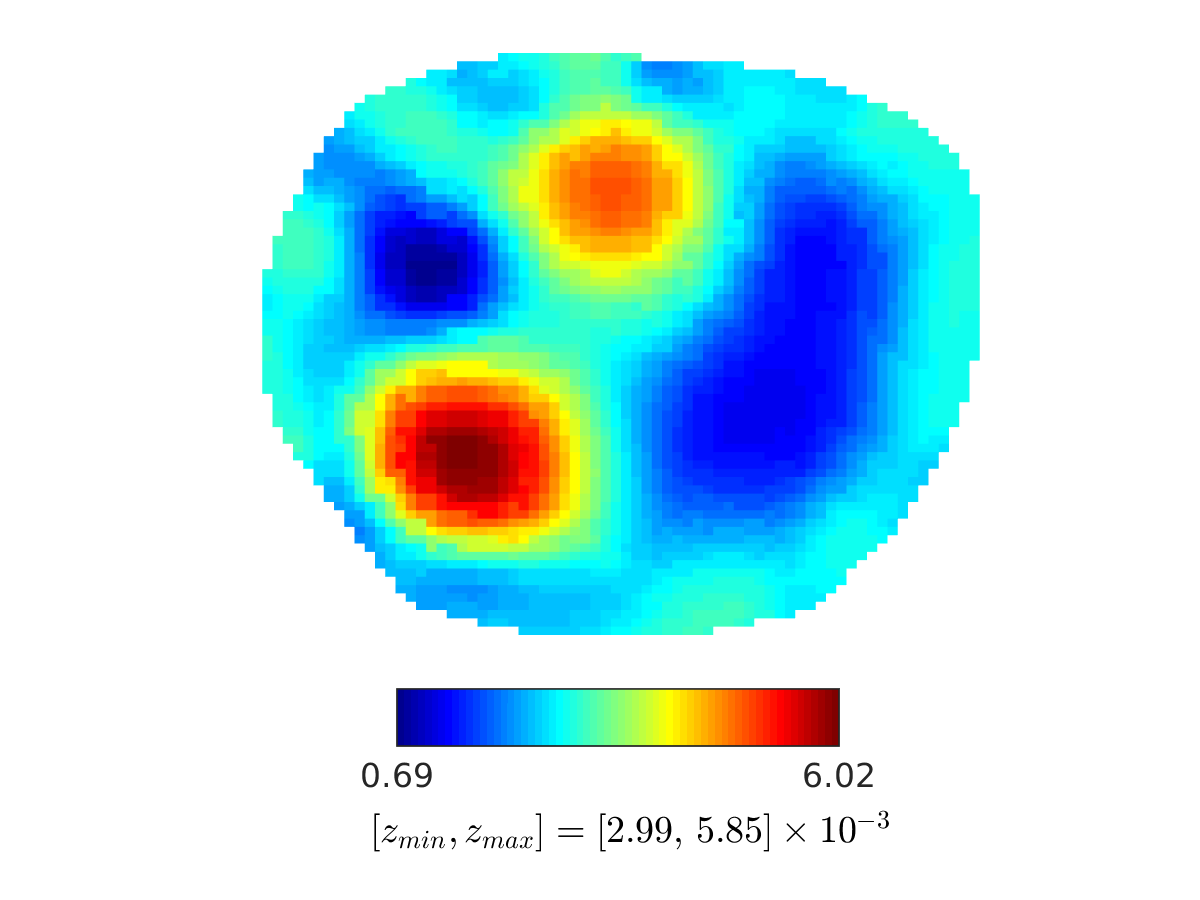}}	
	\end{picture}    
	\caption{\textbf{Case 2:} electrodes were displaced approximately $35\%$ of the true physical length. Top row shows the results of experimental data 1, images in middle row correspond to experimental data 2 and bottom row shows the results of experimental data 3. Left column: the 3 different experimental setups. Middle column: traditional reconstructions of isotropic conductivity in $\Omega$ using incorrect electrode locations. Right column: reconstructions with the proposed method using $\Omega_m = \Omega$ but incorrect electrode locations. The displayed quantity is  $\gamma_c(y) = \gamma(\bar{G}^{-1}(y))$, for $y \in \Omega_c = M (F_i(\Omega_m))$. The minimum and maximum value of the estimated contact impedance $\widehat{z}$ are shown below each figure.}
	\label{fig:case2_eldis_50}
\end{figure}

The corresponding results for the Case study 2 are shown in Figure~\ref{fig:case2_eldis_50}. The relative error~\eqref{eq:relative-error} for the recovery domain was $E(\Omega_c) = 2.79\%$ for experimental data 1, $E(\Omega_c) = 2.26\%$ for experimental data 2 and $E(\Omega_c) = 2.57\%$ for experimental data 3.   

\subsection{Checking the conformal maps $\bar{G}$ }\label{subsec:mapG}

Given a vector of $V \in \R^N$ of voltage measurement corresponding to a conductivity $\gamma$ in $\Omega$, the proposed method reconstructs an approximate conductivity
\[ \gamma_c(y)  = \gamma( \bar{G}^{-1} (y) ),\quad \text{for}\; y \in \Omega_c,\]
where $\bar{G} = M \circ F_i \circ F_e : \Omega \to \Omega_c$ is a conformal map.  That is, $\gamma_c$ is a conformally deformed image of the true conductivity $\gamma$.

Then, given a vectors of $V_1, V_2 \in \R^N$ of voltage measurements corresponding to conductivities $\gamma_1$ and $\gamma_2$ respectively, we have the conformal maps $\bar{G}_1: \Omega \to \Omega_{c,1}$ and $\bar{G}_2: \Omega \to \Omega_{c,2}$. Theory predicts that the map $\bar{G}_1 \circ \bar{G}_2^{-1}: \Omega_{c,2} \to \Omega_{c,1}$ should be the identity map. However, due to numerical errors we could only expect that this map is close to the identity map. 

\begin{figure}[ht!]
	\centering
	\setlength{\unitlength}{1cm}
	\begin{picture}(20,9.7) 
	\put(-0.7,4.7){\includegraphics[width=7.0cm,height=5.0cm]{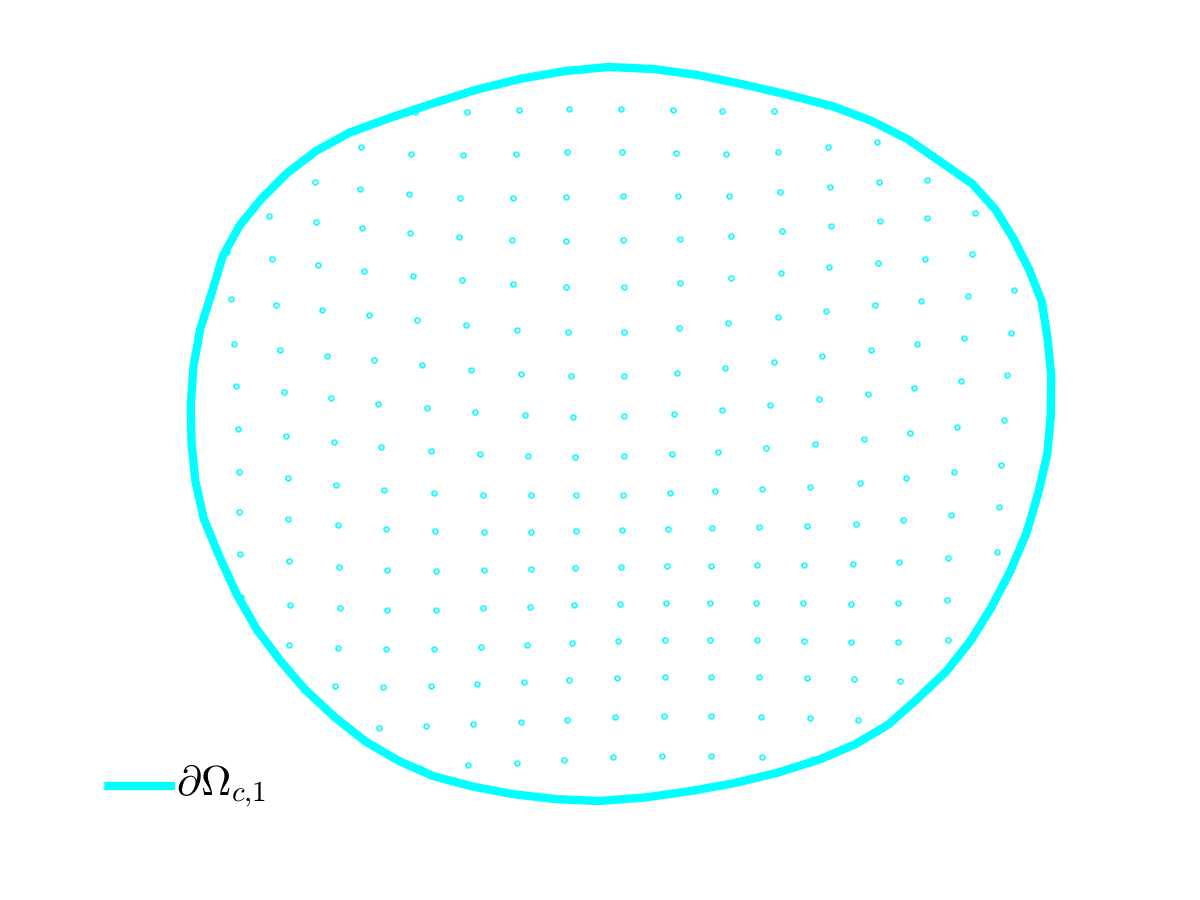}}
	\put(5.9,4.7){\includegraphics[width=7.0cm,height=5.0cm]{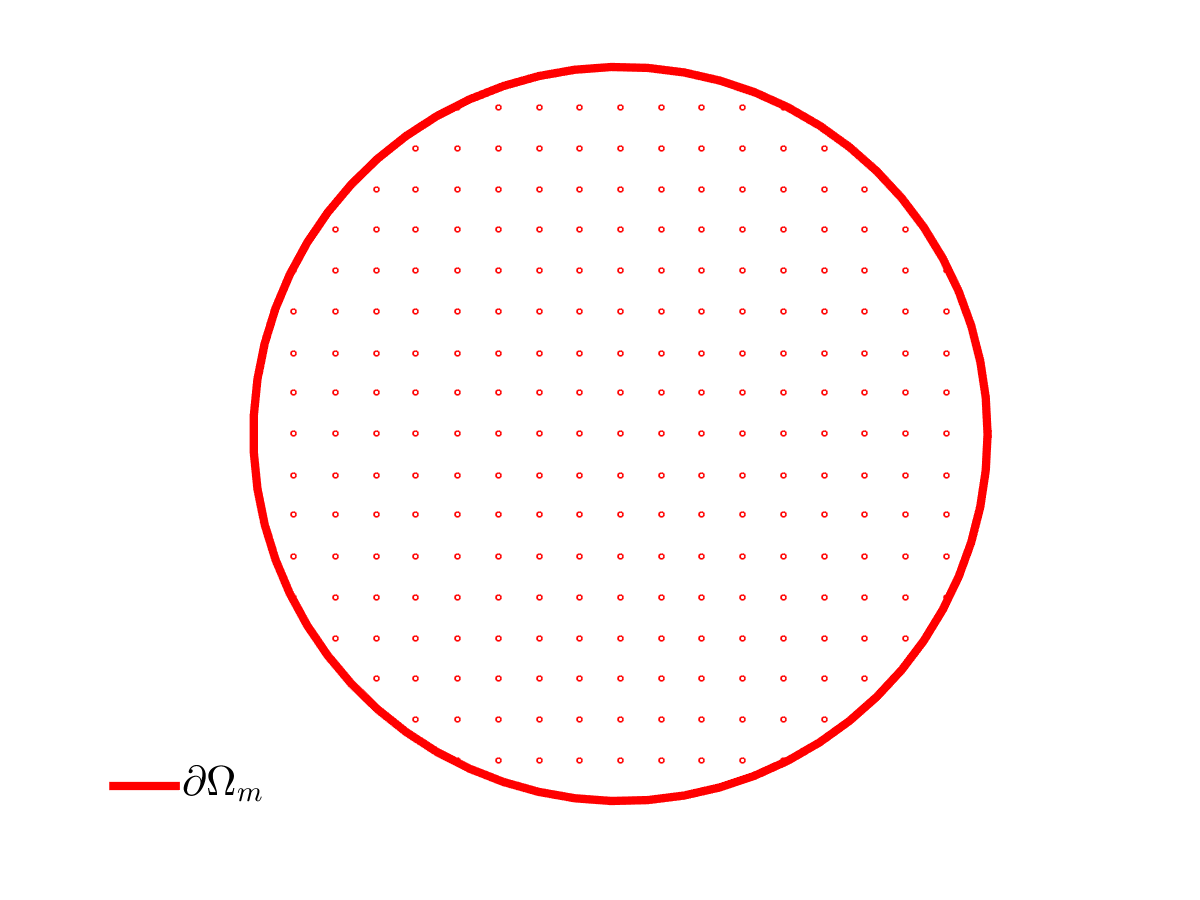}}
	\put(-0.7,-0.2 ){\includegraphics[width=7.0cm,height=5.0cm]{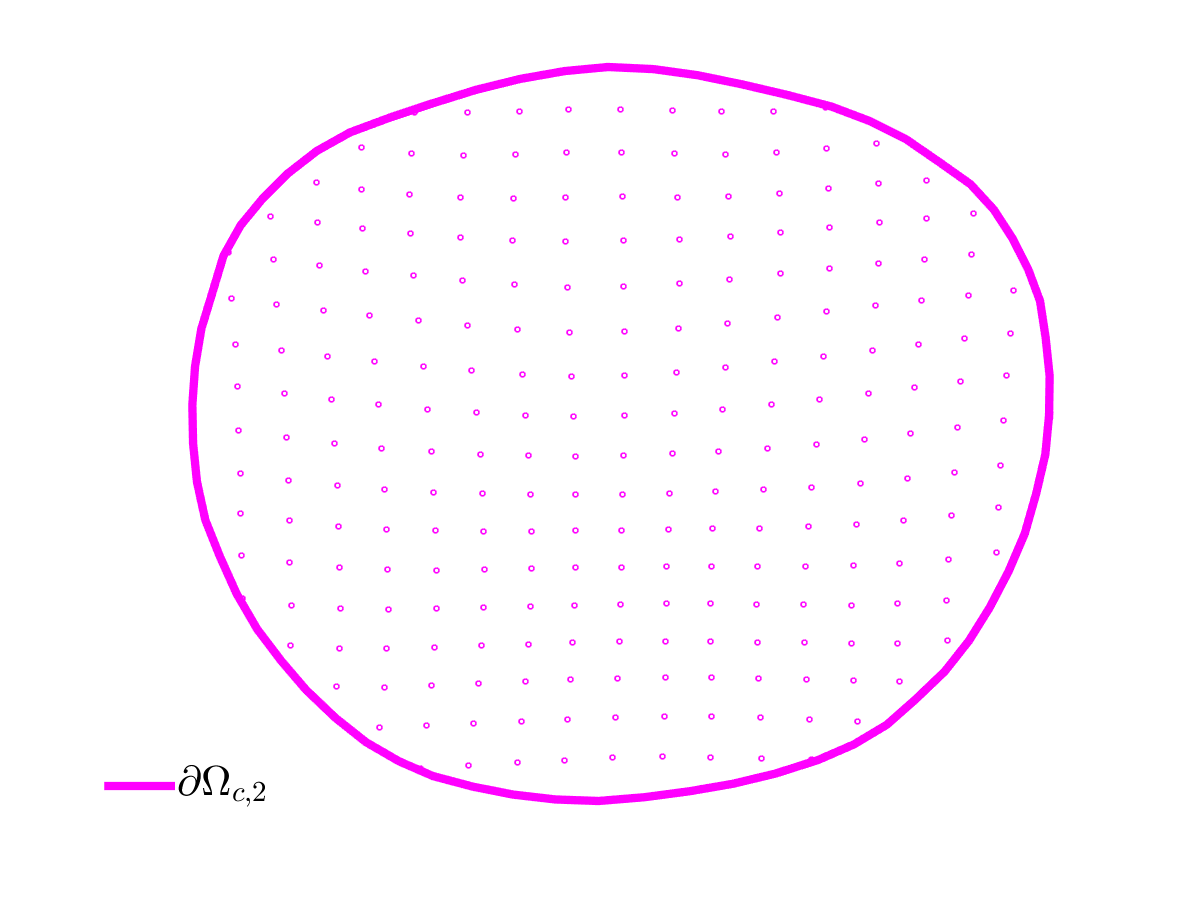}}
	\put(6.0,-0.2){\includegraphics[width=6.9cm,height=5.1cm]{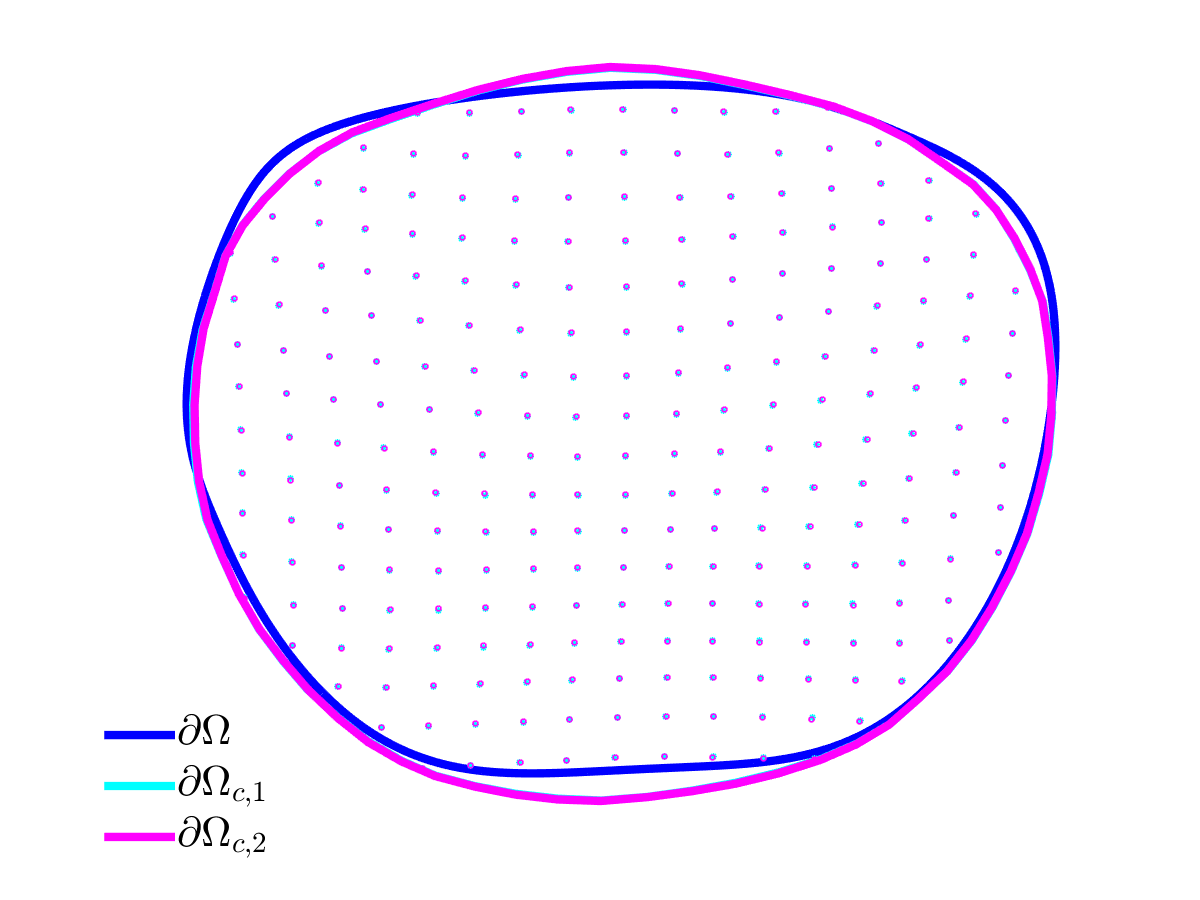}}
	\end{picture}    
	\caption{Map $\bar{G}_2 \circ \bar{G}_1^{-1}$. 
		Top left: reconstructed domain $\Omega_{c,1} = \bar{G}_1(\Omega)$. Top right: model domain $\Omega_m = F_{i,1}^{-1} \circ M_{1}^{-1}(\Omega_{c,1})$, the cyan points in  $\Omega_{c,1}$ are mapped to red points in $\Omega_m$. 
		Bottom left: domain $\Omega_{c,2} = M_2 \circ F_{i,2} (\Omega_m) = \bar{G}_2 \circ \bar{G}_1^{-1}(\Omega_{c,1})$, the red points in $\Omega_m$ are mapped to magenta points in $\Omega_{c,2}$. Bottom right: boundaries of domains $\Omega_{c,j}$, $j=1,2$ and true domain $\Omega$, the cyan points in $\Omega_{c,1}$ are mapped to magenta points in $\Omega_{c,2}$.}
	\label{fig:chek_G_exp_1_2}
\end{figure}

\begin{figure}[ht!]
	\centering
	\setlength{\unitlength}{1cm}
	\begin{picture}(20,9.7) 
	\put(-0.7,4.7){\includegraphics[width=7.0cm,height=5.0cm]{Fig112-a}}
	\put(5.9,4.7){\includegraphics[width=7.0cm,height=5.0cm]{Fig112-b}}
	\put(-0.7,-0.2 ){\includegraphics[width=7.0cm,height=5.0cm]{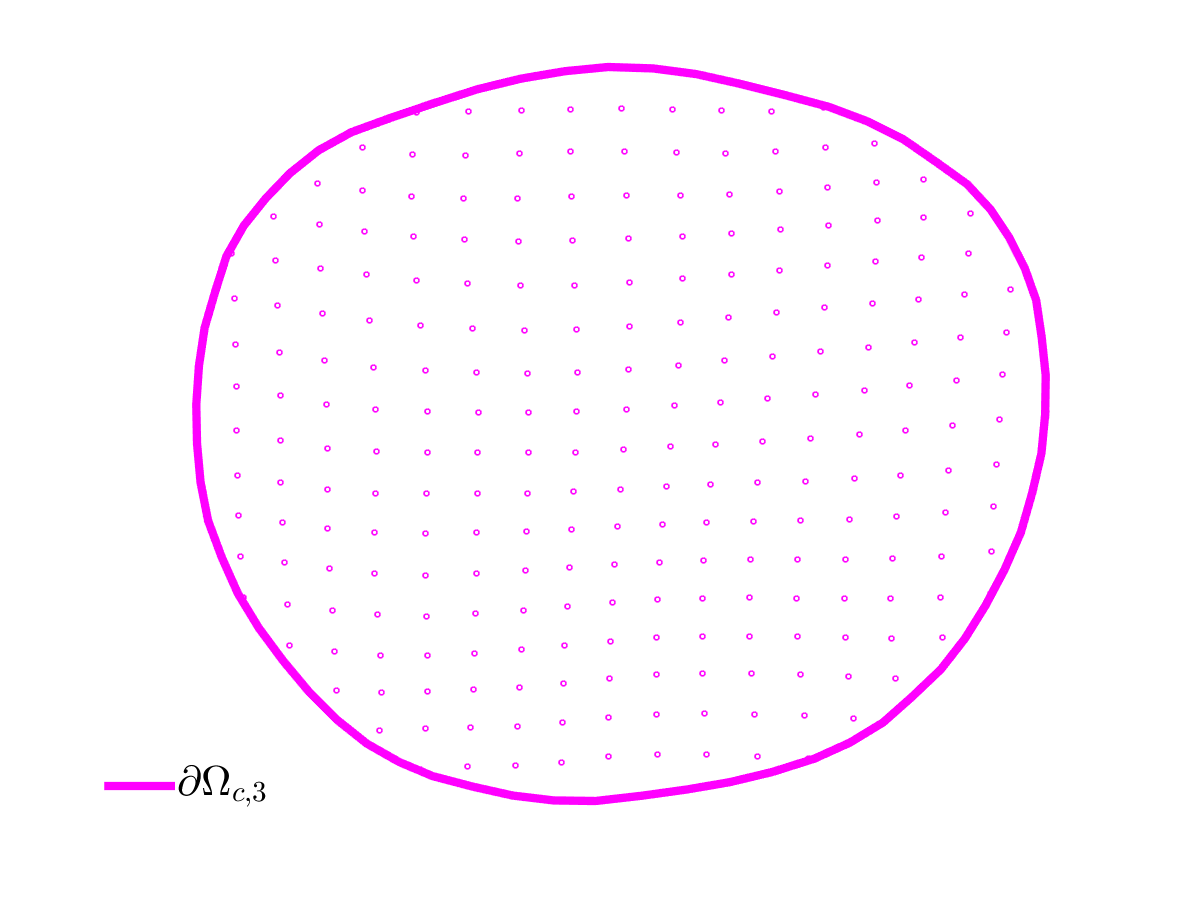}}
	\put(6.0,-0.2){\includegraphics[width=6.9cm,height=5.1cm]{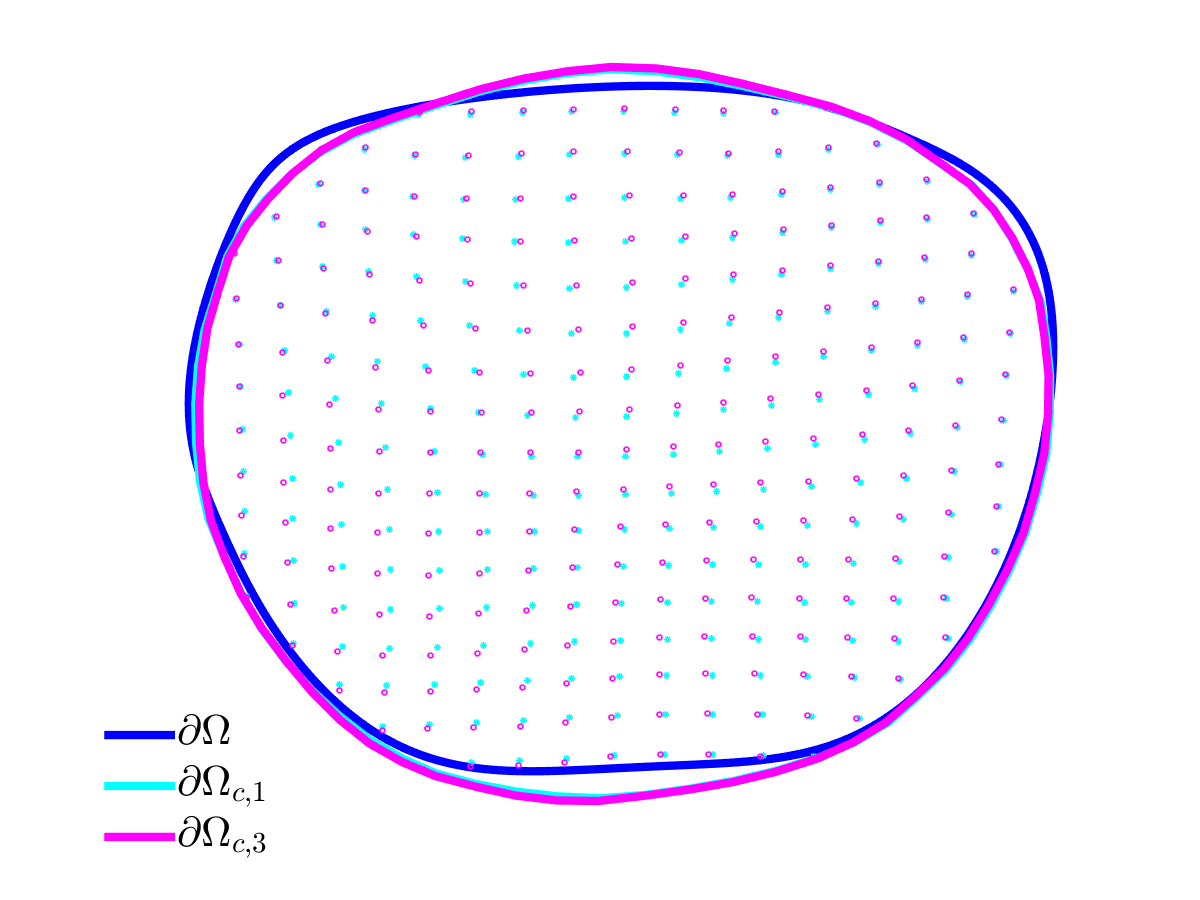}}
	\end{picture}    
	\caption{Map $\bar{G}_3 \circ \bar{G}_1^{-1}$. 
		Top left: reconstructed domain $\Omega_{c,1} = \bar{G}_1(\Omega)$. Top right: model domain $\Omega_m = F_{i,1}^{-1} \circ M_{1}^{-1}(\Omega_{c,1})$, the cyan points in  $\Omega_{c,1}$ are mapped to red points in $\Omega_m$. 
		Bottom left: domain $\Omega_{c,3} = M_3 \circ F_{i,3} (\Omega_m) = \bar{G}_3 \circ \bar{G}_1^{-1}(\Omega_{c,1})$, the red points in $\Omega_m$ are mapped to magenta points in $\Omega_{c,3}$. Bottom right: boundaries of domains $\Omega_{c,j}$, $j=1,3$ and true domain $\Omega$, the cyan points in $\Omega_{c,1}$ are mapped to magenta points in $\Omega_{c,3}$.}
	\label{fig:chek_G_exp_1_3}
\end{figure}

Recall that 
\begin{align*}
\bar{G}_j \circ \bar{G}_k^{-1} 
&= M_j \circ F_{i,j} \circ F_{e} \circ F_{e}^{-1} \circ F_{i,k}^{-1}  \circ M_{k}^{-1} \\
&= M_j \circ F_{i,j} \circ F_{i,k}^{-1} \circ M_{k}^{-1} 
\end{align*}
where $F_{i,k}^{-1} \circ M_{k}^{-1}: \Omega_{c,k} \to \Omega_{m}$ and $M_j \circ F_{i,j}: \Omega_{m} \to \Omega_{c,j}$. 
Figures~\ref{fig:chek_G_exp_1_2} and~\ref{fig:chek_G_exp_1_3} show the action of the map $\bar{G}_j \circ \bar{G}_k^{-1}: \Omega_{c,k} \to \Omega_{c,j}$ for $j=2,3$ and $k=1$. 
Both figures correspond to the different experimental data sets presented in section~\ref{subsec:disc-model-domain}. Observe that the recovered boundaries $\partial \Omega_{c,k}$ and $\partial \Omega_{c,j}$ from noisy voltage measurements from the three different conductivities are highly similar and they can hardly be distinguished from each other in the figures.

\begin{figure}[ht!]
	\centering
	\setlength{\unitlength}{1cm}
	\begin{picture}(20,10.5) 
	\put(-1.0,5.3){\includegraphics[width=7.5cm,height=5.4cm]{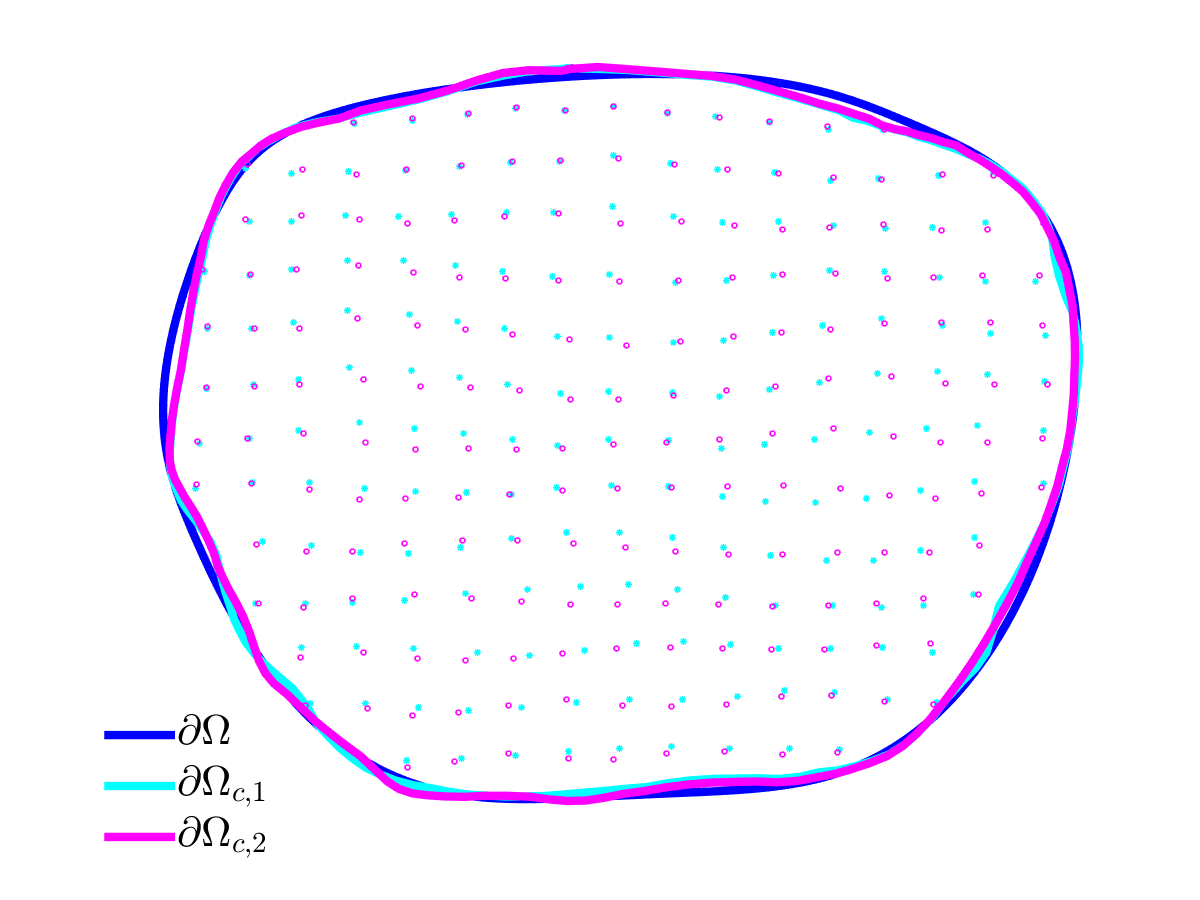}}
	\put(5.9,5.3){\includegraphics[width=7.5cm,height=5.4cm]{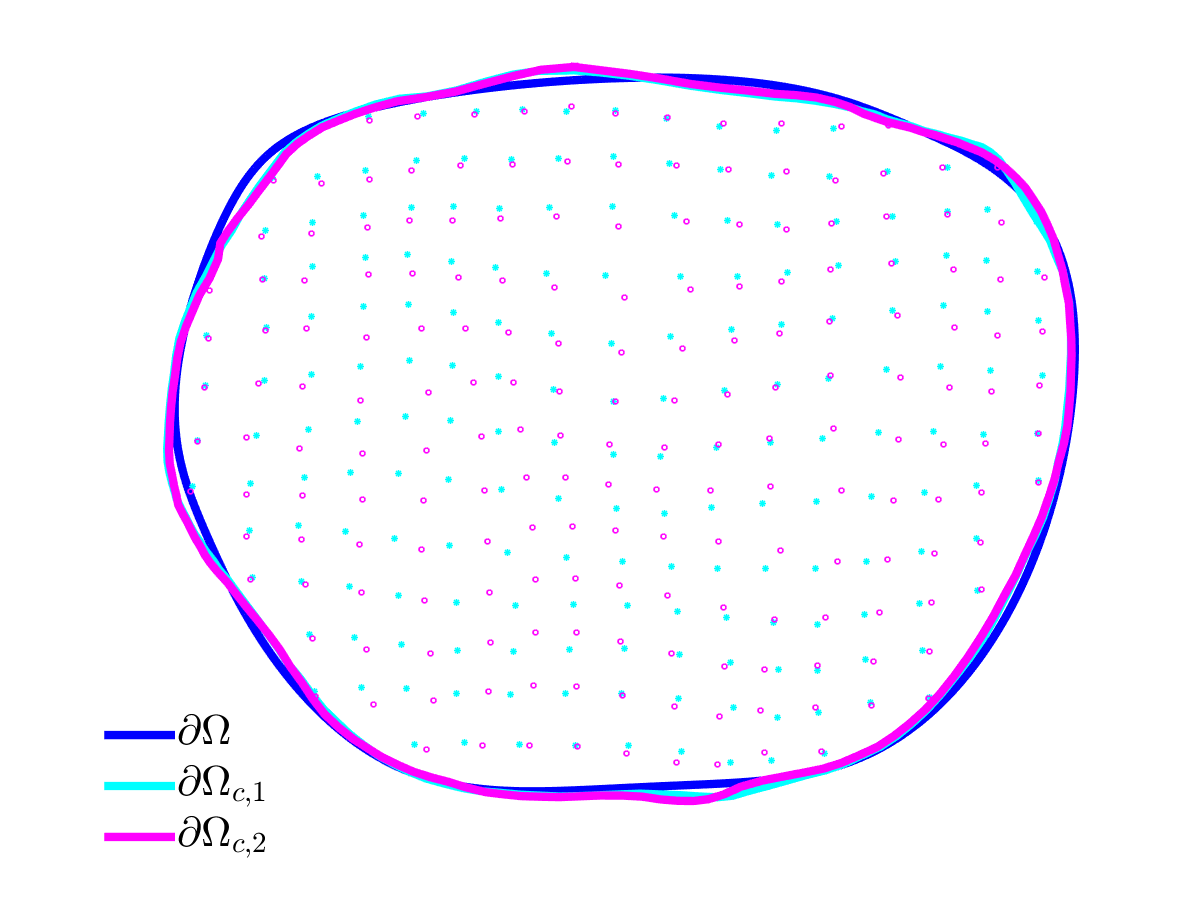}}
	\put(-1.0,-0.2 ){\includegraphics[width=7.5cm,height=5.4cm]{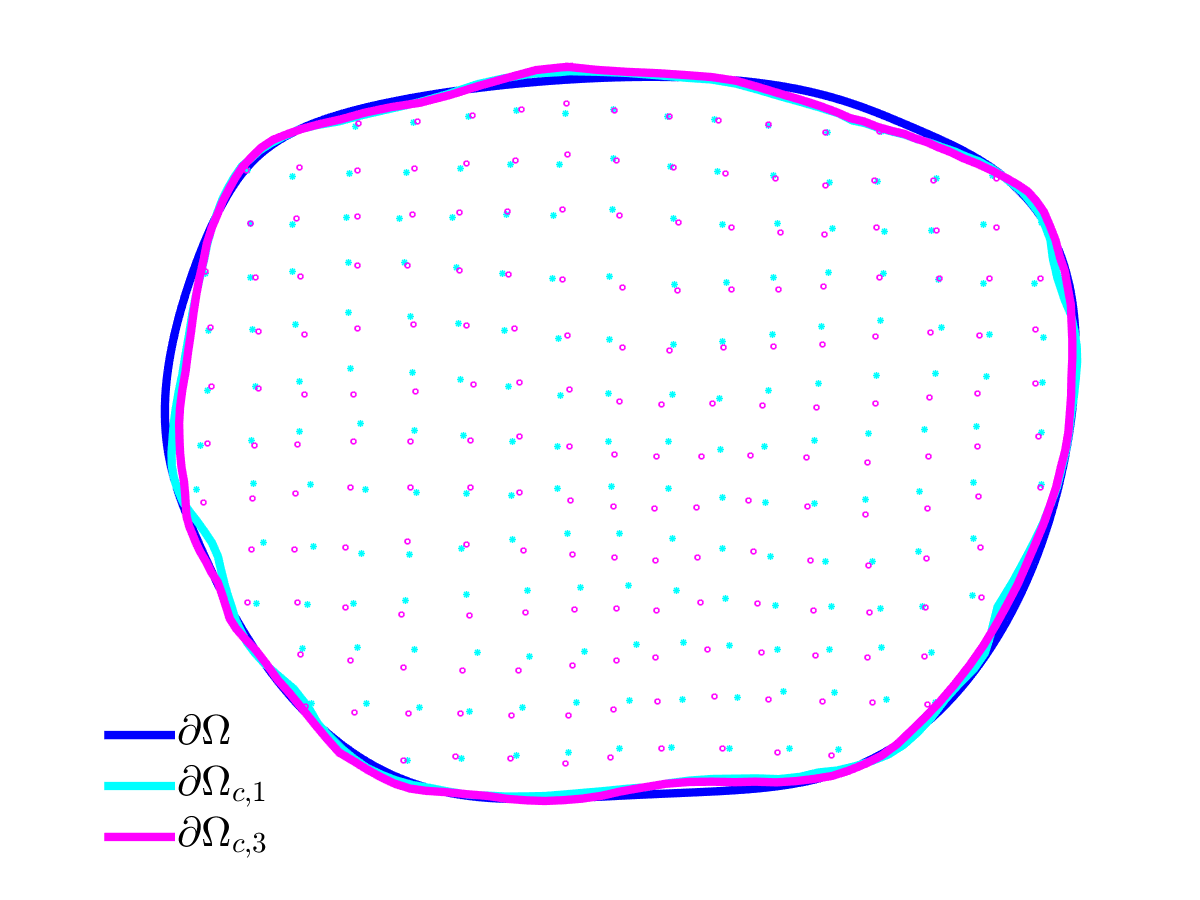}}
	\put(5.9,-0.2){\includegraphics[width=7.5cm,height=5.5cm]{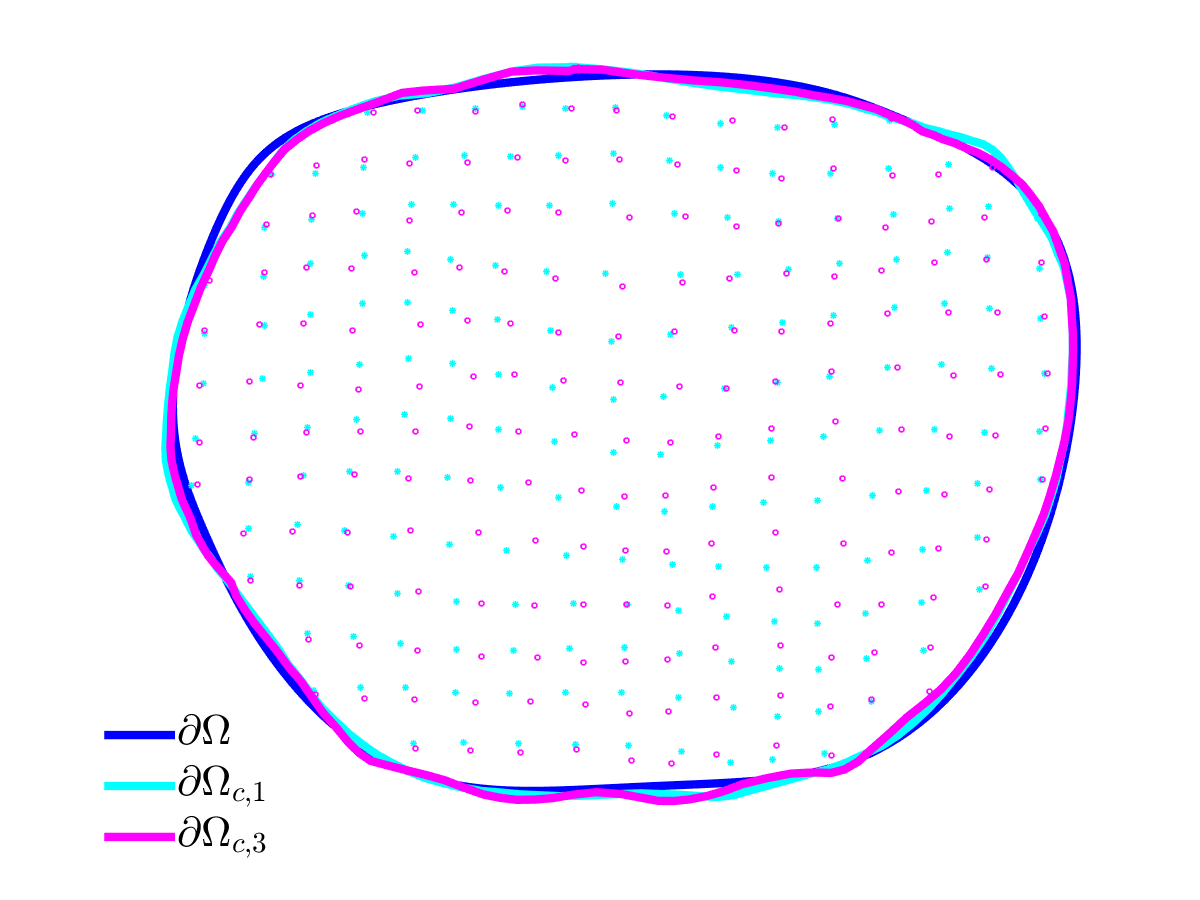}}
	\end{picture}    
	\caption{Maps $\bar{G}_k \circ \bar{G}_1^{-1}, k=2,3,$ for the situation where electrodes are displaced. Left column shows the results for the Case 1 ($25\%$ of displacement) and right column shows the results for Case 2 ($35\%$ of displacement). 
		Top row: boundaries of domains $\Omega_{c,1}$, $\Omega_{c,2} = \bar{G}_2 \circ \bar{G}_1^{-1}(\Omega_{c,1})$ and true domain $\Omega$.
		Bottom row: boundaries of domains $\Omega_{c,1}$, $\Omega_{c,3} = \bar{G}_3 \circ \bar{G}_1^{-1}(\Omega_{c,1})$ and true domain $\Omega$.
		In all figures the cyan points in $\Omega_{c,1}$ are mapped to magenta points in $\Omega_{c,k}$ by the map $\bar{G}_k \circ \bar{G}_1^{-1},$ for $k=2,3$. }
	\label{fig:chek_G_eldis}
\end{figure}

In Figure~\ref{fig:chek_G_eldis} we show the action of the map $\bar{G}_j \circ \bar{G}_k^{-1}: \Omega_{c,k} \to \Omega_{c,j}$ for $j=2,3$ and $k=1$,  
but corresponding to the different experimental data sets presented in section~\ref{subsec:eldis}. 

\subsection{Dynamic imaging in unknown domain}

In dynamic EIT the objective is to monitor changes in the conductivity in time or between measurements with different frequencies. In this section we consider the application of the proposed approach to dynamic imaging.

\subsubsection{ Two conductivities $\gamma_1$ and $\gamma_2$ defined in the same domain $\Omega$}

Given an isotropic conductivity $\gamma$ in $\Omega$, from Proposition~\ref{prop:anisotropic}, there is a unique map $F_e:\Omega \to \Omega_m$ such that $F_e \vert_{\partial \Omega} = f_m$ and 
\begin{equation}\label{eq:eta-det-1}
\eta(x) = \text{det}(\gamma_a(x))^{1/2} = \gamma( F_e^{-1}(x)), \qquad \text{for } \; x  \in \Omega_m,
\end{equation}
where $ \gamma_a $ is the unique solution of the minimization problem~\eqref{eq:constrained-minimization-1}. 
As we noted before, this result says that there exits a unique anisotropic conductivity $\gamma_a$ in $\Omega_m$ such that 
$ \text{det}(\gamma_a(x))^{1/2}$ 
gives a deformed image in the model domain of the original conductivity $\gamma$ defined in the true domain. The deformation $F_e$ of the conductivity image depends on the map $f_m$, 
but not on the original conductivity $\gamma$. This result is very useful because it implies that local perturbations of conductivity remain local in the reconstruction as we show below.  

Let $f_m : \partial \Omega \to \partial \Omega_m$ be a (fixed) boundary modeling map and let $\gamma_1$ and $\gamma_2 = \gamma_1  + \delta \gamma$ be two isotropic conductivities in $\Omega$. Then, by equation~\eqref{eq:eta-det-1} we have
\begin{equation}\label{eq:diff-imag-Om}
\eta_1(x) - \eta_2(x) = \delta \gamma( F_e^{-1}(x)), \qquad \text{for } \; x  \in \Omega_m.
\end{equation}
Note that the theory states that maps $F_e$,  $F_i$ and $M$ are the same for the two conductivities $\gamma_1$ and $\gamma_2$ when the model map $ f_m $ is the same. From the numerical point of view, we observed in section~\ref{subsec:mapG} that although these maps are not exactly the same, they are highly similar.

\begin{figure}
	\centering
	\setlength{\unitlength}{1cm}
	\begin{picture}(20,14.5) 	
	\put(-1.0,10.8){\includegraphics[width=3.9cm,height=3.1cm]{Fig_ct-1.jpg}}
	\put(3.2,9.5){\includegraphics[width=6.4cm,height=4.8cm]{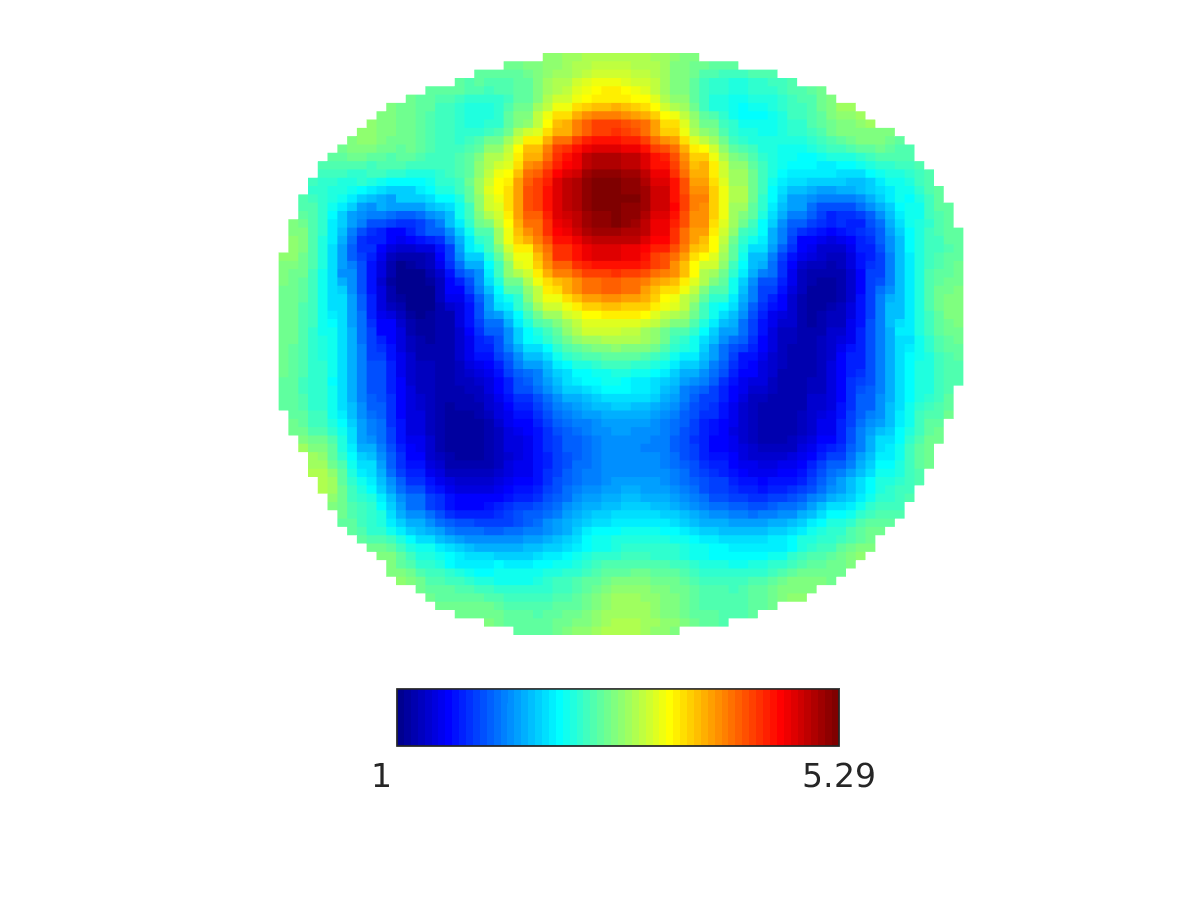}}
	\put(-1.0,5.8){\includegraphics[width=3.9cm,height=3.1cm]{Fig_ct-2.jpg}}
	\put(3.2,4.5){\includegraphics[width=6.4cm,height=4.8cm]{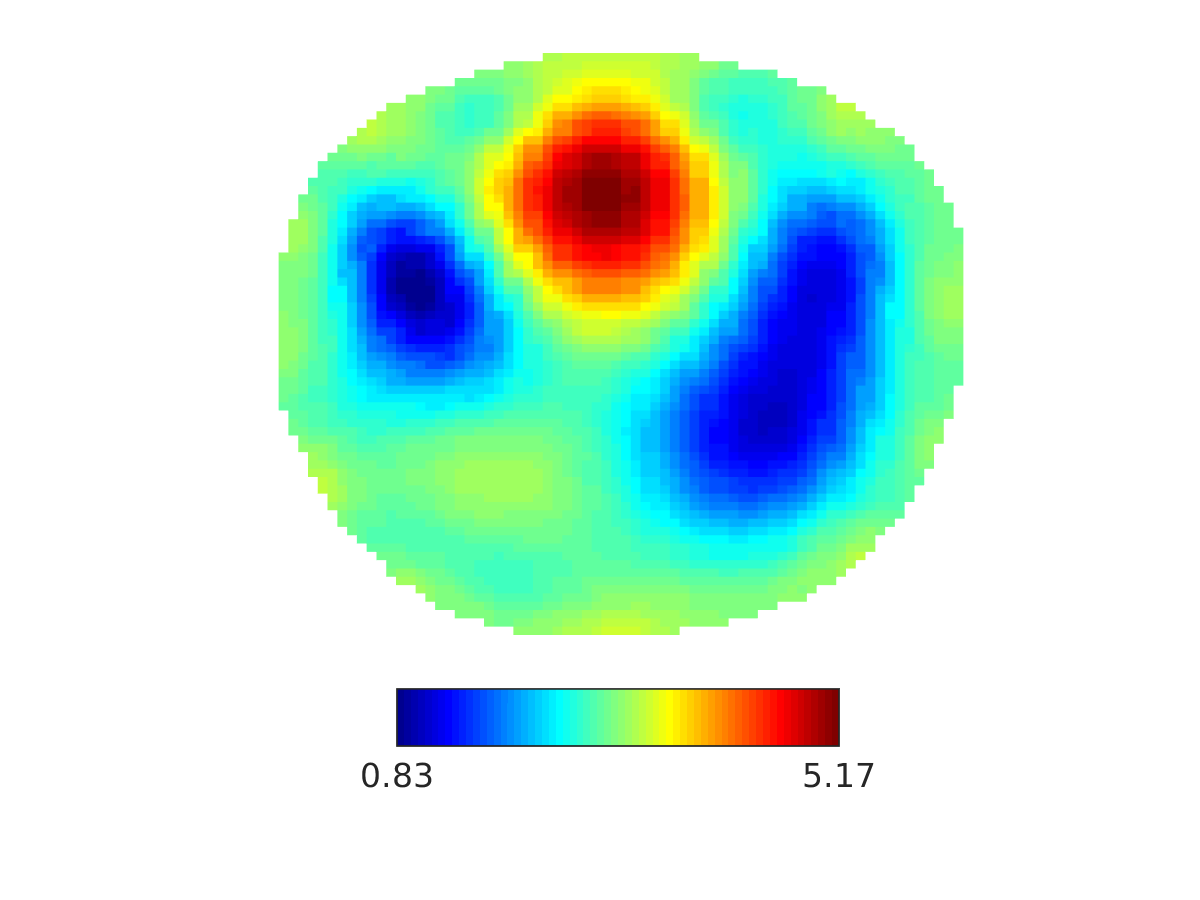}}
	\put(8.4,4.5){\includegraphics[width=6.4cm,height=4.8cm]{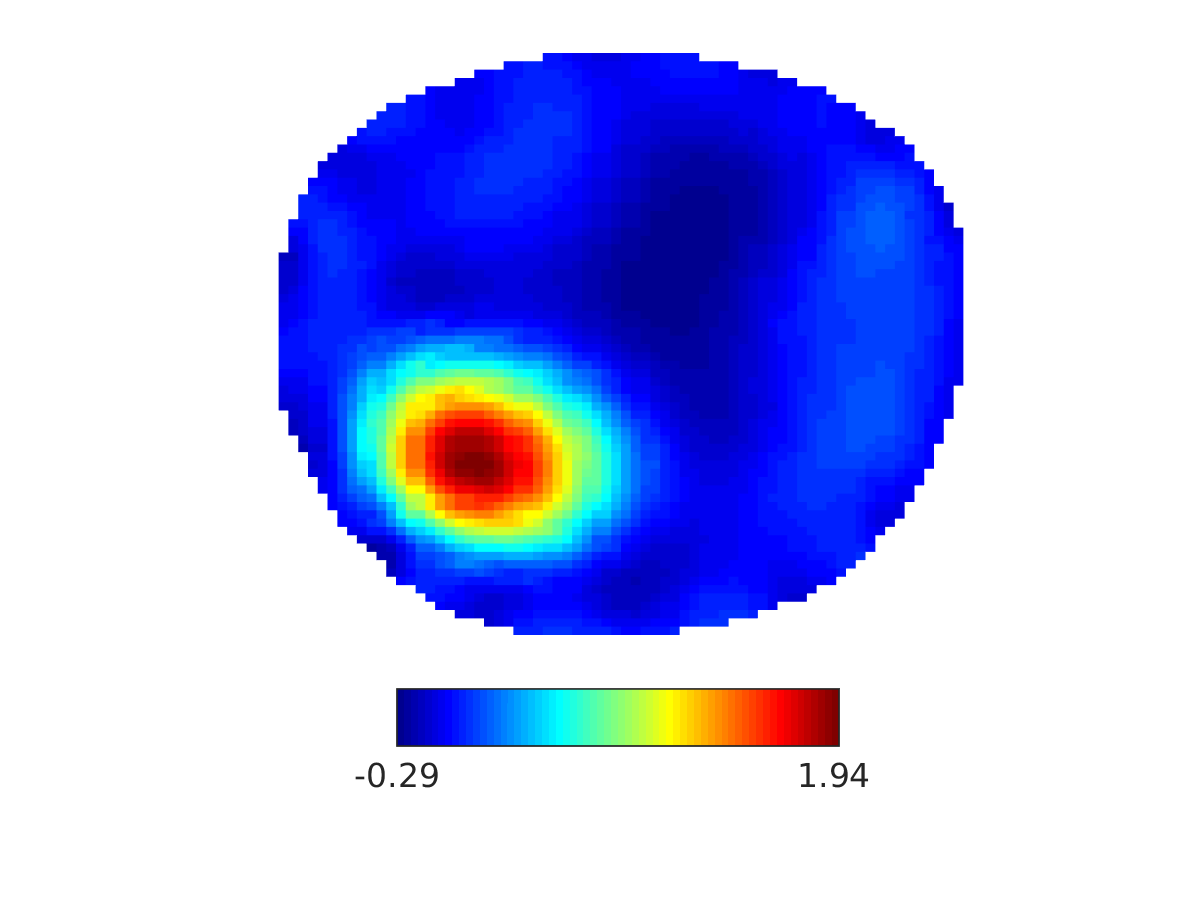}}
	\put(-1.0,0.8){\includegraphics[width=3.9cm,height=3.1cm]{Fig_ct-3.jpg}}
	\put(3.2,-0.5){\includegraphics[width=6.4cm,height=4.8cm]{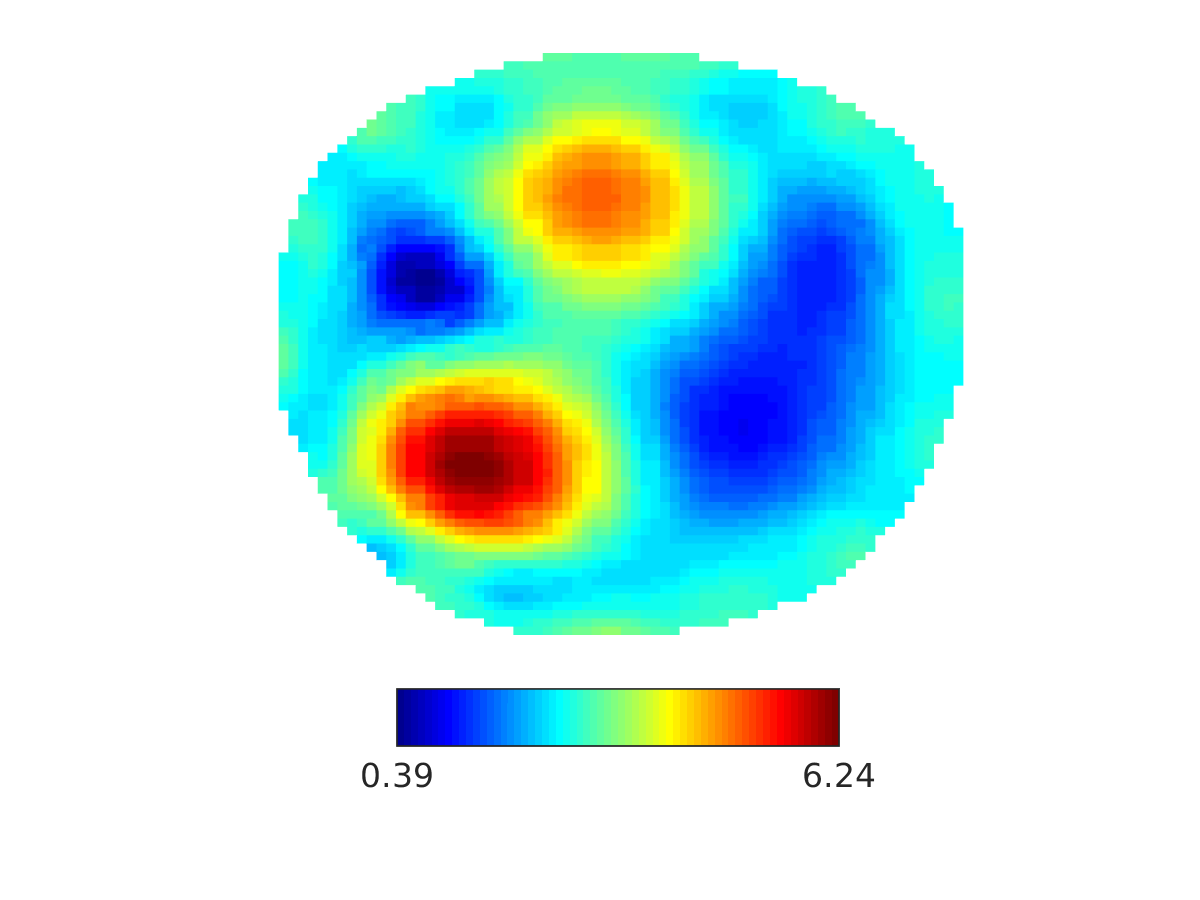}}
	\put(8.5,-0.5){\includegraphics[width=6.4cm,height=4.8cm]{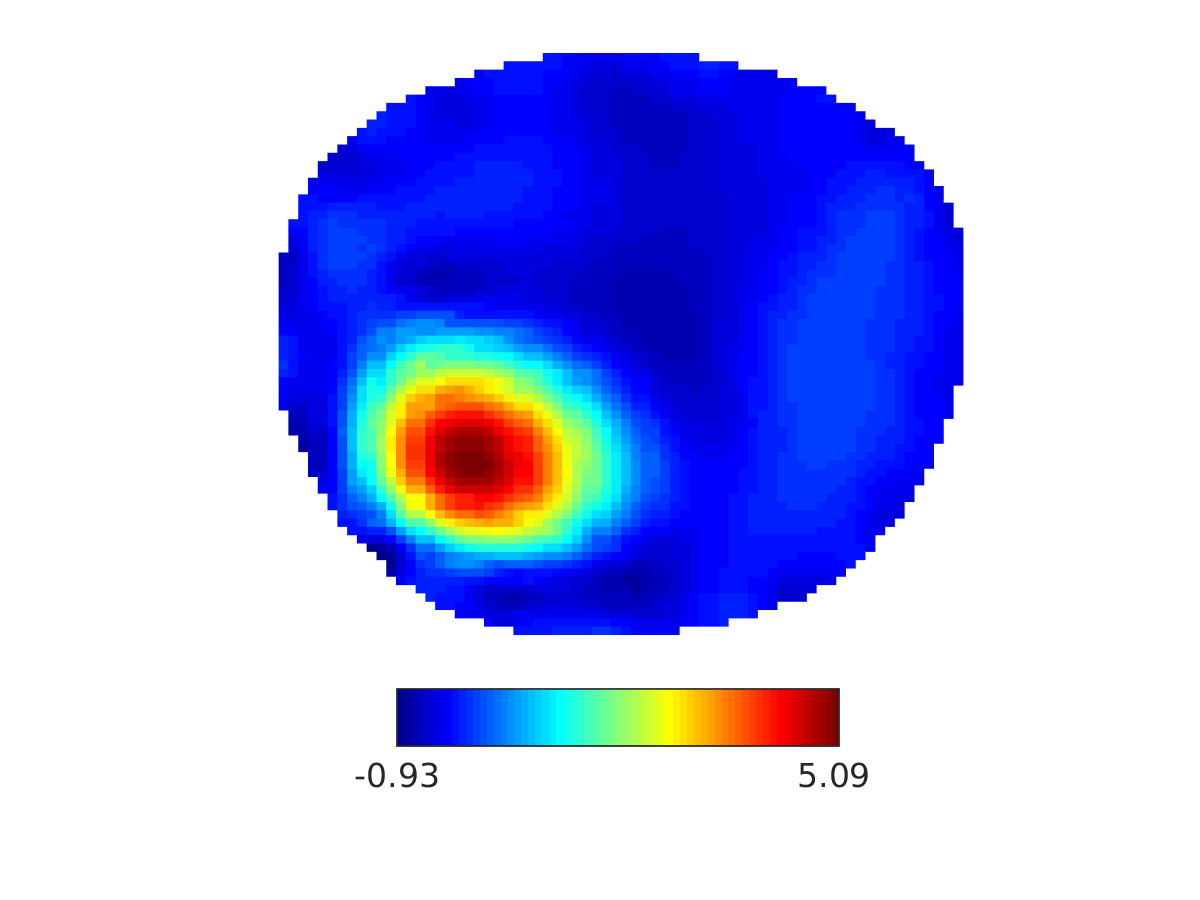}}
	\put(0.7,14.3){$\gamma_1$}
	\put(0.7,9.3){$\gamma_2$}
	\put(1.7,4.3){$\gamma_3$}
	\put(4.5,14.3){$\gamma_{c,1}(y)=\eta_1(F_{i,1}^{-1} \circ M_{1}^{-1}(y))$}
	\put(5.0,9.3){$\eta_2(F_{i,1}^{-1} \circ M_1^{-1}(y))$}
	\put(5.0,4.3){$\eta_3(F_{i,1}^{-1} \circ M_1^{-1}(y))$}
	\put(10.3,9.3){$ \left( \gamma_2 -  \gamma_1 \right) (\bar{G}_1^{-1}(y))$}
	\put(10.3,4.3){$ \left( \gamma_3 -  \gamma_1 \right) (\bar{G}_1^{-1}(y))$}
	\end{picture}    
	\caption{\textbf{Difference imaging: $\gamma_1$ vs $\gamma_k, \, k=2,3$ in $\Omega_{c,1}$.} 
		Top row shows the results of experimental data 1,  images in middle row correspond to experimental data 2 and bottom row shows the results of experimental data 3.
		Left column: the 3 different experimental setups.
		Middle column: functions $\eta_1 (F_{i,1}^{-1} \circ M_{1}^{-1}(y))$, $\eta_2 (F_{i,1}^{-1} \circ M_{1}^{-1}(y))$ and   $\eta_3 (F_{i,1}^{-1} \circ M_{1}^{-1}(y))$ for $y  \in \Omega_{c,1}$.  
		Right column: difference  $ \left( \eta_2 - \eta_1 \right) (F_{i,1}^{-1} \circ M_{1}^{-1}(y)) = \left( \gamma_2 -  \gamma_1 \right)(\bar{G}_1^{-1}(y))$ and  $ \left(\eta_3  - \eta_1 \right) (F_{i,1}^{-1} \circ M_{1}^{-1}(y)) = \left( \gamma_3 -  \gamma_1 \right) (\bar{G}_1^{-1}(y))$, respectively (see equation~\eqref{eq:diff-imag}).
	}
	\label{fig:diff_imag_exp_1_2_3_v1}
\end{figure}

Let $F_{i,1}$ be the solution of problem~\eqref{eq:beltrami-1}-\eqref{eq:beltrami-3} with coefficient $\mu$ associated to conductivity $\gamma_{a,1}$ and $M_{1}$ the Möbius map from $\Omega_{i,1} = F_{i,1}(\Omega_m)$ to $\Omega_{c,1}$, then
\begin{equation}\label{eq:diff-imag}
\eta_1(F_{i,1}^{-1} \circ M_{1}^{-1}(y) )- \eta_2(F_{i,1}^{-1} \circ M_{1}^{-1}(y)) = \delta \gamma( \bar{G}_1^{-1}(y)), \qquad \text{for } \; y  \in \Omega_{c,1}.
\end{equation}
That is, we can obtain a conformally deformed image of the difference between conductivities $\gamma_1$ and $\gamma_2$. Observe that in equation~\eqref{eq:diff-imag}, it is not necessary to compute $F_{i,2}$ and $M_2$ since we only need function $\eta_2$.
In order to obtain function $\eta_2$ we consider the minimization problem~\eqref{eq:min-problem-discrete} but fixing $z$, $\lambda$ and $\theta$ to be the same as with experimental data 1, that is $V_1 \in \R^N$, but use the experimental measurements $V_2 \in \R^N$ to compute $\eta_2$. In Figure~\ref{fig:diff_imag_exp_1_2_3_v1} we show the results corresponding to the different experimental data sets presented in section~\ref{subsec:disc-model-domain}.

\begin{figure}
	\centering
	\setlength{\unitlength}{1cm}
	\begin{picture}(20,14.5) 	
	\put(-1.0,10.6){\includegraphics[width=4.0cm,height=3.3cm]{Fig_ct-1.jpg}}
	\put(3.2,9.5){\includegraphics[width=6.4cm,height=4.8cm]{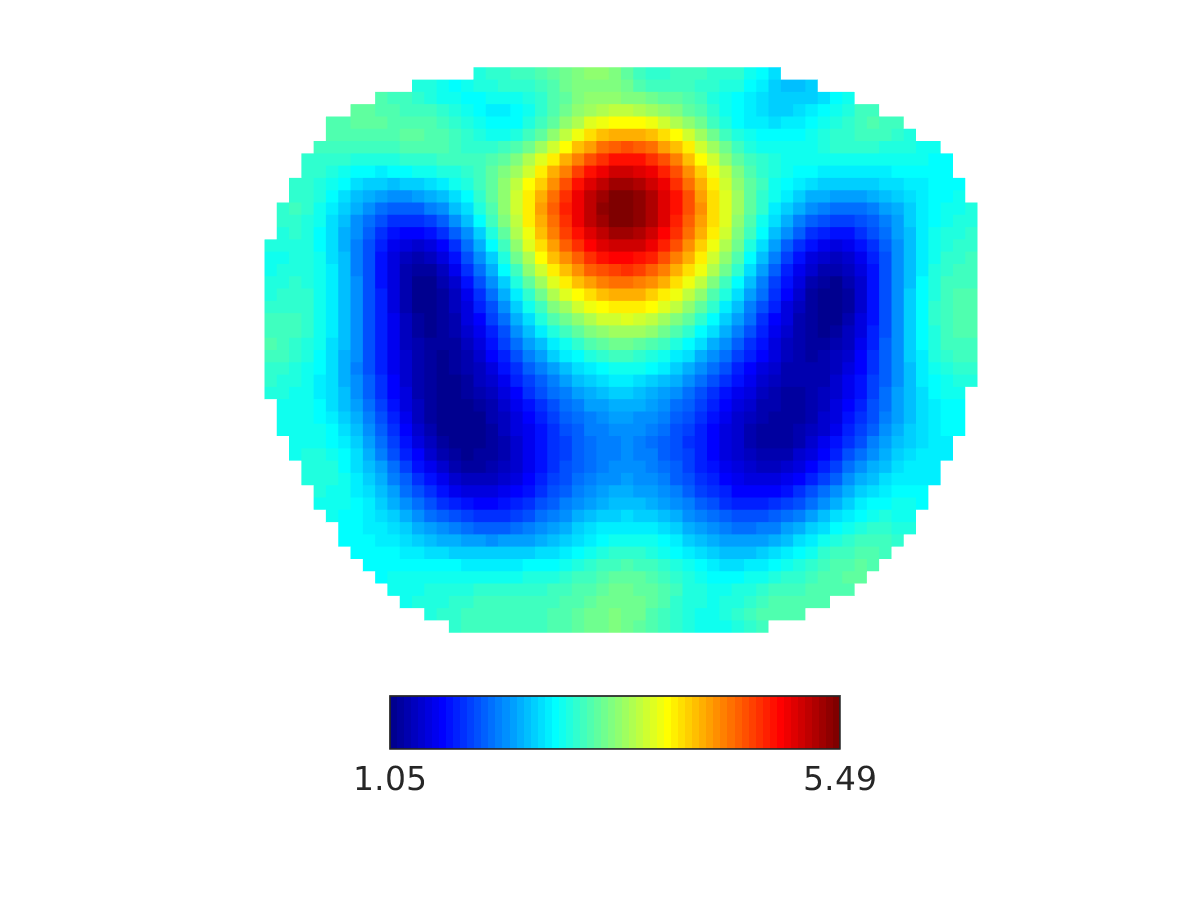}}
	\put(-1.0,5.6){\includegraphics[width=4.0cm,height=3.3cm]{Fig_ct-2.jpg}}
	\put(3.2,4.5){\includegraphics[width=6.4cm,height=4.8cm]{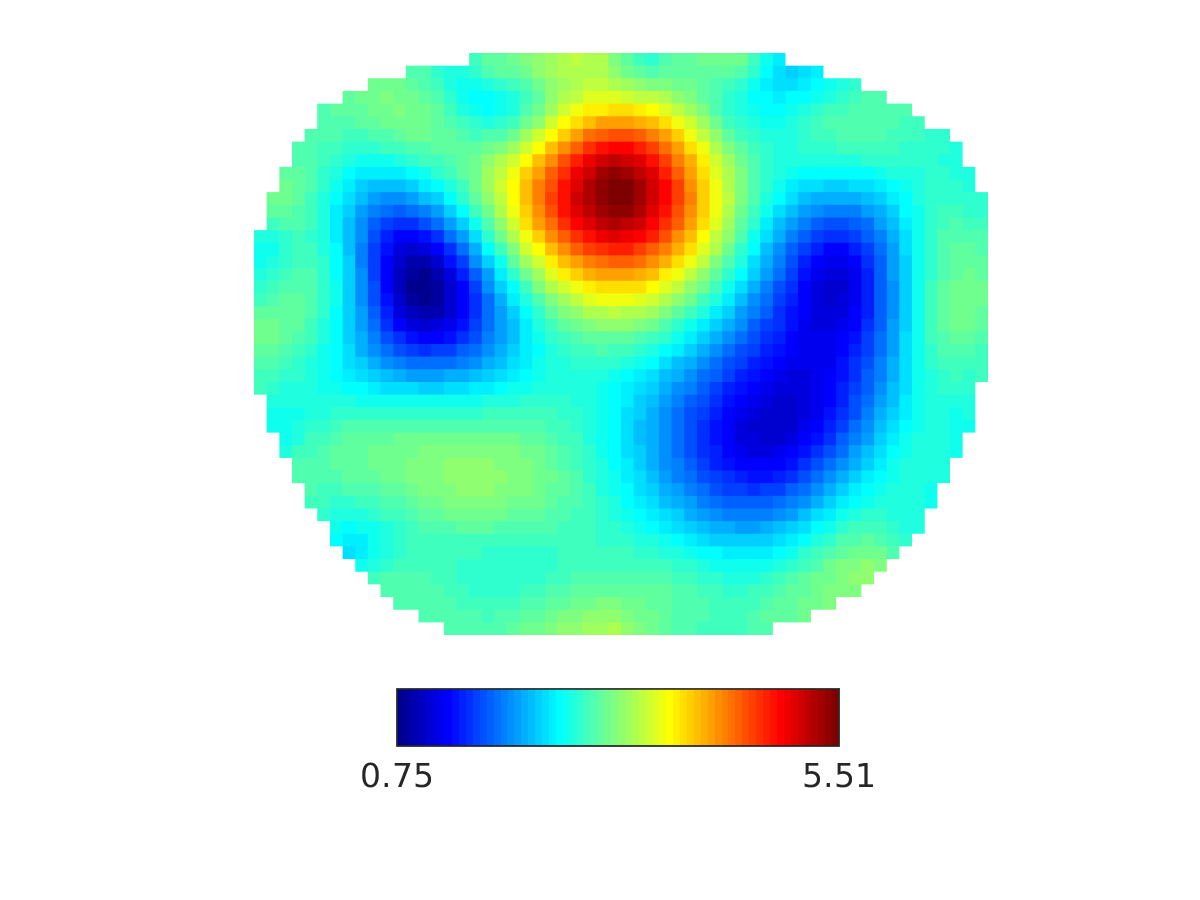}}
	\put(8.5,4.5){\includegraphics[width=6.4cm,height=4.8cm]{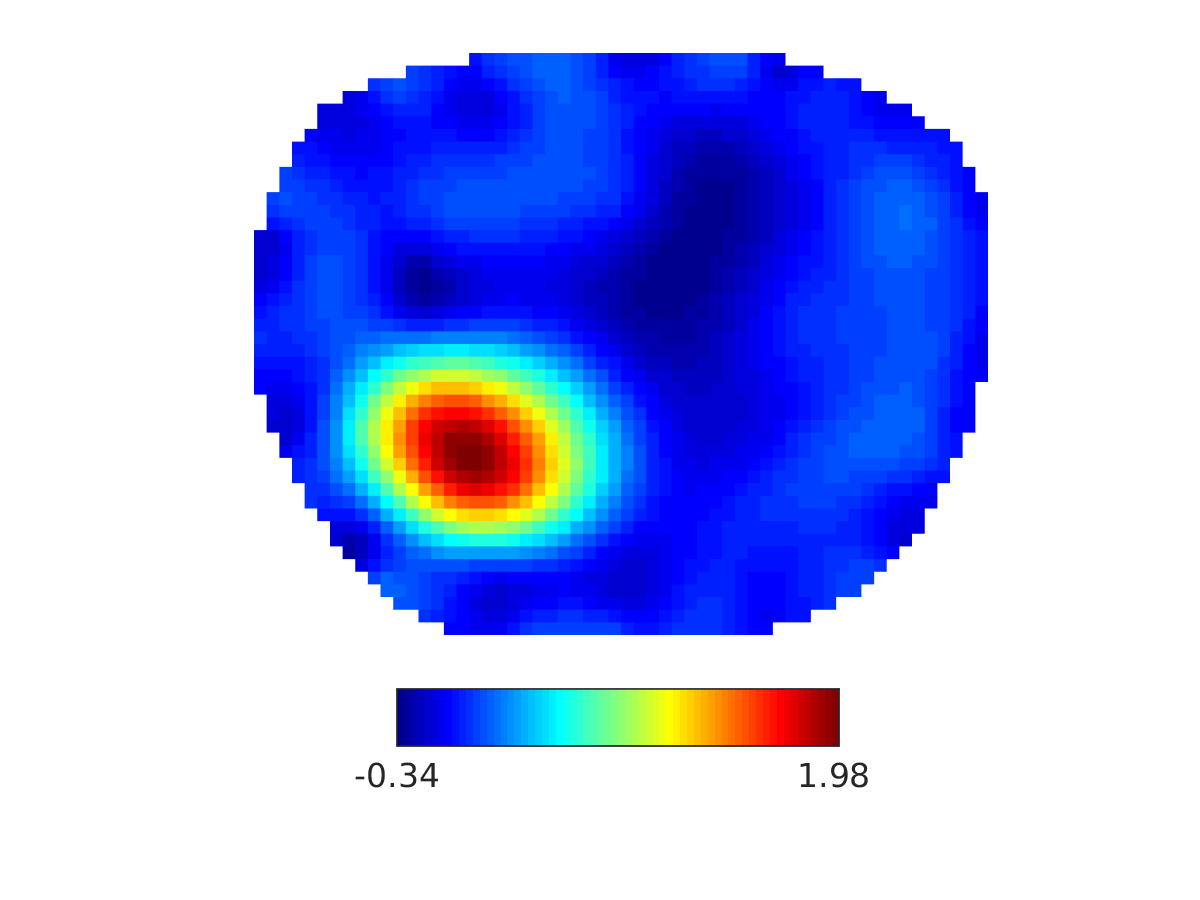}}
	\put(-1.0,0.6){\includegraphics[width=4.0cm,height=3.3cm]{Fig_ct-3.jpg}}
	\put(3.2,-0.5){\includegraphics[width=6.4cm,height=4.8cm]{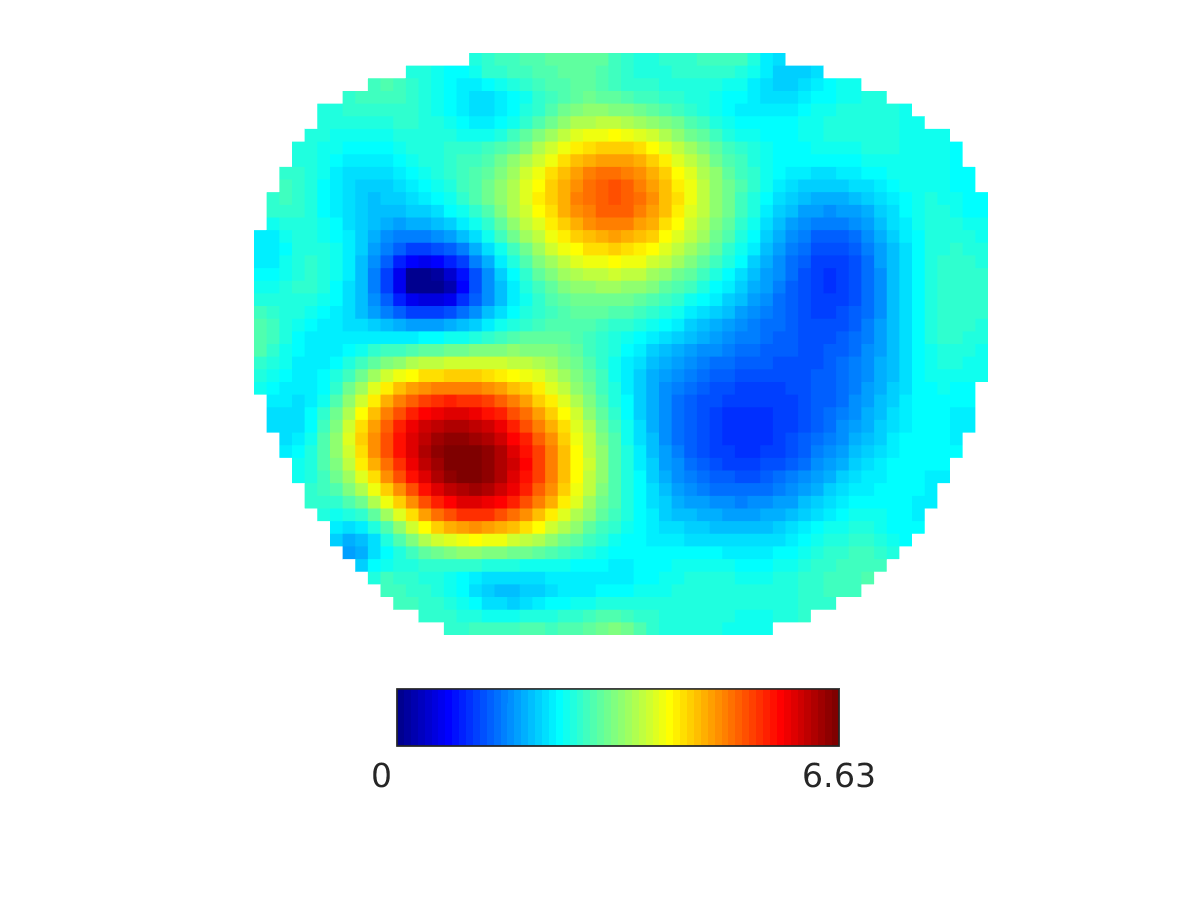}}
	\put(8.5,-0.5){\includegraphics[width=6.4cm,height=4.8cm]{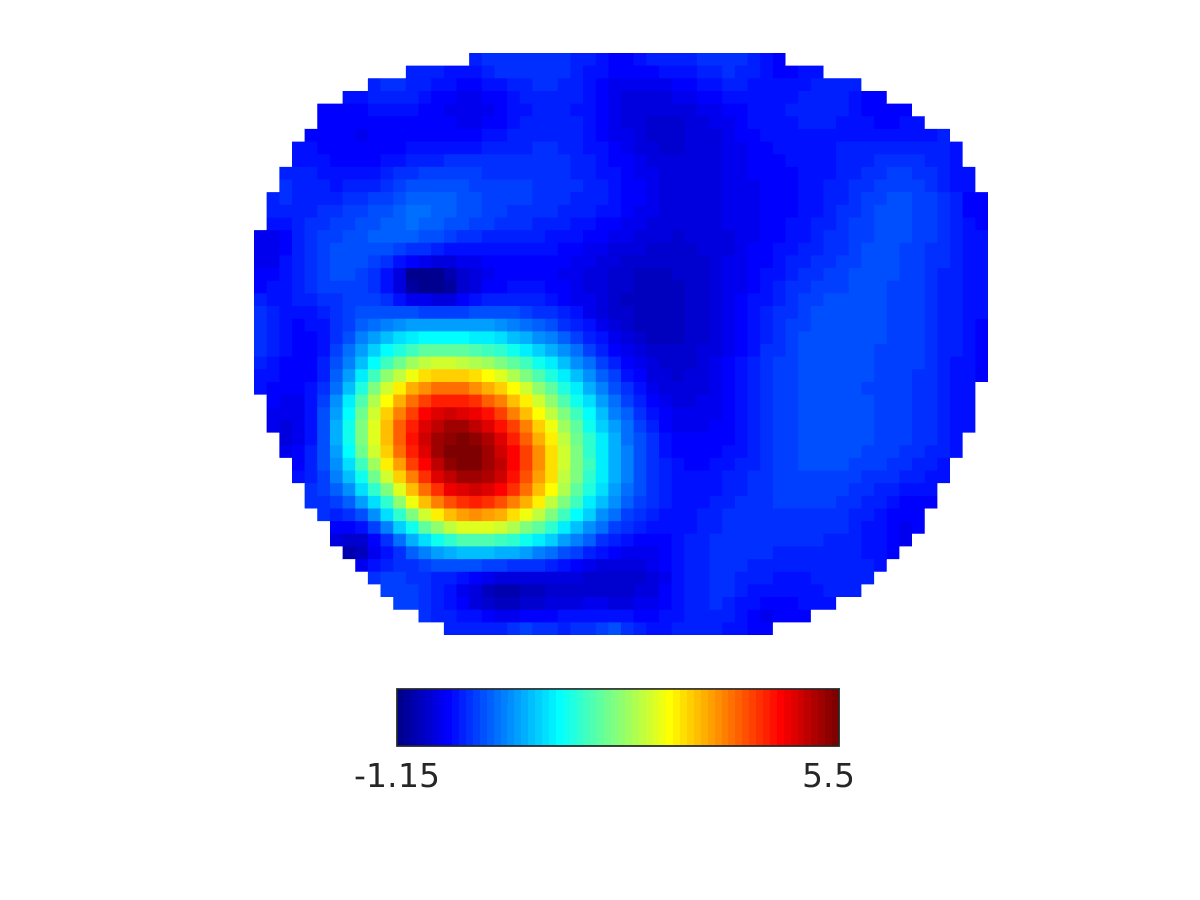}}
	\put(0.7,14.3){$\gamma_1$}
	\put(0.7,9.3){$\gamma_2$}
	\put(0.7,4.3){$\gamma_3$}
	\put(6.0,14.3){$\eta_1(x)$}
	\put(6.0,9.3){$\eta_2(x)$}
	\put(6.0,4.3){$\eta_3(x)$}
	\put(10.3,9.3){$\left( \gamma_2 -  \gamma_1 \right) (F_e^{-1}(x))$}
	\put(10.3,4.3){$\left( \gamma_3 -  \gamma_1 \right) (F_e^{-1}(x))$}
	\end{picture}    
	\caption{\textbf{Difference imaging electrodes displaced $25 \%$ (Case 1) $\gamma_1$ vs $\gamma_k, \, k=2,3$ in $\Omega_m = \Omega$.} Top row shows the results of experimental data 1,  images in middle row correspond to experimental data 2 and bottom row shows the results of experimental data 3.
		Left column: the 3 different experimental setups. 
		Middle column: functions $\eta_1 (x)$, $\eta_2 (x)$ and $\eta_3 (x)$, for  $x  \in \Omega_m = \Omega$.
		Right column: difference  $\left (\eta_2 - \eta_1 \right) (x) = \left( \gamma_2 -  \gamma_1 \right) (F_e^{-1}(x))$ and  $\left( \eta_3 - \eta_1 \right)(x) = \left( \gamma_3 -  \gamma_1 \right) (F_e^{-1}(x))$, respectively (see equation~\eqref{eq:diff-imag-Om}).}
	\label{fig:diff_imag_eldis_35_exp_1_2_3}
\end{figure}

\begin{figure}
	\centering
	\setlength{\unitlength}{1cm}
	\begin{picture}(20,14.5) 	
	\put(-1.0,10.6){\includegraphics[width=4.0cm,height=3.3cm]{Fig_ct-1.jpg}}
	\put(3.2,9.5){\includegraphics[width=6.4cm,height=4.8cm]{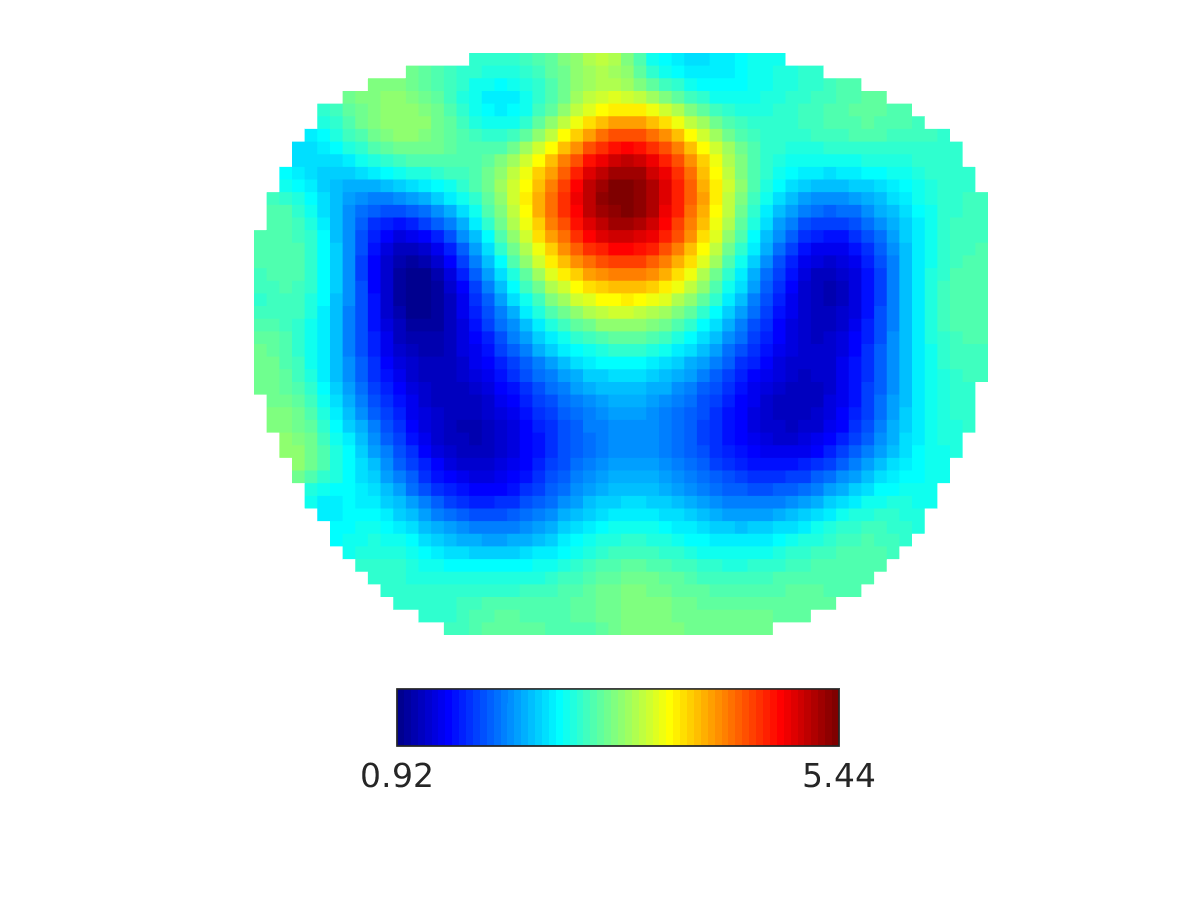}}
	\put(-1.0,5.6){\includegraphics[width=4.0cm,height=3.3cm]{Fig_ct-2.jpg}}
	\put(3.2,4.5){\includegraphics[width=6.4cm,height=4.8cm]{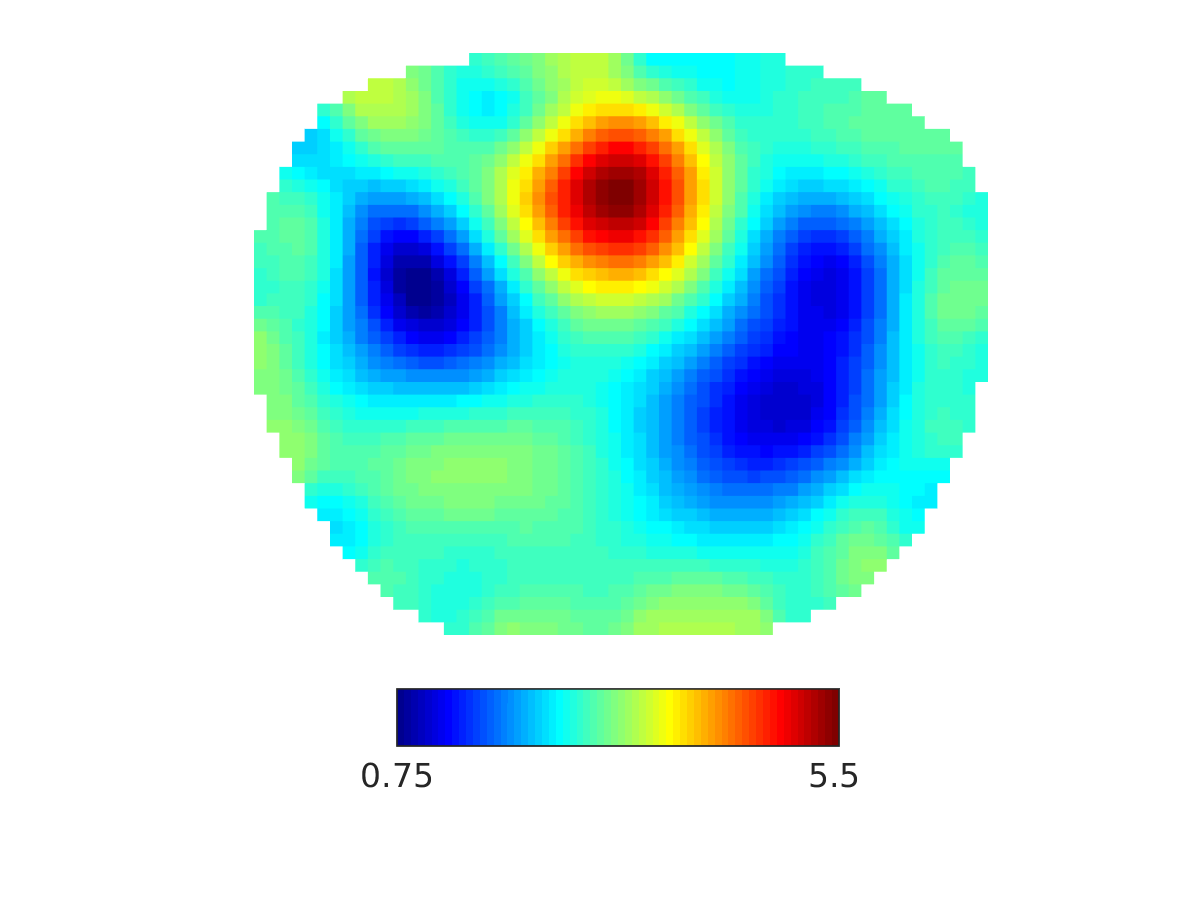}}
	\put(8.5,4.5){\includegraphics[width=6.4cm,height=4.8cm]{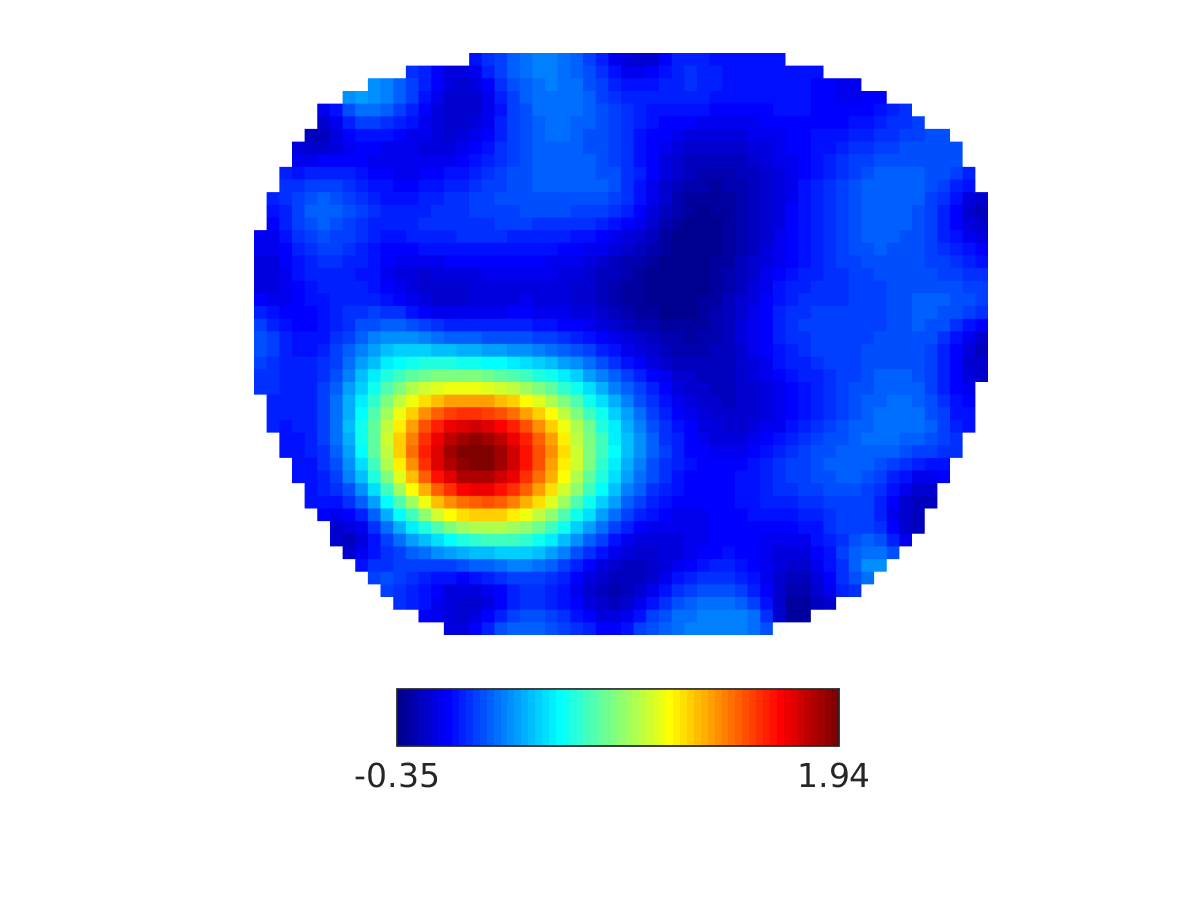}}
	\put(-1.0,0.6){\includegraphics[width=4.0cm,height=3.3cm]{Fig_ct-3.jpg}}
	\put(3.2,-0.5){\includegraphics[width=6.4cm,height=4.8cm]{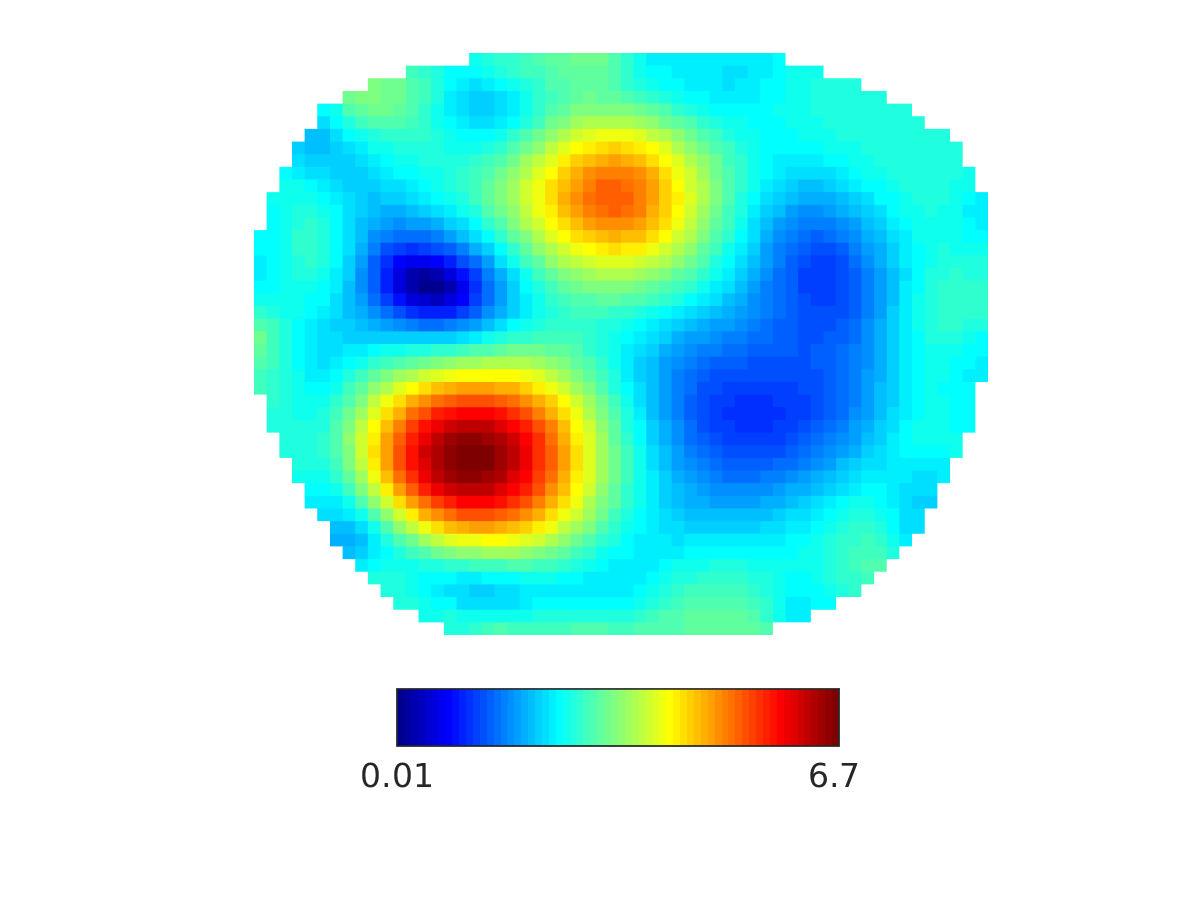}}
	\put(8.5,-0.5){\includegraphics[width=6.4cm,height=4.8cm]{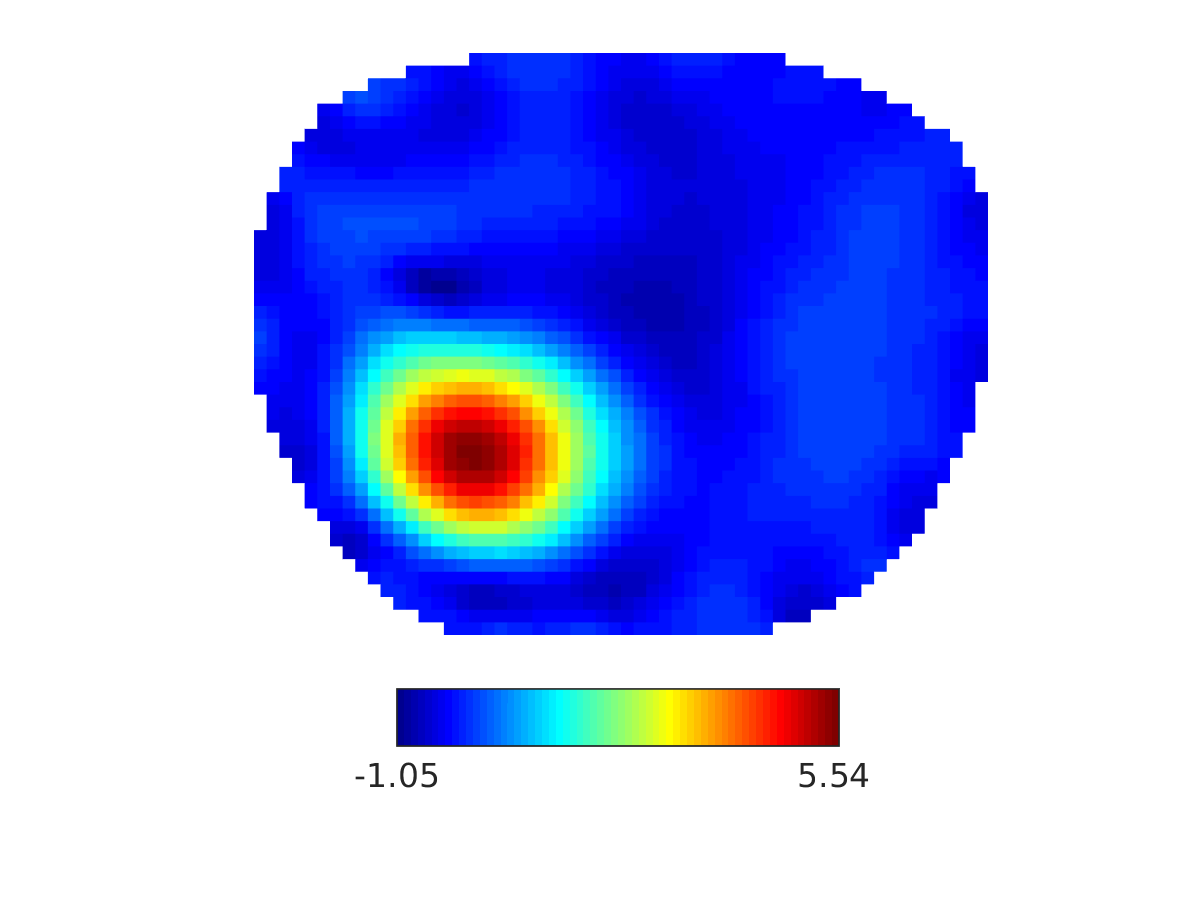}}
	\put(0.7,14.3){$\gamma_1$}
	\put(0.7,9.3){$\gamma_2$}
	\put(0.7,4.3){$\gamma_3$}
	\put(6.0,14.3){$\eta_1(x)$}
	\put(6.0,9.3){$\eta_2(x)$}
	\put(6.0,4.3){$\eta_3(x)$}
	\put(10.3,9.3){$\left( \gamma_2 -  \gamma_1 \right) (F_e^{-1}(x))$}
	\put(10.3,4.3){$\left( \gamma_3 -  \gamma_1 \right) (F_e^{-1}(x))$}
	\end{picture}    
	\caption{\textbf{Difference imaging electrodes displaced $35 \%$ (Case 2) $\gamma_1$ vs $\gamma_k, \, k=2,3$ in $\Omega_m = \Omega$.}
		Top row shows the results of experimental data 1,  images in middle row correspond to experimental data 2 and bottom row shows the results of experimental data 3.
		Left column: the 3 different experimental setups. 
		Middle column: functions $\eta_1 (x)$, $\eta_2 (x)$ and $\eta_3 (x)$, for  $x  \in \Omega_m = \Omega$.
		Right column: difference  $\left (\eta_2 - \eta_1 \right) (x) = \left( \gamma_2 -  \gamma_1 \right) (F_e^{-1}(x))$ and  $\left( \eta_3 - \eta_1 \right)(x) = \left( \gamma_3 -  \gamma_1 \right) (F_e^{-1}(x))$, respectively (see equation~\eqref{eq:diff-imag-Om}).}
	\label{fig:diff_imag_eldis_50_exp_1_2_3}
\end{figure}

Note that in the situation of electrodes displaced we have $\Omega_{m} = \Omega$, then we can directly compare the reconstructions in the true domain using equation~\eqref{eq:diff-imag-Om}.
In Figure~\ref{fig:diff_imag_eldis_35_exp_1_2_3} we present the results corresponding to the situation where electrodes were  displaced  approximately $25\%$ of the true physical length (Case 1) and  in Figure~\ref{fig:diff_imag_eldis_50_exp_1_2_3} the results corresponding to the situation where electrodes were  displaced  approximately  $35\%$  of  the  true physical length (Case 2).

\subsubsection{ Two conductivities $\gamma_1$ and $\gamma_2$ defined in two different domains}

\begin{figure}
	\centering
	\setlength{\unitlength}{1cm}
	\begin{picture}(20,14.5) 	
	\put(-1.0,10.8){\includegraphics[width=3.9cm,height=3.1cm]{Fig_ct-1.jpg}}
	\put(3.2,9.5){\includegraphics[width=6.4cm,height=4.8cm]{Fig_anisotropic_exp1}}
	\put(-1.0,5.8){\includegraphics[width=3.9cm,height=3.1cm]{Fig_ct-2.jpg}}
	\put(3.2,4.5){\includegraphics[width=6.4cm,height=4.8cm]{Fig_anisotropic_exp2}}
	\put(8.4,4.5){\includegraphics[width=6.4cm,height=4.8cm]{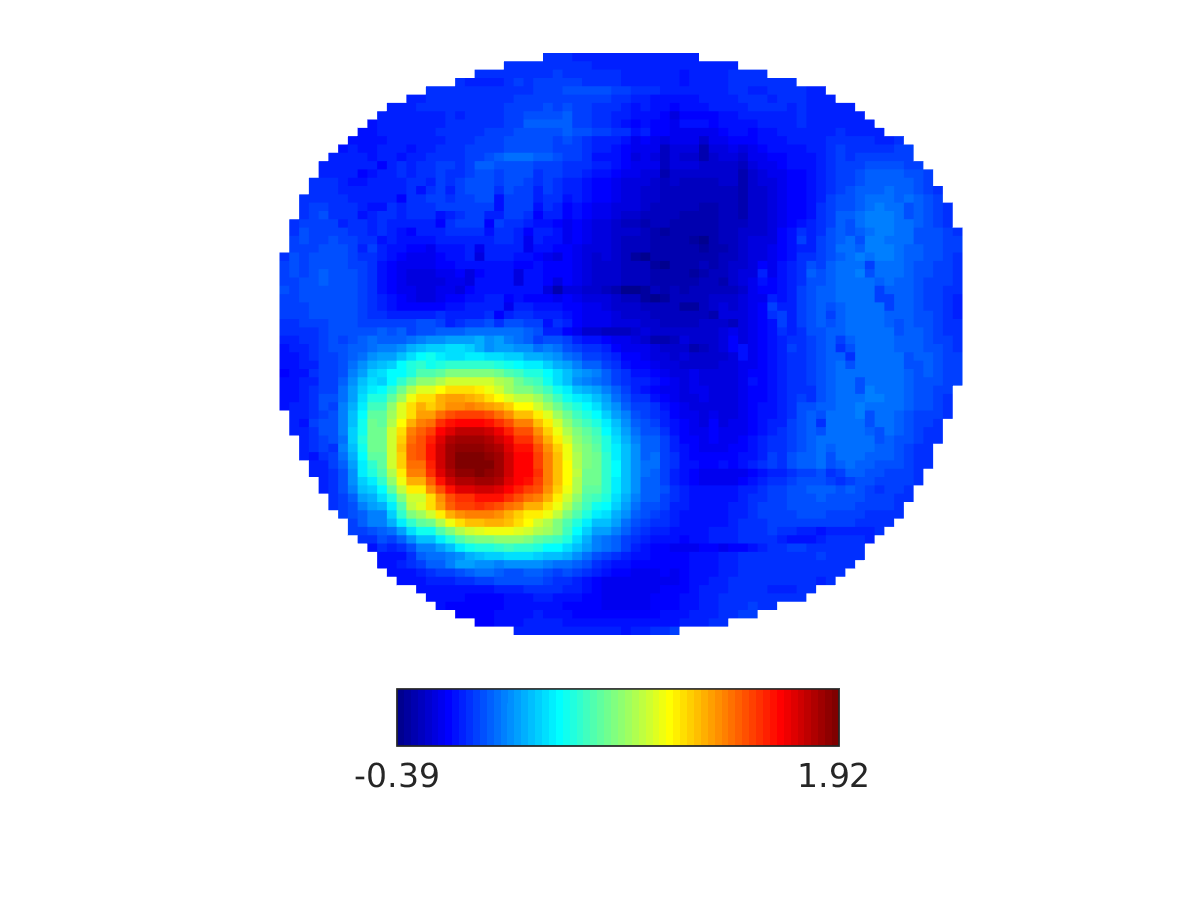}}
	\put(-1.0,0.8){\includegraphics[width=3.9cm,height=3.1cm]{Fig_ct-3.jpg}}
	\put(3.2,-0.5){\includegraphics[width=6.4cm,height=4.8cm]{Fig_anisotropic_exp3}}
	\put(8.5,-0.5){\includegraphics[width=6.4cm,height=4.8cm]{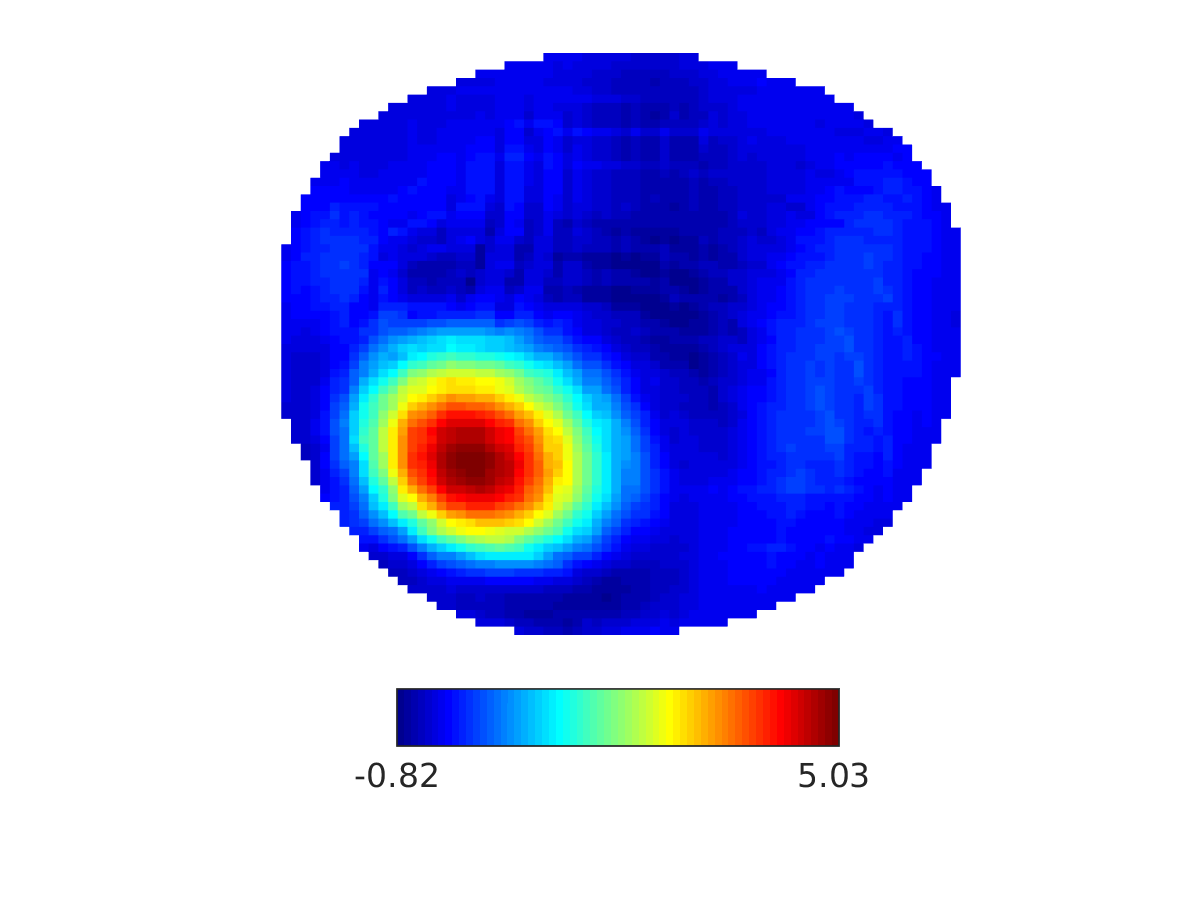}}
	\put(0.7,14.3){$\gamma_1$}
	\put(0.7,9.3){$\gamma_2$}
	\put(0.7,4.3){$\gamma_3$}
	\put(6.1,14.3){$\gamma_{c,1}$}
	\put(6.1,9.3){$\gamma_{c,2}$}
	\put(6.1,4.3){$\gamma_{c,3}$}
	\put(11.0,9.3){$ \gamma_{2,c} - \gamma_{1,c}$}
	\put(11.0,4.3){$ \gamma_{3,c} - \gamma_{1,c}$}
	\end{picture}    
	\caption{\textbf{Difference imaging: $\gamma_1$ vs $\gamma_k, \, k=2,3$ in $\Omega_{c,1} \cap \Omega_{c,k}$.} 
		Top row shows the results of experimental data 1,  images in middle row correspond to experimental data 2 and bottom row shows the results of experimental data 3.
		Left column: the 3 different experimental setups.
		Middle column: conductivities $\gamma_{c,k}(y)$, for $y  \in \Omega_{c,k}$.  
		Right column: difference  $\left( \gamma_{c,2} - \gamma_{c,1} \right) (y) = \left( \gamma_2 \circ \bar{G}_2^{-1} - \gamma_1 \circ \bar{G}_1^{-1} \right) (y) $ for $y \in \Omega_{c,2} \cap \Omega_{c,1}$ and $\left( \gamma_{c,3} - \gamma_{c,1} \right) (y) = \left( \gamma_3 \circ \bar{G}_3^{-1} - \gamma_1 \circ \bar{G}_1^{-1} \right) (y) $ for $y \in \Omega_{c,3} \cap \Omega_{c,1}$, respectively (see equation~\eqref{eq:diff-imag-2}).}
	\label{fig:diff_imag_exp_1_2_3_v2}
\end{figure}

In \eqref{eq:diff-imag} we assume that the unknown boundary does not change between the measurements. Nevertheless, that assumption is often infeasible. For instance, when imaging a human chest during a breathing cycle, the thorax shape and the contact impedances could vary between breathing states. If $\gamma_1$ and $\gamma_2$ are isotropic conductivities of the same object but corresponding to different measurement times and defined in slightly different domains, namely $\Omega_1$ and $\Omega_2$, then changes in the conductivity between the two measurement times could be determined by considering
\begin{equation}\label{eq:diff-imag-2}
\gamma_{c,1}(y) - \gamma_{c,2}(y) =  \gamma_1( \bar{G}_1^{-1}(y)) - \gamma_2 (\bar{G}_2^{-1}(y)), \qquad \text{for } \; y  \in \Omega_{c,1} \cap \Omega_{c,2}.
\end{equation} 
In Figure~\ref{fig:diff_imag_exp_1_2_3_v2} we show the results corresponding to the different experimental data sets presented in section~\ref{subsec:disc-model-domain}.

\subsection{Discussion on the results}\label{sec:discussion}

From the reconstructions in the top right of Figures~\ref{fig:exp1_Oc} to~\ref{fig:exp3_Oc}, we observe that the traditional approach
that ignores the modeling errors produces poor quality reconstructions and in consequence leads to the loss of useful information. On the contrary, utilizing the same incorrect model setting, the proposed approach produces good quality reconstructions carrying useful information about the target conductivity (see bottom left of Figures~\ref{fig:exp1_Oc} to~\ref{fig:exp3_Oc}). We can observe that both, the shapes of the domain and the inclusions are accurately reconstructed. Moreover, while the relative error~\eqref{eq:relative-error} between the true domain and the model domain was over 20\%, the error for the recovery domain is reduced to about $5\%$ in all cases.

In the situation when the true domain is known but the locations of the electrodes are inaccurately known the traditional approach gives reconstruction having severe artifacts, mainly close to the tank boundary (see middle column of Figures~\ref{fig:case1_eldis_35} and \ref{fig:case2_eldis_50}).
On the contrary, the proposed approach eliminates these artifacts and also reproduces the shapes of the inclusions more accurately (see last column of Figures~\ref{fig:case1_eldis_35} and~\ref{fig:case2_eldis_50}). These results indicate that the proposed method tolerates well inaccurate knowledge of the electrode locations.

Finally, as can be seen from Figures~\ref{fig:diff_imag_exp_1_2_3_v1} to~\ref{fig:diff_imag_exp_1_2_3_v2}, the proposed approach seems to be a  feasible tool for dynamic imaging both in the case when the shape boundary is known but electrodes positions are inaccurately known (see Figures~\ref{fig:diff_imag_eldis_35_exp_1_2_3} and \ref{fig:diff_imag_eldis_50_exp_1_2_3}) and in the case when the shape boundary and electrode positions are unknown (see Figures~\ref{fig:diff_imag_exp_1_2_3_v1} and \ref{fig:diff_imag_exp_1_2_3_v2}).

\section{Conclusion}\label{sec:conclusion}
A typical difficulty in practical EIT is that in most measurement situations the knowledge about the boundary of the body, the electrode locations and the contact impedances are usually uncertain.
It is widely known that these modeling errors can cause severe artifacts that ruin the quality of the image reconstruction and in consequence diagnostically relevant information is lost. In this paper, we have introduced a method that is capable of producing good quality reconstructions of the conductivity in settings where the boundary of the body, the electrode locations and the contact impedances are inaccurately known.
The method was evaluated via experimental studies with water tank data. The obtained results indicate that the proposed approach is an effective tool to overcome the difficulties caused by the uncertainties in the measurement configuration usually present in practical EIT.

\section*{Acknowledgements}
The work of JPA was supported by the National Scientific and Technical Research Council of Argentina Grant PIP 11220150100500CO and by Secyt (UNC) Grant 33620180100326CB.
The work of of VK was supported by the Academy of Finland (Projects 312343 and 336791, Finnish Centre of Excellence in Inverse Modelling and Imaging) and the Jane and Aatos Erkko Foundation.
The work of ML, PO and SS was supported by the  Academy of Finland (Finnish Centre of Excellence in Inverse Modelling and Imaging and projects 273979, 284715, and 312110) and the Jane and Aatos Erkko Foundation.

\bibliographystyle{plain}
\bibliography{boundary-project.bib}

\begin{thebibliography}{10}

\bibitem{AGB96}
A.~{Adler}, R.~{Guardo}, and Y.~{Berthiaume}.
\newblock Impedance imaging of lung ventilation: Do we need to account for
  chest expansion?
\newblock {\em IEEE Trans. Biomed. Eng.}, 43:414--420, 1996.

\bibitem{ACLMSS20}
J.P. {Agnelli}, A.~{{\c{C}}öl}, M.~{Lassas}, R.~{Murthy}, M.~{Santacesaria},
  and S.~{Siltanen}.
\newblock Classification of stroke using neural networks in electrical
  impedance tomography.
\newblock {\em Inverse Problems}, 36(11):115008, 2020.

\bibitem{A66}
L.~Ahlfors.
\newblock {\em Lectures on Quasiconformal Mappings}.
\newblock Van Nostrand Mathematical Studies 10. D. Van Nostrand, 1966.

\bibitem{APKSMH16}
K.Y. {Aristovich}, B.C. {Packham}, H.~{Koo}, G.S.D. {Santos}, A.~{McEvoy}, and
  D.S. {Holder}.
\newblock Imaging fast electrical activity in the brain with electrical
  impedance tomography.
\newblock {\em NeuroImage}, 124:204--213, 2016.

\bibitem{AP06}
K.~{Astala} and L.~{Päivärinta}.
\newblock {C}alderón’s inverse conductivity problem.
\newblock {\em Ann. of Math.}, 163:265--299, 2006.

\bibitem{BB84}
D.C. {Barber} and B.H. {Brown}.
\newblock Applied potential tomography.
\newblock {\em J. Phys. E: Sci. Instrum.}, 17:723--733, 1984.

\bibitem{BKKK08}
G.~{Boverman}, T.J. {Kao}, R.~{Kulkarni}, B.S. {Kim}, D.~{Isaacson}, G.J.
  {Saulnier}, and J.C. {Newell}.
\newblock Robust linearized image reconstruction for multifrequency eit of the
  breast.
\newblock {\em IEEE Trans. Med. Imag.}, 27:1439--1448, 2008.

\bibitem{BLMM98}
V.~{Bozin}, N.~{Lakic}, V.~{Markovic}, and M.~{Mateljevic}.
\newblock Unique extremality.
\newblock {\em J. Anal. Math}, 75:299--338, 1998.

\bibitem{CHH19}
V.~{Candiani}, A.~{Hannukainen}, and N.~{Hyvönen}.
\newblock Computational framework for applying electrical impedance tomography
  to head imaging.
\newblock {\em SIAM J. Sci. Comput.}, 41(5):B1034--B1060, 2019.

\bibitem{DHSS13A}
J.~{Dardé}, N.~{Hyvönen}, A.~{Seppänen}, and S.~{Staboulis}.
\newblock Simultaneous reconstruction of outer boundary shape and admittivity
  distribution in electrical impedance tomography.
\newblock {\em SIAM J. Imaging Sci.}, 1(6):176--198, 2013.

\bibitem{DHSS13b}
J.~{Dardé}, N.~{Hyvönen}, A.~{Seppänen}, and S.~{Staboulis}.
\newblock Simultaneous recovery of admittivity and body shape in electrical
  impedance tomography: {A}n experimental evaluation.
\newblock {\em Inverse Problems}, 29(8):085004, 2013.

\bibitem{GHO96}
E.~{Gersing}, B.~{Hoffman}, and M.~{Osypka}.
\newblock Influence of changing peripheral geometry on electrical impedance
  tomography measurements.
\newblock {\em Medical and Biological Engineering and Computing}, 34:359--361,
  1996.

\bibitem{HVWKV02}
L.M. {Heikkinen}, T.~{Vilhunen}, R.M. {West}, J.P. {Kaipio}, and
  M.~{Vauhkonen}.
\newblock Simultaneous reconstruction of electrode contact impedances and
  internal electrical properties: 2. laboratory experiments.
\newblock {\em Meas. Sci. Tech.}, 13:1855--1861, 2002.

\bibitem{HWWT93}
P.~{Hua}, E.J. {Woo}, J.G. {Webster}, and W.J. {Tompkins}.
\newblock Using compund electrodes in electrical impedance tomography.
\newblock {\em IEEE Trans. Biomed. Eng.}, 40(1):29--34, 1993.

\bibitem{HKMS17}
N.~{Hyvönen}, V.~{Kaarnioja}, L.~{Mustonen}, and S.~{Staboulis}.
\newblock Polynomial collocation for handling an inaccurately known measurement
  configuration in electrical impedance tomography.
\newblock {\em SIAM J. Appl. Math.}, 77(1):202--202, 2017.

\bibitem{HMS17}
N.~{Hyvönen}, H.~{Majander}, and S.~{Staboulis}.
\newblock Compensation for geometric modeling errors by positioning of
  electrodes in electrical impedance tomography.
\newblock {\em Inverse Problems}, 33:035006, 2017.

\bibitem{HM17}
N.~{Hyvönen} and L.~{Mustonen}.
\newblock Smoothened complete electrode model.
\newblock {\em SIAM J. Appl. Math.}, 77:2250--2271, 2017.

\bibitem{IMNS04}
D.~{Isaacson}, J.~{Mueller}, J.C.{Newell}, and S.~{Siltanen}.
\newblock Reconstructions of chest phantoms by the d-bar method for electrical
  impedance tomography.
\newblock {\em IEEE Trans. Med. Im.}, 23:821--S28, 2004.

\bibitem{IMNS06}
D.~{Isaacson}, J.~{Mueller}, J.C.{Newell}, and S.~{Siltanen}.
\newblock Imaging cardiac activity by the d-bar method for electrical impedance
  tomography.
\newblock {\em Physiological Measurement}, 27:843--S50, 2006.

\bibitem{KKSV00}
J.P. {Kaipio}, V.~{Kolehmainen}, E.~{Somersalo}, and M.~{Vauhkonen}.
\newblock Statistical inversion and monte carlo sampling methods in electrical
  impedance tomography.
\newblock {\em Inverse Problems}, 16:1487--1522, 2000.

\bibitem{KV87}
R.~Kohn and M.~Vogelius.
\newblock Relaxation of a variational method for impedance computed tomography.
\newblock {\em Communications on Pure and Applied Mathematics}, 40(6):745--777,
  1987.

\bibitem{KLO05}
V.~{Kolehmainen}, M.~{Lassas}, and P.~{Ola}.
\newblock The inverse conductivity problem with an imperfectly known boundary.
\newblock {\em SIAM J. App. Math.}, 66(2):365--383, 2005.

\bibitem{KLO08}
V.~{Kolehmainen}, M.~{Lassas}, and P.~{Ola}.
\newblock Electrical impedance tomography problem with inaccurately known
  boundary and contact impedances.
\newblock {\em IEEE Transactions on Medical Imaging}, 27(10):1404--1414, 2008.

\bibitem{KLO10}
V.~{Kolehmainen}, M.~{Lassas}, and P.~{Ola}.
\newblock Calderon's inverse problem with an imperfectly known boundary and
  reconstruction up to a conformal deformation.
\newblock {\em SIAM J. Math. Anal.}, 42(3):1371--1381, 2010.

\bibitem{KLOS13}
V.~{Kolehmainen}, M.~{Lassas}, P.~{Ola}, and S.~{Siltanen}.
\newblock Recovering boundary shape and conductivity in electrical impedance
  tomography.
\newblock {\em Inverse Problems and Imaging}, 7(1):217--242, 2013.

\bibitem{KVKK97}
V.~{Kolehmainen}, M.~{Vauhkonen}, P.A. {Karjalainen}, and J.P. {Kaipio}.
\newblock Assessment of errors in static electrical impedance tomography with
  adjacent and trigonometric current patterns.
\newblock {\em Physiological Measurement}, 18:289--303, 1997.

\bibitem{KIT4}
J.~{Kourunen}, {Savolainen} J., A.~{Lehikoinen}, M.~{Vauhkonen}, and
  L.~{Heikkinen}.
\newblock Suitability of a pxi platform for an electrical impedance tomography
  system.
\newblock {\em Meas. Sci. Technol.}, 20:015503, 2008.

\bibitem{LKSS15}
D.~{Liu}, V.~{Kolehmainen}, S.~{Siltanen}, and A.~{Seppänen}.
\newblock A nonlinear approach to difference imaging in eit; assessment of the
  robustness in the presence of modelling errors.
\newblock {\em Inverse Problems}, 31:035012, 2015.

\bibitem{N96}
A.~{Nachman}.
\newblock Global uniqueness for a two-dimensional inverse boundary value
  problem.
\newblock {\em Ann. of Math.}, 143:71--96, 1996.

\bibitem{NKKa}
A.~{Nissinen}, V.~{Kolehmainen}, and J.P. {Kaipio}.
\newblock Compensation of modelling errors due to unknown domain boundary in
  electrical impedance tomography.
\newblock {\em IEEE Trans. Med. Imag.}, 30:231–242, 2011.

\bibitem{NKKb}
A.~{Nissinen}, V.~{Kolehmainen}, and J.P. {Kaipio}.
\newblock Reconstruction of domain boundary and conductivity in electrical
  impedance tomography using the approximation error approach.
\newblock {\em International Journal for Uncertainty Quantification},
  1:203–222, 2011.

\bibitem{NW06}
J.~Nocedal and S.J. Wright.
\newblock {\em Numerical Optimization}.
\newblock Springer Series in Operations Research and Financial Engineering.
  Springer, 2006.

\bibitem{R97}
M.~{Raydan}.
\newblock The {B}arzilai and {B}orwein gradient method for the large scale
  unconstrained minimization problem.
\newblock {\em SIAM J. Optim.}, 7(1):26--23, 1997.

\bibitem{SGA06}
M.~{Soleimani}, C.~{Gómez-Laberge}, and A.~{Adler}.
\newblock Imaging of conductivity changes and electrode movement in eit.
\newblock {\em Physiol. Meas.}, 27:S103--S113, 2006.

\bibitem{SChI92}
E.~{Somersalo}, M.~{Cheney}, and D.~{Isaacson}.
\newblock Existence and uniqueness for electrode models for electric current
  computed tomography.
\newblock {\em SIAM J. Appl. Math.}, 52(4):1023--1040, 1992.

\bibitem{S76}
K.~{Strebel}.
\newblock On the existence of extremal teichmüller mappings.
\newblock {\em J. Anal. Math}, 30:464--480, 1976.

\bibitem{sylvester1990anisotropic}
John Sylvester.
\newblock An anisotropic inverse boundary value problem.
\newblock {\em Communications on Pure and Applied Mathematics}, 43(2):201--232,
  1990.

\bibitem{GV97}
G.~{Vainikko}.
\newblock Fast solvers of the {L}ippmann-{S}chwinger equation.
\newblock {\em in “Direct and Inverse Problems of Mathematical Physics”
  (Newark, DE, 1997) Int. Soc. Anal. Appl. Comput., Kluwer Acad. Publ,
  Dordrecht}, 5:423--440, 2000.

\bibitem{W05}
M.H. {Wright}.
\newblock The interior-point revolution in optimization: History, recent
  developments, and lasting consequences.
\newblock {\em Bull. Amer. Math. Soc.}, 42:39--56, 2005.

\end{thebibliography}

\end{document}